\documentclass[10pt,reqno]{amsart}

\usepackage{amsmath}
\usepackage{amsthm}
\usepackage{amssymb}
\usepackage{amsfonts}
\usepackage{amsxtra}
\usepackage{multicol}
\usepackage{fullpage}
\usepackage{amscd}
\usepackage{color} 
\usepackage{bm} 
\usepackage{mathrsfs}
\usepackage{tikz}
\usepackage{tikz-cd}
\usetikzlibrary{cd}
\usepackage{enumerate}
\usepackage[all]{xy}
\usepackage[mathscr]{eucal}
\usepackage{verbatim}

\let\nc\newcommand
\let\renc\renewcommand

\theoremstyle{plain}
\newtheorem{thm}{Theorem}
\newtheorem{prop}[thm]{Proposition}
\newtheorem{cor}[thm]{Corollary}
\newtheorem{lem}[thm]{Lemma}

\theoremstyle{definition}
\newtheorem{defn}[thm]{Definition}
\newtheorem{example}[thm]{Example}

\theoremstyle{remark}
\newtheorem{remark}[thm]{Remark}
\newtheorem{rem}[thm]{Remark}
\numberwithin{thm}{section}

\nc{\bdm}{\begin{displaymath}}
\nc{\edm}{\end{displaymath}}
\nc{\bthm}{\begin{thm}}
\nc{\ethm}{\end{thm}}
\nc{\blem}{\begin{lem}}
\nc{\elem}{\end{lem}}
\nc{\bcor}{\begin{cor}}
\nc{\ecor}{\end{cor}}
\nc{\bprop}{\begin{prop}}
\nc{\eprop}{\end{prop}}
\nc{\bdef}{\begin{defn}}
\nc{\eddef}{\end{defn}}

\makeatletter
\renewcommand{\subsection}{\@startsection{subsection}{2}{0pt}{-3ex
plus -1ex minus -0.2ex}{-2mm plus -0pt minus
-2pt}{\normalfont\bfseries}} \makeatother

\numberwithin{equation}{section}

\newcommand{\Gr}{\mathrm{Gr}}

\DeclareMathOperator{\Ext}{\mathrm{Ext}}

\DeclareMathOperator{\Ker}{\mathrm{Ker}}

\DeclareMathOperator{\End}{\mathrm{End}}

\DeclareMathOperator{\gr}{\mathrm{gr}}

\DeclareMathOperator{\CC}{\mathcal{C}}

\DeclareMathOperator{\Tr}{\mathrm{Tr}}

\DeclareMathOperator{\Rep}{\mathrm{Rep}}

\newcommand{\bv}{\boldsymbol{v}}

\newcommand{\beq}{\begin{equation}\label}
\newcommand{\eeq}{\end{equation}}

\DeclareMathOperator{\Spec}{\mathrm{Spec}}
\DeclareMathOperator{\Id}{\mathrm{Id}}

\newcommand{\iso}{{\;\stackrel{_\sim}{\to}\;}}

\DeclareMathOperator{\Hom}{\mathrm{Hom}}
\DeclareMathOperator{\GL}{\mathrm{GL}}
\DeclareMathOperator{\Sp}{\mathrm{Sp}}

\nc{\Z}{\mathbb{Z}}
\newcommand{\N}{\mathbb{N}}
\newcommand{\Q}{\mathbb{Q}}
\newcommand{\R}{\mathbb{R}}
\newcommand{\C}{\mathbb{C}}
\newcommand{\h}{\mathfrak{h}}

\nc{\rank}{\textrm{rank} \,}
\nc{\ds}{\dots}

\let\mc\mathcal
\let\mf\mathfrak

\nc{\HW}{\bar{H}_{\mathbf{c}}(W)}
\nc{\HK}{\bar{H}_{\mathbf{c}}(K)}
\nc{\HtK}{\widetilde{H}_{\mathbf{c}}(K)}
\nc{\CMW}{\textsf{CM}_{\mbf{c}}(W)}
\nc{\CMK}{\textsf{CM}_{\mbf{c}}(K)}
\nc{\mbf}{\mathbf}
\nc{\LK}{\textsf{Irr}(K)}
\nc{\LW}{\textsf{Irr}(W)}
\nc{\Res}{\mathsf{Res} \, }
\nc{\Ind}{\mathsf{Ind} \, }

\nc{\cont}{\textrm{cont}}

\nc{\eWb}{\mathbf{e}_{W_b}}
\nc{\eW}{\mathbf{e}_{W}}
\nc{\msf}{\mathsf}
\nc{\Ui}{\mc{U}_{i,+}}
\nc{\Uone}{\mc{U}_{1,+}}
\nc{\Utwo}{\mc{U}_{2,+}}

\nc{\minusone}{-1}
\nc{\minustwo}{-2}
\nc{\Mod}{\mathrm{Mod} \,}
\nc{\ms}{\mathscr}
\nc{\Frac}{\mathrm{Frac} \,}
\nc{\ra}{\rightarrow}
\nc{\hra}{\hookrightarrow}
\nc{\lab}{\label}
\renc{\O}{\mc{O}}
\nc{\Tan}{\mc{T}}
\nc{\ul}{\underline}
\nc{\s}{\mathfrak{S}}
\nc{\g}{\mf{g}}
\nc{\pa}{\partial}
\nc{\tit}{\textit}
\nc{\Maxspec}{\mathrm{Maxspec} \, }
\nc{\gldim}{\mathrm{gl.dim}}
\nc{\rkm}{\mathrm{rk} \, (\mf{m})}
\nc{\sm}{\mathrm{sm}}
\nc{\PD}{\mathbb{PD}}
\nc{\hilb}{\textrm{Hilb}}
\nc{\mr}{\mathrm}
\renc{\>}{\rangle}
\nc{\T}{\mathbb{T}}
\nc{\X}{\mathbb{X}}
\nc{\F}{\mathbb{F}}
\nc{\kt}{\mbf{k}}
\nc{\ko}{\mbf{k}(0)}
\nc{\Ok}{\mc{O}_G \boxtimes \kt_X}
\nc{\Oko}{\mc{O}_G \boxtimes \ko_X}
\nc{\OYk}{\mc{O}_Y \boxtimes \kt_X}
\nc{\id}{\msf{id}}
\nc{\A}{\mathbb{A}}
\nc{\Grat}{\mc{Grat}}
\nc{\Squo}[1]{\A^{(#1)}}
\nc{\twist}{\mathrm{twist}}
\nc{\Cd}{\mc{C}}
\nc{\Span}{\mathrm{Span}}
\nc{\Grass}{\mathrm{Gr}}
\nc{\Fr}{\mathrm{Fr}}
\nc{\pco}[1]{k[V]^{p\mathrm{co} #1}}
\nc{\Irr}{\mathrm{Irr}}
\renc{\o}{\otimes}
\renc{\gr}{\mathsf{gr}}
\nc{\fin}{\mathrm{fin}}
\nc{\bk}{\mathbf{k}}
\nc{\bn}{\mathbf{n}}
\nc{\algD}{\mf{D}}
\nc{\hr}{\mf{h}_{\textrm{reg}}}
\nc{\D}{\mathscr{D}}
\nc{\PIdeg}{\mathrm{P.I.-degree}}
\nc{\ch}{\mathrm{ch}}
\nc{\ev}{\mathsf{ev}}
\nc{\Stab}{\mathrm{Stab}}
\nc{\Der}{\mathrm{Der}}
\nc{\rightsim}{\stackrel{\sim}{\longrightarrow}}
\nc{\HZ}{H_{\mbf{h},\Z}(\Z_m)}
\nc{\sing}{\mathrm{sing}}
\nc{\dd}{\mathscr{D}}
\nc{\bc}{\mathbf{c}}
\nc{\vc}{\underline{\mathbf{c}}}
\nc{\ba}{\mathbf{a}}
\nc{\reg}{\mathrm{reg}}
\nc{\Amp}{\mathrm{Amp}}
\nc{\Nef}{\mathrm{Nef}}
\nc{\SL}{\mathrm{SL}}
\nc{\Sym}{\mathrm{Sym}}
\nc{\Mov}{\mathrm{Mov}}
\nc{\Pic}{\mathrm{Pic}}
\nc{\Cs}{\C^{\times}}
\nc{\Nak}[3]{\mf{M}_{{#1}} ({#2},{#3}) }
\nc{\Naka}[2]{\mf{M}({#1},{#2}) }
\nc{\Mtheta}[1]{\mc{M}_{#1}}

\DeclareMathOperator{\aff}{aff}

\newcommand{\B}{\mathsf{B}}
\nc{\Qu}{\mathsf{Q}}
\nc{\bw}{\mathbf{w}}
\nc{\Nq}{\underline{N}}
\nc{\qu}{\varrho}
\nc{\QuGinfty}{\Qu(\Gamma,\infty)}
\renc{\SS}[1]{\Sigma \Sigma_{#1}}

\newcommand{\git}{\ensuremath{/\!\!/\!}}

\nc{\Qulinf}{\Qu^\ell(\infty)}
\nc{\oQulinf}{\overline{\Qu^\ell(\infty)}}
\nc{\oQuii}{\overline{\Qu^\infty(\infty)}}
\nc{\smod}{\,\,\mathrm{mod}\,\,}

\usepackage[vcentermath]{youngtab} 
\usepackage{todonotes}

\DeclareMathOperator{\Core}{{\mathrm{Core}}}
\DeclareMathOperator{\diag}{{\mathrm{diag}}}
\renewcommand{\Ref}{\mathrm{Ref}}

\newcommand{\longiso}{\stackrel{\sim}{\longrightarrow}}

\def\a{\alpha}

\def\e{\varepsilon}

\def\CM{{\mathbb{C}}}
\def\ZM{{\mathbb{Z}}}

\def\PC{{\mathcal{P}}}
\def\XC{{\mathcal{X}}}
\def\YC{{\mathcal{Y}}}

\def\Hb{{\mathbf H}}
\def\Zb{{\mathbf Z}} 
\def\Pb{{\mathbf P}}

\def\aff{\mathrm{aff}}

\def\Qov{{\overline{\Qu}}}

\def\Irm{{\mathrm{I}}}
\def\Trm{{\mathrm{T}}}

\def\vide{\varnothing}
\def\DS{\displaystyle}
\def\ZCB{{\boldsymbol{\mathcal{Z}}}}
\def\longto{\longrightarrow}
 
\def\ZCB{{\boldsymbol{\mathcal{Z}}}}
\def\rdim{\operatorname{\dim^{\rm reg}}}
\def\LG{{\mathfrak L}}

\def\lexp#1#2{\kern\scriptspace\vphantom{#2}^{#1}\kern-\scriptspace#2}

\definecolor{shadecolor}{gray}{0.90}

\def\equat{\refstepcounter{thm}\begin{equation}}
\def\endequat{\end{equation}}

\def\fonction#1#2#3#4#5{\begin{array}{rccc}
{#1} : & {#2} & \longto & {#3}  \\
& {#4} & \longmapsto & {#5} 
\end{array}}

\begin{document}

\title{Singularities in Calogero--Moser varieties}

\author{Gwyn Bellamy}
\address{School of Mathematics and Statistics, University of Glasgow, 15 University Gardens, Glasgow, G12 8QW, UK.}
\email{gwyn.bellamy@glasgow.ac.uk}
\urladdr{http://www.maths.gla.ac.uk/~gbellamy/}

\author{Ruslan Maksimau}
\address{Laboratoire Analyse, Géométrie et Modélisation, CY Cergy Paris Université, 2 av. Adolphe Chauvin, 95302 Cergy-Pontoise, France}
\email{ruslan.maksimau@cyu.fr, ruslmax@gmail.com}
\urladdr{https://maksimau.perso.math.cnrs.fr/}

\author{Travis Schedler} 
\address{Mathematics Department, Imperial College London, London, SW7 2AZ, UK.}
\email{t.schedler@imperial.ac.uk}
\urladdr{https://www.imperial.ac.uk/people/t.schedler/}

\begin{abstract}
In this article we describe completely the singularities appearing in Calogero--Moser varieties associated (at any parameter) to the wreath product symplectic reflection groups. We do so by parameterizing the symplectic leaves in the variety, describing combinatorially the resulting closure relation and computing a transverse slice to each leaf. We also show that the normalization of the closure of each symplectic leaf is isomorphic to a Calogero--Moser variety for an associated (explicit) subquotient of the symplectic reflection group. This confirms a conjecture of Bonnaf\'e for these groups. 

We use the fact that the Calogero--Moser varieties associated to wreath products can be identified with certain Nakajima quiver varieties. In particular, our result identifying the normalization of the closure of each symplectic leaf with another quiver variety holds for arbitrary quiver varieties.  
\end{abstract}

\maketitle

\section{Introduction}

Calogero--Moser varieties first appeared as the phase space of the Calogero--Moser Hamiltonian \cite{Calogero}. Much later it was shown by Etingof--Ginzburg \cite{EG} that one can associate to any finite subgroup of the symplectic linear group a flat family of Calogero--Moser varieties, where the original Calogero--Moser phase space is recovered as the Calogero--Moser variety associated to the symmetric group. 

In the world of symplectic representation theory, where interesting representation theory is viewed as modules over quantization algebras of singular symplectic varieties,  the representation theory of semisimple Lie algebras is encoded as modules over quantizations of the nilpotent cone. Replacing the nilpotent cone by the quotient of a symplectic vector space by a finite subgroup of the symplectic linear group, the quantizations one obtains are  spherical symplectic reflection algebras.  Parallel to this, Calogero--Moser varieties appear as deformations of the finite quotient, in the same way that regular (co)adjoint orbit closures in a semisimple Lie algebra deform the nilpotent cone. Just as in the case of regular orbit closures, the geometry of the Calogero--Moser varieties reflect representation theory of spherical symplectic reflection algebras.

Since they are symplectic singularities, Calogero--Moser varieties have a finite stratification by symplectic leaves; this can also be seen in terms of their construction by deforming symplectic quotient singularities, where the leaves correspond to the conjugacy class of stabilizer subgroup.  Their analogues in the case of the nilpotent cone and its deformations are the ubiquitous coadjoint orbits.  In this article we describe both the (\'etale local) singularities transverse to each leaf (``going up'' in the Hasse diagram of the stratification) and the normalization of the closure of each leaf (``going down'' in the aforementioned Hasse diagram). 

Let $\Gamma \subset SL(2,\C)$ be a finite group and $\Gamma_n = \Gamma^n \rtimes \s_n$ the associated wreath product group, acting as a symplectic reflection group on the symplectic vector space $V := (\C^{2})^{n}$. To this one associates a family of symplectic reflection algebras
$\Hb_{c}(\Gamma_n)$ deforming $\C[V] \rtimes \Gamma_n$,
and of 
Calogero--Moser varieties $\ZCB_{c}(\Gamma_n)$, defined as the spectrum of the centres of $\Hb_{c}(\Gamma_n)$. See Section~\ref{sec:symplc-refl-alg} below for details. For each $c$, we:
\begin{itemize}
    \item Parameterize the symplectic leaves $\LG$ of $\ZCB_{c}(\Gamma_n)$.
    \item Describe combinatorially the closure order on leaves (i.e., the Hasse diagram). 
    \item Identify the normalization of each leaf closure with another Calogero--Moser variety. 
    \item Describe the transverse slices to each leaf as a (framed) finite type Nakajima quiver variety. 
\end{itemize}
In \cite{Cuspidal,LosevSRAComplete} it was shown that, to each symplectic leaf in $\ZCB_{c}(\Gamma_n)$, one can attach a ``label'' of a  conjugacy class $(\mr{P})$ of parabolic subgroups $\mr{P}$ of $\Gamma_n$; here ``parabolic'' means the subgroup is the stabilizer subgroup for some vector in $V$. There can be multiple (or no) leaves attached to each conjugacy class of parabolic subgroups. The leaves labeled by $(\mr{P})$ are those induced, in a precise sense, from zero-dimensional leaves in a Calogero--Moser variety for $\mr{P}$. Let $N_{\Gamma_n}(\mr{P})$ be the normalizer of $\mr{P}$ in $\Gamma_n$. Then $\Nq(\mr{P}) = N_{\Gamma_n}(\mr{P})/\mr{P}$  acts symplectically on $V^{\mr{P}}$ and one can consider the corresponding Calogero--Moser variety. Motivated by the conjectural relationship between rational Cherednik algebras at $t = 0$ and Harish-Chandra theory for finite groups of Lie type, it has been conjectured by Bonnaf\'e \cite[Conjecture B]{BonnafeAuto} that, given a (finite) complex reflection group  $W < \GL_n$ and any class function on the set of reflections in $W$, the normalization of the closure of any symplectic leaf on the associated Calogero--Moser variety is isomorphic to a Calogero--Moser variety for the group $N_{W}(\mr{P})/\mr{P}$. The statement of the conjecture makes sense for the more general situation of finite subgroups $W < \Sp_{2n}$ of the symplectic group. 
In the case $W = \Gamma_n < \Sp_{2n}$ of wreath product groups (which are complex reflection groups in the subcase where $\Gamma$ is cyclic), we prove this:

\begin{thm}\label{thm:mainclosureAll}
	If the leaf $\LG$ is labeled by the class $(\mr{P})$ then there exists $w \in W^{\aff}$ such that 
	$$
	\widetilde{\LG} \cong \ZCB_{w^*(c)}(\Nq(\mr{P}),V^{\mr{P}}),
	$$
 where $	\widetilde{\LG}$ is the normalization of the closure of $\LG$ in $\ZCB_{c}(\Gamma_n)$.
\end{thm}

It is difficult to work out explicitly which affine Weyl group element $w$ appears in Theorem~\ref{thm:mainclosureAll}. However, one can read off the parameter $w^*(c)$ (up to the action of the non-affine Weyl group, which does not change the isomorphism class of the variety) directly from the leaf $\LG$; see Proposition~\ref{prop:wbccomp} and Remark~\ref{rem:param-c'-cyclic}.

Let $R$ denote the set of roots in the affine root system associated to $\Gamma$ via the McKay correspondence. Let $\delta \in R^+$ denote the minimal imaginary positive root. One may think of $c$ as a linear functional $\bc$ on the space spanned by the roots $R$. Then there is a natural notion of "level'', which is the scalar $\bc (\delta)$. In order to both state and to prove our results, we must treat the non-zero level $\bc(\delta) \neq 0$ and zero level $\bc(\delta) = 0$ separately. Since level zero is more degenerate, we focus in the introduction on non-zero level. 

\subsection{Symplectic leaves}

At the non-zero level, the root system $R_{\bc} = \{ \alpha \in R \, | \, \bc(\alpha) = 0 \}$ is finite and we denote by $\Delta(\bc) \subset R^+_{\bc}$ the set of simple roots. We define a function $\varrho \colon \N \Delta(\bc) \to \Q$ by $\varrho(\beta) = \beta_0 + \frac{1}{2} (\beta,\beta)$, where $\beta_0$ is the coefficient at the extending simple root and $(-,-)$ is the Cartan pairing, and a partial ordering $\eta \succ \beta$ if and only if $\eta - \beta \in \N \Delta(\bc)$ and $\eta \neq \beta$. Then $\Xi(\bc)$ denotes the set of all $\beta$ in $\N \Delta(\bc)$ such that $\eta \in \N \Delta(\bc)$ with $\eta \succ \beta$ implies that $\varrho(\eta) > \varrho(\beta)$. 

\begin{thm}\label{thm:nonzerolevelleavesI}
There is a bijection between $\{ \beta \in \Xi(\bc) \, | \, \qu (\beta) \le n \}$ and the symplectic leaves of $\ZCB_{c}(\Gamma_n)$, $\beta \mapsto \LG(\beta)$, such that 
	\begin{enumerate}
		\item[(i)] $\dim \LG(\beta) = 2m$, and
		\item[(ii)] the leaf $\LG(\beta)$ is labeled by the parabolic conjugacy class $(\Gamma_m)$,
	\end{enumerate}    
where $m := n - \qu(\beta)$. 
\end{thm}

Next, we consider the closure ordering on the set of symplectic leaves in $\ZCB_{c}(\Gamma_n)$.  

\begin{prop}\label{prop:leafclosure1I}
	For $\beta, \eta \in \Xi(\bc)$ with $\beta \neq \eta$, we have $\LG(\eta) \subset \overline{\LG(\beta)}$
    if and only if $\eta \succ \beta$. 
\end{prop}
Observe that, since $\overline{\LG(\beta)}$ is automatically closed under the Hamiltonian flow, 
the condition $\LG(\eta) \subset \overline{\LG(\beta)}$ is equivalent to  $\LG(\eta) \cap \overline{\LG(\beta)} \neq \emptyset$.

A natural question arising from Theorem~\ref{thm:mainclosureAll} is whether the closure of a leaf is a normal variety. We show: 

\begin{prop}
When $\bc(\delta) \neq 0$, each leaf closure $\overline{\LG} \subset \ZCB_{c}(\Gamma_n)$ is normal.  
\end{prop}

When $\bc(\delta) = 0$, it is easily seen that leaf closures in $\ZCB_{c}(\Gamma_n)$ are not generally normal. 

\subsection{Transverse slices}

Under the assumption $\bc(\delta) \neq 0$, one can show (Lemma~\ref{lem:wparabolicrootsystem}) that there exists an element $w$ in the affine Weyl group $W^{\aff}$ such that $R_0 := w(R_{\bc})$ is a parabolic root subsystem of $R$; this means that $R_0$ is a root system generated by a subset of the simple roots $\Delta \subset R^+$. The latter corresponds to a (finite type, possibly disconnected) subgraph $\mathsf{G}(c)$ of the affine Dynkin diagram $\mathsf{G}(\Gamma)$. 

\begin{thm}\label{thm:mainsingAll}
	Let $\LG$ be a symplectic leaf of $\ZCB_{c}(\Gamma_n)$. There exist dimension vectors $\bw$ and $\bv$ such that the transverse slice to $\LG$ is isomorphic to the (framed) quiver variety $\mf{M}_0(\mathsf{G}(c),\bw,\bv)$. 
\end{thm}

We explain in Section~\ref{sec:nonzerotransverse} how to compute $\mathsf{G}(c)$, $\bw$, and $\bv$ directly from the root $\beta$ labeling the leaf $\LG$. In most cases, the quiver variety $\mf{M}_0(\mathsf{G}(c),\bw,\bv)$ is not isomorphic to any Calogero--Moser variety.

\begin{rem}\label{rem:manyslices}
	As a converse to Theorem~\ref{thm:mainsingAll} we note the following. 
	\begin{enumerate}
		\item Let $\mathsf{G}$ be a finite type simply-laced graph (i.e. ADE graph) and $\Gamma$ the finite subgroup of $\SL(2,\C)$ whose McKay graph is the affine Dynkin graph $\widetilde{\mathsf{G}}$. For any pair $(\bw,\bv)$ of dimension vectors for $\mathsf{G}$, we show that there exists $n, c$ and a leaf $\LG$ such that the framed Nakajima quiver variety $\mf{M}_0(\mathsf{G},\bw,\bv)$ (of finite type) can be realized as the transverse slice to $\LG$ in $\ZCB_{c}(\Gamma_n)$. 
  
		\item As a special case, if we take $\Gamma = \Z / {\ell} \Z$ and let $S(\mu,\nu)$ be any Slodowy slice of type $\mathsf{A}$, then for $n,\ell$ sufficiently large one can always find a parameter $c$ and a leaf $\LG \subset \ZCB_{c}(\Gamma_n)$ such that the singularity transverse to $\LG$ is isomorphic to $S(\mu,\nu)$. Thus, each type $\mathsf{A}$ Slodowy slice occurs as a singularity in some Calogero--Moser variety. Equivalently, one could talk about slices in the affine Grassmaniann of type $\mathsf{A}$ occurring in Calogero--Moser varieties of type $\mathsf{A}$.  
	\end{enumerate}
 See Proposition~\ref{prop:transversearbitrary} and Remark~\ref{rem:rangeofslices} for details. 
\end{rem}

\subsection{The cyclic group}

Among the finite subgroups of $SL(2,\C)$, the cyclic group is distinguished in that the associated Calogero--Moser variety has an important additional symmetry. Namely, there is a Hamiltonian $\Cs$-action with finitely many fixed points. This means that we can give a more explicit combinatorial description of the leaves in this case. The combinatorics that arise are important in the representation theory of restricted rational Cherednik algebras and the conjectural links to certain finite groups of Lie type. 

Let $\Gamma=\ZM/\ell\ZM$ be a cyclic subgroup of $\CM^\times\subset SL(2)$. In this case, the Calogero--Moser variety $\ZCB_{c}(\Gamma_n)$ can be described as the quiver variety $\mf{M}_\bc(\Qu^\ell, \Lambda_0, n\delta)$, where $\Qu^\ell$ is the cyclic quiver of length $\ell$. To have a good combinatorial description of this variety, we would prefer to replace $\bc$ by a parameter whose stabilizer in $\widehat\s_\ell$ is a parabolic subgroup. By applying a sequence of admissible reflections to $\bc$, we get a different parameter $\bc'$ whose stabilizer in $\widehat\s_\ell$ is a parabolic subgroup $W_J$ for some parabolic type $J\subset \ZM/\ell\ZM$. This sequence of reflections replaces the dimension vector $n\delta$ with some other dimension vector $\alpha$. Finally, our Calogero--Moser variety can be described as the quiver variety $\mf{M}_{\bc'}(\Qu^\ell, \Lambda_0, \alpha)$. We can write $\alpha = n\delta + \Res_\ell(\nu)$, where $\nu$ is an $\ell$-core (i.e., $\nu$ is a partition without removable $\ell$-hooks), and $\Res_\ell(\nu)$ the $\ell$-residue of $\nu$.

It is known due to \cite[Prop.~8.3~(i)]{GordonQuiver} that the $\CM^\times$-fixed points of this variety are labeled by the $J$-cores of elements of $\PC_\nu(n\ell+|\nu|)$, where $\PC_\nu(n\ell+|\nu|)$ is the set of partitions of $n\ell+|\nu|$ whose $\ell$-core is $\nu$; a $J$-core of a partition is the partition obtained from it by removing all possible removable boxes with residues in $J$. For each $J$-core as above, we give an explicit construction of a quiver representation giving the corresponding $\CM^\times$-fixed point.

Moreover, we provide a combinatorial construction of the symplectic leaves of the Calogero--Moser variety, showing that they are parameterized by $\ell$-cores of $J$-cores of elements of $\coprod_{n'=0}^n\PC (n'\ell + |\nu|)$. We show that the map sending the $\CM^\times$-fixed point to the symplectic leaf containing it corresponds combinatorially to the $\ell$-core map. This shows in particular that the symplectic leaves containing at least one $\CM^\times$-fixed point are those coming from $n' = n$.

\subsection{The hyperoctahedral group}

Calogero--Moser varieties associated to finite Coxeter groups are important 
because it is expected that their geometry is related to Harish-Chandra theory for finite groups of Lie type \cite{BonnafeConjectures}. Two infinite families are particularly interesting: dihedral groups and Weyl groups of type $\mathsf{B}$, since these are the infinite families for which there exist "unequal parameters''; equivalently, for which there are both long and short roots. The dihedral groups are studied in \cite{BonnafeDihedral1,BonnafeDihdral3}. As an extended example, we explain in greater detail what our results mean for the infinite family of Weyl groups of type $\mathsf{B}$. Thus, $\Gamma = \Z /2 \Z$ and $\Gamma_n = W(\mathsf{B}_n)$. Here, the parameters are $c = (c_1,c_{\gamma}) \in \C^2$.

As in the general situation, the singularities appearing in $\ZCB_{c}(\B_n)$ fall into two distinct families. When $\bc(\delta) = c_1 \neq 0$, they are certain nilpotent orbit closures in $\mf{gl}(N)$ and when $\bc(\delta) = c_1 = 0$ they are symplectic quotient singularities for products of symmetric groups. Based on the underlying ext-quiver that describes the local singularities, one may think of these as finite and affine type $\mathsf{A}$ situations respectively. 

Let $\Gr (k,N)$ denote the Grassmanian of $k$-planes in $\C^N$. Then $T^* \Gr (k,N)$ is a symplectic resolution of 
$$
\overline{\mc{O}(k,N)} = \Spec \C[T^* \Gr(k,N)].
$$
where $\mc{O}(k,N)$ is the nilpotent orbit in $\mf{gl}(N)$ consisting of all matrices $X$ of rank $k$ with $X^2 = 0$. When $2 k \le N$, it is the nilpotent orbit labeled by the partition $(2^k,1^{N- 2 k})$. Moreover,  
$$
\overline{\mc{O}(k,N)} = \mc{O}(k,N) \sqcup \mc{O}(k-1,N) \sqcup \cdots \sqcup \mc{O}(0,N)
$$
shows that the closure relation is a total order. 

We denote set of all partitions of $n$ by $\mc{P}(n)$. If $\lambda=(\lambda_1,\lambda_2,\ds)$ is such a partition, set $\s(\lambda) := \prod_{i \ge 1} \s_{n_i}$ with $n_i := | \{ k \, | \, \lambda_k = i \} |$.
Denote by $\s_{\lambda} \subset \s_n$ is the parabolic subgroup of the symmetric group $\s_n$ labeled by $\lambda$. The length of the partition $\lambda$ is denoted $\ell(\lambda)$.

\begin{thm}\label{thm:mainsing} 
	Let $c = (c_1,c_{\gamma})$ be non-zero.
	\begin{enumerate}
		\item[(i)] If $c_1 = 0, c_{\gamma} \neq 0$, and $p \in \LG_{\lambda}$ for $\lambda \in \mc{P}(n)$, then 
		$$
		(\ZCB_{c}(\B_n),p) \cong \left( (\h \times \h^*) / \s_{\lambda} \times T^* \C^{\ell(\lambda)}, (0,0) \right),
		$$
		where $\h$ is the reflection representation for $\s_{\lambda}$.
		\item[(ii)] If $c_{\gamma} = m c_1$, with $c_1 \neq 0$ and $1 \le m \le n-1$, and $p \in \LG_k$ for $k(k + m) \le n$, then
		$$
		(\ZCB_{c}(\B_n),p) \cong \left( \overline{\mc{O}(k,2k + m)} \times T^*( \C^{n - k(k+ m)}), (0,0) \right). 
		$$

        \item[(iii)] If $c_{\gamma} = -m c_1$, with $c_1 \neq 0$ and $1 \le m \le n-1$, and $p \in \LG_k$ for $k(k + m) \le n$, then
		$$
		(\ZCB_{c}(\B_n),p) \cong \left( \overline{\mc{O}(k,2k + m)} \times T^*( \C^{n - k(k+ m)}), (0,0) \right). 
		$$
        
	\end{enumerate}
    For all other values of $c$, $\ZCB_{c}(\B_n)$ is smooth. 
\end{thm}

The proof is, to a large extent, an application of Martino's thesis \cite{MartinoThesis}. It will be given in Section~\ref{sec:mainsingproof} below. 

\begin{thm}\label{thm:mainclosureB} 
	Let $c = (c_1,c_{\gamma})$ be non-zero.
	\begin{enumerate}
		\item[(i)] If $c_1 = 0$, and $\rho \in \mc{P}(n)$, then 
		$$
		\widetilde{\LG}_{\rho} \cong \ZCB_{c'}(\Z_2 \wr \s(\rho)),
		$$
		where $c'(t) = 1$ for all reflections $t \in (\Z / 2 \Z)^{\ell}$ and $c'(s) = 0$ for all reflections $s$ in $\s(\rho)$.
		\item[(ii)] If $c_{\gamma} = m c_1$, with $c_1 \neq 0$, $1 \le m < n$ and $k \ge 0$ with $k(k + m) \le n$ then
		$$
		\widetilde{\LG}_{k} \cong \ZCB_{c'}(\B_{n - k(k+m)}).
		$$
		where $c_1' = 1$ and $c_{\gamma}' = m+2k$. 
		\item[(iii)] If $c_{\gamma} = -m c_1$, with $c_1 \neq 0$, $0 \le m < n$ and $k \ge 0$ with $k(k + m) \le n$ then
		$$
		\widetilde{\LG}_{k} \cong \ZCB_{c'}(\B_{n - k(k+m)}),
		$$
        where $c_1' = 1$ and $c_{\gamma}' = m+2k$. 
	\end{enumerate}
\end{thm}

In part (i), the Calogero--Moser variety $\ZCB_{c'}((\Z / 2 \Z) \wr \s(\rho))$ comes from the natural action of $(\Z / 2 \Z)^{\ell(\rho)} \rtimes \s(\rho)$ on $\C^{\ell(\rho)}$. 

\subsection{Method of proof}

As the reader can already see from the above theorems, describing the singularities of the Calogero--Moser varieties $\ZCB_{c}$ is tractable because these varieties are isomorphic to certain (framed) affine quiver varieties. Such isomorphisms only appear to exist for the wreath product groups. They were first constructed when $\ZCB_c(\Gamma_n)$ is smooth in \cite{EG} and it was shown in \cite{MarsdenWeinsteinStratification} that the isomorphisms extend to the the case where the Calogero--Moser variety is singular. The isomorphisms are Poisson and hence identify symplectic leaves. Then we can leverage the finite stratification of quiver varieties by representation type constructed by Crawley-Boevey \cite{CBmomap} (and identified with the stratification by symplectic leaves in \cite{MarsdenWeinsteinStratification,BellSchedQuiver}). This gives the parameterization and ordering on leaves in $\ZCB_c(\Gamma_n)$. Finally, we use Crawley-Boevey's \'etale local normal form \cite{CBnormal} using ext-quivers (the hyperk\"ahler version of this result is described \cite[Section~6]{Nak1994}), to describe the transverse slices to leaves.   

\subsection{Leaf closures in quiver varieties}

Since we are using throughout the fact that $\ZCB_{c}(\Gamma_n)$ is a quiver variety, it is useful to consider more generally the geometry of leaf closures in arbitrary quiver varieties.

Let $\mf{M}_{\bc}(\alpha)$ be the (unframed) affine quiver variety with dimension vector $\alpha$ and deformation parameter $\bc$ associated to the graph $\mathsf{G}$. The closed points of $\mf{M}_{\bc}(\alpha)$ parameterize isomorphism classes of semi-simple representation of the deformed preprojective algebra $\Pi^{\bc}$ of dimension $\alpha$.

Let $\Sigma_{\bc}$ denote the set of all dimension vectors for which there exists a simple $\Pi^{\bc}$-module. As recalled in Section~\ref{sec:sympleavesdefn}, the symplectic leaves $\LG_{\tau}$ of $\mf{M}_{\bc}(\alpha)$ are the representation type strata and are labeled by decompositions of $\alpha$:
$$
\tau = (n_1,\alpha^{(1)}; \ds ; n_k,\alpha^{(k)})
$$
where $\alpha^{(i)} \in \Sigma_{\bc}$, $\sum_i n_i \alpha^{(i)} = \alpha$ and each real root occurs at most once. Each closed point in the stratum $\mf{M}_{\bc}(\alpha)_{\tau}$ corresponds to a semisimple representation $M = M_1^{\oplus n_1} \oplus \cdots \oplus M_k^{\oplus n_k}$, where the $M_i$ are pairwise non-isomorphic simple $\Pi^{\bc}$-modules and $\dim M_i = \alpha^{(i)}$. We think of $\tau$ as a function from the set $\Sigma_{\bc}$ into $\mc{P}$, the set of all partitions, such that 
$$
\sum_{\beta \in \Sigma_{\bc}} |\tau(\beta)| \, \beta = \alpha.
$$ 
To any such function we associate a product of symmetric groups 
$$
\s_{\tau} = \prod_{\beta \in \Sigma_{\bc}} \s(\tau(\beta)).
$$

\begin{thm}\label{thm:quiverleafclosure}
	There is a morphism $\prod_{\beta \in \Sigma_{\bc}} \mf{M}_{\bc}(\beta)^{\ell(\tau(\beta))}/\s(\tau(\beta)) \to \mf{M}_{\bc}(\alpha)$ whose image equals $\overline{\LG}_{\tau}$. The resulting map 
 \[
 \prod_{\beta \in \Sigma_{\bc}} \mf{M}_{\bc}(\beta)^{\ell(\tau(\beta))}/\s(\tau(\beta)) \to \overline{\LG}_{\tau}
 \]
 is the normalization of $\overline{\LG}_{\tau}$. 
\end{thm}

\begin{remark}
	We remark that, technically, the space $\prod_{\beta \in \Sigma_{\bc}} \mf{M}_{\bc}(\beta)^{\ell(\beta)}/\s(\tau(\beta))$ is not a quiver variety, although it is if all of the $\beta$ that occur are either real or isotropic imaginary roots (geometrically meaning that $\mf{M}_{\bc}(\beta)$ has dimension $0$ or $2$), since for isotropic imaginary $\beta$ we have $\mf{M}_{\bc}(\beta)^n / \s_n = \mf{M}_{\bc}(n\beta)$ by \cite{CBdecomp}.
\end{remark}

\subsection{Outline of the article}

In Section~\ref{sec:quivervarieties} we recall the background results on quiver varieties that we require later. We also prove the result describing the normalization of a leaf closure in an arbitrary quiver variety. Section~\ref{sec:symplc-refl-alg} then introduces Calogero--Moser varieties and explains the labeling of leaves by parabolic subgroups and the identification with quiver varieties. Sections~\ref{sec:nonzero} and \ref{sec:zerolevel} consider the leaf closures (ordering and identification with another Calogero--Moser varieties) in the non-zero and zero levels respectively. Then, we recall in Section~\ref{sec:combinatoricsCM} some background material on combinatorics before turning to the case of $\Gamma$ a cyclic group in Section 7. Section 8 describes the transverse slices to each leaf. Finally, Section 9 describes our results in greater detail in the case of the Weyl group of type $\mathsf{B}$.

\subsection*{Acknowledgments}

We would like to thank C\'edric Bonnaf\'e and Daniel Juteau for encouraging us to write this article and for many interesting conversations. 

The first author was supported in part by Research Project Grant
RPG-2021-149 from The Leverhulme Trust and EPSRC grants
EP-W013053-1 and EP-R034826-1. 

\section{Quiver varieties}\label{sec:quivervarieties}

\subsection{Notation}\label{sec:notation}

 Throughout, $\N := \{ 0,1,2,\ds \} = \Z_{\ge 0}$.

\begin{defn}
A \emph{partition} is a tuple $\lambda=(\lambda_1,\lambda_2,\dots,\lambda_r)$ of positive integers 
	(with no fixed length) such that $\lambda_1\geqslant\lambda_2\geqslant\dots\geqslant\lambda_r > 0$, $r\geqslant 0$. 
	Set $|\lambda|=\sum_{i=1}^r \lambda_i$ and $\ell(\lambda)=r$. If $|\lambda|=n$, we say that $\lambda$ is a partition of $n$. Denote by $\PC$ (resp. $\PC(n)$) the set of all partitions (resp. the set of all partitions of $n$). By convention, $\PC(0)$ contains one (empty) partition, with $\ell(\emptyset) = 0$. For $\lambda\in\PC(n)$, let $\s_\lambda=\prod_{t=1}^r\s_{\lambda_t}$ be the corresponding parabolic subgroup of the symmetric group $\s_n$. We also set $\s(\lambda) := \prod_{i \ge 1} \s_{n_i}$ where $n_i := | \{ k \, | \, \lambda_k = i \} |$. 
 \end{defn}

\subsection{Graphs}
\label{subs:graphs}

We fix a finite unoriented graph $\mathsf{G}$ with vertex set $I = \{ 0,1, \ds, r \}$. Let $H$ be the set of pairs consisting of an edge together with a choice of orientation of that edge. For $h \in H$, write $s(h)$ for the source vertex of $h$ and $t(h)$ for the target vertex. To this graph we associate a \textit{root lattice} $Q = \bigoplus_{i \in I} \Z e_i$ and a \textit{weight lattice} $P = \bigoplus_{i \in I} \Z \Lambda_i$. We set 
\[
P^+ = \sum_{i = 0}^r \Z_{\ge 0} \Lambda_i, \quad Q^+ = \sum_{i = 0}^r \Z_{\ge 0} e_i.
\]
For $\alpha \in Q$, we write throughout $\alpha = \sum_{i \in I} \alpha_i e_i$, where $\alpha_i \in \Z$. If $a_{i,j}$ denotes the number of edges between vertices $i,j \in I$ then we define a symmetric bilinear form $( - , - )$ on $Q$, called the Cartan pairing, by 
\[
(e_i,e_j) = 2 \cdot 1_{i,j} - a_{i,j},
\]
where $1_{i,j}$ is the Kronecker delta. If $i$ is a loopfree vertex then there is a reflection $s_i \colon Q \to Q$ defined by $s_i(\alpha) = \alpha - (\alpha, e_i) e_i$. 

There is also the dual reflection, 
\[
s_i^* \colon \Hom_{\Z}(Q,\C) \to \Hom_{\Z}(Q,\C).
\]
For an element $\bc\in\Hom_{\Z}(Q,\C)$, we abbreviate $\bc_i=\bc(e_i)$. In particular, we identify $\Hom_{\Z}(Q,\C)$ with $\CM^I$. The dual reflection $s_i^*$ is given by $(s_i^* \bc)_j = \bc_j - (e_i,e_j)\bc_i$.  The subgroup $W_{\mathsf{G}}$ of the group of automorphisms of the abelian group $Q$ generated by all reflections $s_i$ is the Weyl group associated to the graph. 
For each $\bw=\sum_{i\in I}\bw_i\Lambda_i \in P$, there is also a twisted action of $W_{\mathsf{G}}$ on $Q$ given by 
\begin{equation}\label{eq:twistedaction}
    s_i \star_{\bw} \alpha = s_i(\alpha) + \bw_i e_i. 
\end{equation}

\begin{remark}\label{rem:noloopaffineLie}
When the graph contains no loops, one can associate to $\mathsf{G}$ a Kac--Moody 
algebra $\mf{g}$ \cite{KacBook} with Cartan subalgebra $\mf{h}$. Then the $e_i$ are the simple roots and the $\Lambda_i$ are the fundamental weights, so that $Q$ is the root lattice, and $P$ is the lattice of weights of the derived subalgebra $\mf{g}':=[\mf{g},\mf{g}]$.  As a result we have a natural inclusion $Q \subset \mathfrak{h}^*$, and if we make a choice of complement $\mf{d}$ of $\mathfrak{h} \cap \mf{g}'$ in $\mathfrak{h}$, this induces another natural inclusion $P \subset \mf{d}^\perp \subset \mathfrak{h}^*$, with $P+Q$ a full lattice in $\mathfrak{h}^*$. 
Moreover,  the (possibly degenerate) Cartan pairing on $Q$ extends to a nondegenerate pairing on $\mathfrak{h}^*$ such that $(\Lambda_i, e_j)=\delta_{ij}$, and the action of $W_{\mathsf{G}}$ extends to $\mathfrak{h}^*$, defined by
\[
s_i (\bw) = \bw - (\bw, e_i) e_i=\bw - \bw_i e_i.
\]
Then $s_i(\bw - \alpha) = \bw - s_i \star_{\bw} \alpha$. Alternatively, the action $\star_{\bw}$ can be viewed as coming from the usual action on dimension vectors of the Weyl group of the deframed graph; see Section~\ref{sec:CBtrick} and \cite[Definition~2.3]{NakajimaWeyl}. 
\end{remark}

\subsection{Root systems}\label{root-sys}

Set $p(\alpha) = 1 - (1/2) (\alpha, \alpha)$. The support of a vector $\alpha \in Q$ is the full subgraph whose vertices are $\{ i \in I \, | \, \alpha_i \neq 0 \}$. The fundamental region $F$ is the set of all $\alpha \in Q^+$ with connected support and with $(\alpha, e_i) \le 0$ for all $i$. The real roots (respectively imaginary roots) are the elements of $Q$ which can be obtained from the coordinate vector at a loopfree vertex (respectively $\pm$ an element of the fundamental region) by applying some sequence of reflection at loopfree vertices. Recall that a root $\beta$ is isotropic imaginary if $p(\beta) = 1$ (equivalently $(\beta, \beta) = 0$) and anisotropic imaginary if $p(\beta) > 1$. Abusing terminology, we will simply say that a root $\alpha$ is (a) real if $p(\alpha) = 0$, (b) isotropic if $p(\alpha) = 1$, and (c) anisotropic if $p(\alpha) > 1$.

We recall the following important set, first defined by Crawley-Boevey. Fix $\bc \in \Hom(Q,\C)$. Let $\Sigma_{\bc}$ denote the set of all positive roots $\alpha$ such that $\bc(\alpha) = 0$ and 
\[
p(\alpha) > p(\beta^{(0)}) + \cdots + p(\beta^{(k)})
\]
for all proper decompositions $\alpha = \beta^{(0)} + \cdots + \beta^{(k)}$ into a sum of positive roots $\beta^{(i)}$ also satisfying $\bc(\beta^{(i)}) = 0$. Then \cite[Theorem~1.2]{CBmomap} says that there exists a simple representation of dimension $\alpha$ for the deformed preprojective algebra $\Pi^{\bc}$ associated to $\mathsf{G}$ if and only if $\alpha \in \Sigma_{\bc}$. 

\subsection{Quiver varieties}\label{sec:quivernotation} 

We fix a pair $(\bw,\alpha) \in P^+ \times Q^+$. Let $(V_i)_{i \in I}$ and $(W_i)_{i \in I}$ be tuples of complex vector spaces with $\dim V_i = \alpha_i$, $\dim W_i = \bw_i$. Consider the action of the group $G(\alpha) := \prod_{i \in I} GL(\alpha_i)$ on the space  
\begin{equation}\label{eq:repspace}
    \Rep(\bw,\alpha) := \bigoplus_{h \in H} \Hom_{\C}(V_{s(h)}, V_{t(h)}) \bigoplus_{i \in I} \Hom_{\C}(W_i,V_i) \bigoplus_{j \in I}  \Hom_{\C}(V_j,W_j). 
\end{equation}
There is an involution $(-)^*$ on $H$ such that $h^*$ has the same underlying edge as $h$ but opposite orientation. Fix $\epsilon \colon H \to \{ \pm 1 \}$ such that $\epsilon(h^*) = - \epsilon(h)$. The subset of oriented edges (arrows) $h$ in $H$ with $\epsilon(h) = 1$ form a quiver $\Qu$ whose double $\overline{\Qu}$ has arrow set $H$. The function $\epsilon$ specifies a symplectic form on $\Rep(\bw,\alpha)$, making the action of $G(\alpha)$ Hamiltonian. We identify $\mf{g}(\alpha)$ with its dual using the trace form so that the moment map $\mu \colon \Rep(\bw,\alpha) \to \mf{g}(\alpha)$ for this action, uniquely specified by $\mu(0) = 0$, is given by 
$$
\mu(X,v,w) = \sum_{h \in H} \epsilon(h) X_{h^*} X_h + w v.
$$
We will sometimes write $\mu_\alpha$ instead of $\mu$ to specify $\alpha$. Let $\bc \in \Hom_{\Z}(Q,\C) \cong \C^I$, identified with the tuple of scalar matrices $(\bc_i \Id_{V_i})_{i \in I} \in \mf{g}(\alpha)$. The framed quiver variety is
\[
\mf{M}_{\bc}(\bw,\alpha) = \mu^{-1}(\bc) \git \, G(\alpha). 
\]
In the above, we equip $\mf{M}_{\bc}(\bw,\alpha)$ with the reduced scheme structure. Then $\mf{M}_{\bc}(\bw,\alpha)$ is an irreducible normal affine variety \cite{CBmomap,CBdecomp,CBnormal}. Except for sections~\ref{sec:slices} and~\ref{sec:hyperoctahedral}, we will be exclusively interested in the case where $\bw = \Lambda_0$. In this case, we write
\[
\XC_{\bc}(\alpha) := \mf{M}_{\bc}(\Lambda_0,\alpha) \label{eq:quivernotation1} 
\]
for brevity. When $\bw = 0$, we write $\mf{M}_{\bc}(\alpha) = \mf{M}_{\bc}(0,\alpha)$ and $\Rep(0,\alpha)=\Rep(\alpha)$. We may also sometimes write $\Rep(\overline\Qu,\alpha)$ instead of $\Rep(\alpha)$ to recall that this is the representation space for the double quiver $\overline\Qu$. The definition of $\mf{M}_{\bc}(\bw,\alpha)$ depends on the choice of function $\epsilon$ (equivalently, on the quiver $\Qu$), though there is a canonical isomorphism between the spaces $\mf{M}_{\bc}(\bw,\alpha)$ defined using different $\epsilon$. When we need to indicate the dependence of $\mf{M}_{\bc}(\bw,\alpha)$ on $\Qu$ we write $\mf{M}_{\bc}(\Qu,\bw,\alpha)$ and when we only need to specify $\mathsf{G}$, we'll write $\mf{M}_{\bc}(\mathsf{G},\bw,\alpha)$.

\subsection{Deframing}\label{sec:CBtrick}

If $\bw \neq 0$, one can always use Crawley-Boevey's trick to "deframe" the framed quiver and realize $\mf{M}_{\bc}(\bw,\alpha)$ as an unframed quiver variety $\mf{M}_{{\bc}_{\alpha}}(e_{\infty} + \alpha)$ associated to a new graph $\mathsf{G}'$ where the new vertex set $I'$ has one additional vertex ${\infty}$ and there are $\bw_i$ additional edges between the vertex ${\infty}$ and vertex $i$ and no loops at $\infty$. Here ${\bc}_{\alpha}(e_i) = \bc(e_i)$ for $i \neq \infty$ and ${\bc}_{\alpha}(e_{\infty}) = - \bc (\alpha)$. Moreover, we have defined the action $\star_{\bw}$ in \eqref{eq:twistedaction} precisely so that $s_i(e_{\infty} + \beta) = e_{\infty} + s_i \star_{\bw} \beta$ for $i \in I \subset I'$ and $\beta \in Q \subset Q'$. We define $E^{\bw}_{\bc} \subset Q^+$ to be the set of dimension vectors $\gamma\in Q^+$ such that $e_{\infty} + \gamma \in \Sigma_{{\bc}_{\gamma}}(\mathsf{G}')$.  

When the graph $\mathsf{G}$ has no loops, we can associate to it (as in Remark~\ref{rem:noloopaffineLie}) a Kac--Moody Lie algebra $\mf{g}$ with Cartan subalgebra $\mf{h}$ such that $Q \subset \mf{h}^*$ is the root lattice of $\mf{g}$, and $P$ is the lattice of weights of the derived subalgebra $\mf{g}'=[\mf{g},\mf{g}]$. Up to automorphism of $\mf{g}$, there is a unique way to view $P$ also as a subset of the dual $\mf{h}^*$ of the Cartan subalgebra, such that $P+Q$ spans $\mf{h}^*$.  Thus sums in $P+Q$ can be interpreted as weights of $\mf{g}$.
With this in mind, the dimension vector $\bw \in P^+$ labels an integrable highest weight $\mf{g}$-module $L(\bw)$. For a fixed $\bw \in P^+$, we define the quadratic function 
\[
p_{\bw}(\beta) := p(\beta) + {\bw \cdot \beta} - 1, \quad \forall \, \beta \in Q.
\]
Then $p_{\bw}(\beta) = p(e_{\infty} + \beta)$. It follows from \cite[Theorem~2.15]{NakajimaBranching} that we can equivalently define $E_{\bc}^{\bw}$ to be set of $\gamma \in Q^+$ such that 
\begin{enumerate}
    \item[(i)] $\bw - \gamma$ is a weight of $L(\bw)$;
    \item[(ii)] $p_{\bw}(\gamma) > p_{\bw}(\beta^{(0)}) + p(\beta^{(1)}) + \cdots + p(\beta^{(k)})$, where $\gamma = \beta^{(0)} + \beta^{(1)} + \cdots + \beta^{(k)}$ is a proper decomposition with $\beta^{(0)} \in Q^+$, $\bw - \beta^{(0)}$ a weight of $L(\bw)$ and $\beta^{(i)} \in \Sigma_{\bc}(\mathsf{G})$ for $i \ge 1$. 
\end{enumerate}
This is because $e_{\infty} + \gamma$ is a root for $\mathsf{G}'$ if and only if $\bw - \gamma$ is a weight of the representation $L(\bw)$. We note that $0 \in E_{\bc}^{\bw}$ but $0$ never belongs to $\Sigma_{\bc}$. 

When $\bw = \Lambda_0$, we drop $\bw$ from the notation and write $E_{\bc}=E_{\bc}^{\Lambda_0}$.

\subsection{The symplectic leaves}\label{sec:sympleavesdefn}

The quiver variety $\mf{M}_{\bc}(\bw,\alpha)$ has a finite stratification by (locally closed, smooth) symplectic leaves. This stratification equals the stratification by representation type \cite[Theorem~1.9]{BellSchedQuiver}. To explain what this means, we first recall that a \textit{decomposition} $\tau$ of $\alpha$ with respect to $(\bw,\bc)$ is tuple 
\begin{equation}\label{eq:taureptype}
	\tau = (\beta^{(0)},(n_1,\beta^{(1)}; \ds;n_k,\beta^{(k)})),
\end{equation}
where $\beta^{(0)} \in E^{\bw}_{\bc}$, $\beta^{(i)} \in \Sigma_{\bc}$ for $i \ge 1$ and $\alpha = \beta^{(0)} + \sum_{i = 1}^k n_i \beta^{(i)}$. Note that the $\beta^{(i)}$ need not be pairwise distinct. In the case $\bw = 0$, the term $\beta^{(0)}$ is omitted and we think of $\tau$ as a function $\tau \colon \Sigma_{\bc} \to \mc{P}$ such that $\sum_{\beta \in \Sigma_{\bc}} |\tau(\beta)|\beta = \alpha$. Here $\tau(\beta) := (n_{i_1} \ge n_{i_2} \ge \dots)$ if $i_1,i_2,\dots$ are all indices such that $\beta^{(i_j)} = \beta$. 

The leaf $\LG_{\tau}$ labeled by the decomposition $\tau$ parameterizes all isomorphism classes $[M]$ of representations $M = M_0 \oplus \left(\bigoplus_{i = 1}^k M_i^{\oplus n_i} \right)$, where $M_i$ for $i\ne 0$ is a simple representation of the deformed preprojective algebra $\Pi^{\bc_{\alpha}}(\mathsf{G}')$ of dimension $\beta^{(i)}$, $M_0$ the representation of dimension $e_\infty+\beta^{(0)}$, and $M_i \not\cong M_j$ for $i \neq j$. If a real root occurs more than once in $\tau$ then $\LG_{\tau} = \emptyset$. Otherwise, it is non-empty. 

We note that stratification by symplectic leaves satisfies the frontier condition: $\LG \cap \overline{\mc{M}} \neq \emptyset$ implies that $\LG \subset \overline{\mc{M}}$. 

\subsection{Admissible reflections}
\label{subs:adm-refl}

A reflection $s_i$ is said to be \textit{$\bc$-admissible} if $\bc_i \neq 0$. In this case, it was shown by Maffei \cite{Maffei} (and Nakajima \cite{NakajimaWeyl} in the hyperk\"ahler setting) that there is an isomorphism $\mf{M}_{\bc}(\Qu,\alpha) \cong \mf{M}_{s_i^*(\bc)}(\bw,s_i \star_{\bw} \alpha)$, where $\mf{M}_{s_i^*(\bc)}(\bw,s_i \star_{\bw} \alpha) = \emptyset$ if $s_i \star_{\bw} \alpha \notin Q^+$. Moreover, by \cite[Lemma~6.4.3]{Losev} this isomorphism is Poisson. 

More generally, we say that $w \in W_{\mathsf{G}}$ is $\bc$-admissible if there exists a decomposition $w = s_{i_k} \cdots s_{i_1}$ such that $s_{i_j}$ is $s_{i_{j-1}}^* \cdots s_{i_1}^*(\bc)$-admissible. In this case, $\mf{M}_{\bc}(\bw,\alpha) \cong \mf{M}_{w^*(\bc)}(\bw,w \star_{\bw} \alpha)$. 

If $\bc' \in W_{\mathsf{G}} \cdot \bc$ and $w \in W_{\mathsf{G}}$ an element of minimal length with the property that $w^*(\bc) = \bc'$ then \cite[Corollary~5.2]{CrawleyBoeveyHolland} says that $w$ is $\bc$-admissible.

\subsection{Leaf closures}

In this section we assume $\bw = 0$. Note that, as explained in Section~\ref{sec:CBtrick}, we can always deframe a framed quiver variety. Therefore, our results below also apply to framed quiver varieties. Recall that for a decomposition $\tau$, $\LG_{\tau}$ is the symplectic leaf labeled by $\tau$. 

We give the proof of Theorem~\ref{thm:quiverleafclosure}, which we recall below. 

\begin{thm}\label{thm:quiverleafclosureagain}
	There is a morphism $\prod_{\beta \in \Sigma_{\bc}} \mf{M}_{\bc}(\beta)^{\ell(\tau(\beta))}/\s(\tau(\beta)) \to \mf{M}_{\bc}(\alpha)$ whose image equals $\overline{\LG}_{\tau}$. The resulting map 
 \[
 \prod_{\beta \in \Sigma_{\bc}} \mf{M}_{\bc}(\beta)^{\ell(\tau(\beta))}/\s(\tau(\beta)) \to \overline{\LG}_{\tau}
 \]
 is the normalization of $\overline{\LG}_{\tau}$. 
\end{thm}

We begin by noting that it follows from the main result of \cite{CBnormal} that $\prod_{\beta \in \Sigma_{\bc}} \mf{M}_{\bc}(\beta)^{\ell(\beta)}$ is normal. Hence the product $\prod_{\beta \in \Sigma_{\bc}} \mf{M}_{\bc}(\beta)^{\ell(\beta)}/\s(\tau(\beta))$ is normal. We will apply the following result. 

\begin{lem}\label{lem:normalizationisocodim2}
Let $X,Y$ be irreducible affine varieties, with $X$ normal. If $f \colon X \to Y$ is a morphism which restricts to an isomorphism $U \to V$ on open subsets, with complement of codimension at least two in $X$ and $Y$ respectively, then $X$ can be identified with the normalization of $Y$ such that $f$ is the normalization map. 
\end{lem}

\begin{proof}
If $p \colon \widetilde{Y} \to Y$ is the normalization of $Y$ then $f = p \circ \tilde{f}$ factors through $p$. Moreover, since $p$ is finite $\widetilde{V} := p^{-1}(V)$ will have complement of codimension at least two in $\widetilde{Y}$, with $\tilde{f}$ restricting to an isomorphism $U \to \widetilde{V}$. Therefore, replacing $Y$ by $\widetilde{Y}$, we may assume $Y$ is normal too. By the $(S_2)$ property (a version of Hartogs' Lemma), every global function on $U$ uniquely extends to a global function on $X$, so that $X = \Spec \Gamma(U,\mc{O}_U)$ and similarly for $Y$. As a result, the pullback map $f^*$ induces an isomorphism $\Gamma(V,\mc{O}_V) \iso \Gamma(U, \mc{O}_U)$. 
\end{proof}

We construct the morphism 
$$
\prod_{\beta \in \Sigma_{\bc}} \mf{M}_{\bc}(\beta)^{\ell(\tau(\beta))}/\s(\tau(\beta)) \to \mf{M}_{\bc}(\alpha).
$$
First, there is a closed embedding $\prod_{\beta \in \Sigma_{\bc}} \mr{Rep}(\overline{\Qu},\beta)^{\ell(\beta)} \to \mr{Rep}(\overline{\Qu},\alpha)$ given by taking direct sum of representations:
$$
((M_{\beta,i})_{i = 1}^{\ell(\tau(\beta))})_{\beta} \mapsto \bigoplus_{\beta} \bigoplus_{i = 1}^{\ell(\tau(\beta))} M_{\beta,i}^{\oplus \tau(\beta)_i}.
$$
If $G(\tau) = \prod_{\beta \in \Sigma_{\bc}} G(\beta)^{\times \ell(\tau(\beta))}$ is the group acting on the product $\prod_{\beta \in \Sigma_{\bc}} \mr{Rep}(\overline{\Qu},\beta)^{\ell(\beta)}$ then $G(\tau) \subset G(\alpha)$ and the closed embedding is equivariant for $G(\tau)$. Differentiating this inclusion gives an embedding of Lie algebras $\bigoplus_{\beta \in \Sigma_{\bc}} \mf{g}(\beta)^{\oplus \ell(\tau(\beta))} \hookrightarrow \mf{g}(\alpha)$. This fits into a commutative diagram 
$$
\begin{tikzcd}
	\bigoplus_{\beta \in \Sigma_{\bc}} \mr{Rep}(\overline{\Qu},\beta)^{\ell(\beta)} \ar[r,hook] \ar[d,"\prod \mu_{\beta}^{\ell(\beta)}"] & \mr{Rep}(\overline{\Qu},\alpha) \ar[d,"\mu"] \\
	\bigoplus_{\beta \in \Sigma_{\bc}} \mf{g}(\beta)^{\oplus \ell(\tau(\beta))} \ar[r,hook] & \mf{g}(\alpha)
\end{tikzcd}
$$
implying that there is a closed embedding of schemes $(\prod \mu_{\beta}^{\ell(\beta)})^{-1}(\bc) \hookrightarrow \mu^{-1}(\bc)$. Since this is $G(\tau)$-equivariant it induces a map 
$$
(\prod \mu_{\beta}^{\ell(\beta)})^{-1}(\bc) \git \, G(\tau) \to \mu^{-1}(\bc) \git \, G(\alpha). 
$$
Taking the reduced scheme structure on both sides we get a map $\prod_{\beta \in \Sigma_{\bc}} \mf{M}_{\bc}(\beta)^{\ell(\beta)} \to \mf{M}_{\bc}(\alpha)$. The group $\prod_{\beta \in \Sigma_{\bc}}  \s(\tau(\beta))$ acts naturally on the left hand side. It is clear on the level of points that the above map factors through the action of this group. However, to see this algebraically, we note that $\prod_{\beta \in \Sigma_{\bc}}  \s(\tau(\beta)) \hookrightarrow N_{G(\alpha)}(G(\tau))/G(\tau)$ and so we get an induced map 
$$
\left( (\prod \mu_{\beta}^{\ell(\beta)})^{-1}(\bc) \git \, G(\tau) \right) / \prod_{\beta \in \Sigma_{\bc}}  \s(\tau(\beta)) \to \mu^{-1}(\bc) \git \, G(\alpha).
$$
In other words, this is a morphism $\prod_{\beta \in \Sigma_{\bc}} \mf{M}_{\bc}(\beta)^{\ell(\beta)}/\s(\tau(\beta)) \to \mf{M}_{\bc}(\alpha)$. As noted previously, the left hand side is normal. Inside $\mf{M}_{\bc}(\beta)^{\ell(\beta)}$ we write $\mf{M}_{\bc}(\beta)^{\ell(\beta),\circ}$ for the (dense) open set consisting of pairwise distinct, and simple, representations. The set of points in $\mf{M}_{\bc}(\beta)$ parameterizing simple representations is the open stratum and hence is even dimensional with complement of codimension at least two. This implies that the complement to $\mf{M}_{\bc}(\beta)^{\ell(\beta),\circ}$ in $\mf{M}_{\bc}(\beta)^{\ell(\beta)}$ has codimension at least two as well. Notice that $\s(\tau(\beta))$ acts freely on this open set. The image of $\prod_{\beta \in \Sigma_{\bc}} \mf{M}_{\bc}(\beta)^{\ell(\beta),\circ}$ equals the symplectic leaf $\LG_{\tau}$. Moreover, 
\[
\prod_{\beta \in \Sigma_{\bc}} \mf{M}_{\bc}(\beta)^{\ell(\beta),\circ} /\s(\tau(\beta)) \to \LG_{\tau}
\]
is a bjiection on closed points. Since both domain and image are smooth, this is an isomorphism. Taking affine closures, it induces a morphism 
$$
\prod_{\beta \in \Sigma_{\bc}} \mf{M}_{\bc}(\beta)^{\ell(\beta)}/\s(\tau(\beta))\to \overline{\LG}_{\tau}.
$$
This must necessarily agree with the previous map. 

\begin{lem}\label{lem:smoothlocusleafclosure} 
	The smooth locus of $\overline{\LG}_{\tau}$ equals $\LG_{\tau}$ and the singular locus has codimension at least two. 
\end{lem}

\begin{proof}
	First we note that $\mf{M}_{\bc}(\alpha)$ has symplectic singularities. Therefore it is holonomic in the sense of Kaledin \cite{Kaledinsympsingularities}. Then \cite[Proposition 3.1]{Kaledinsympsingularities} implies that the singular locus of  $\overline{\LG}_{\tau}$ has codimension at least two. Moreover, \cite[Lemma~1.4]{Kaledinsympsingularities} says that the Poisson structure on the smooth locus of  $\overline{\LG}_{\tau}$ is non-degenerate. Since $\overline{\LG}_{\tau}$ is a union of leaves, with unique leaf $\LG_{\tau}$ of top dimension, we deduce that $\LG_{\tau}$ is the smooth locus of $\overline{\LG}_{\tau}$ since the Poisson structure is degenerate along all other leaves.
\end{proof}

Then the main result follows from the following.  

\begin{lem}
\label{lem:normalization-leaf-cl}
	The map $\prod_{\beta \in \Sigma_{\bc}} \mf{M}_{\bc}(\beta)^{\ell(\beta)}/\s(\tau(\beta)) \to \overline{\LG}_{\tau}$ is the normalization of $\overline{\LG}_{\tau}$.
\end{lem}

\begin{proof}
	Let $A = \mc{O}[\overline{\LG}_{\tau}]$ and $B$ the ring of functions on $\prod_{\beta \in \Sigma_{\bc}} \mf{M}_{\bc}(\beta)^{\ell(\beta)}/\s(\tau(\beta))$. Then $\phi$ corresponds to an inclusion $A \hookrightarrow B$. Lemma~\ref{lem:smoothlocusleafclosure} guarantees that the hypothesis of Lemma~\ref{lem:normalizationisocodim2} hold. Therefore the statement of the lemma follows from the latter result.
\end{proof} 

\begin{example}\label{ex:symsn}
	Consider the case where $\mathsf{G}$ is the graph with one vertex and one loop. Then $\mf{M}_0(n) \cong \C^{2n} / \s_n$ and $\Sigma_0 = \{ e_1 \}$. This means that the strata are labeled by $\tau \colon \{ e_1 \} \to \mc{P}$ such that $|\tau(e_1)| = n$. This is precisely the set of partitions of $n$. Therefore we think of $\tau = (1^{\tau_1}, 2^{\tau_2},\ds )$ as a partition of $n$. Then $\mf{M}_0(e_1) = \C^2$ so for a given partition, the normalization of $\overline{\LG}_{\tau}$ is the map 
	$$
	\C^{2 \ell(\tau)} / \prod_{i \ge 1} \s_{\tau_i} \to \overline{\LG}_{\tau}.
	$$
	Notice, in particular, that when each part occurs in $\tau$ at most once (this means $\tau_i \le 1$ for all $i$) then the normalization of $\overline{\LG}_{\tau}$ is just $\C^{2 \ell(\tau)}$. The leaf closure $\overline{\LG}_{\tau}$ is normal if and only if $\tau = (m^k)$ has a rectangular Young diagram. This can be deduced from \cite[Lemma~2.2]{RichardsonNormality}, just as in \cite[Proposition~8.1.2]{RichardsonNormality} which considers instead the closure of strata in $\C^n / \s_n$.
\end{example}

\subsection{Normality of a leaf closure}

Having identified the normalization of a leaf closure, it is natural to ask if the leaf closure itself is normal. This is addressed in detail in \cite{BellSchedlerQuiverNormal}. In the case of Calogero--Moser varieties, we will see below that the answer is no in general. However, we give one situation (which often occurs in the context of Calogero--Moser varieties) where the leaf closure is normal. 

\begin{prop}\label{prop:leafclosurenormal}
	Assume that $\tau = (n_1,\beta;n_2, \alpha^{(2)}; \ds; n_k, \alpha^{(k)})$, where every $\alpha^{(i)}$ for $i \ge 2$ is real. Then $\overline{\LG}_{\tau} \cong \mf{M}_{\bc}(\beta)$ is normal. 
\end{prop}

\begin{proof}
	Note that since $\alpha^{(i)}$ is real, $\mf{M}_{\bc}(\alpha^{(i)})$ is a single point. Therefore, Theorem~\ref{thm:quiverleafclosure} says that $\iota \colon \mf{M}_{\bc}(\beta) \to \overline{\LG}_{\tau}$ is the normalization map. The morphism $\iota^* \colon \C[\overline{\LG}_{\tau}] \to \C[\mf{M}_{\bc}(\beta)]$ is necessarily injective. We must show that it is surjective. By \cite[Theorem~1]{Semisimplequiver}, both the ring $\C[\mf{M}_{\bc}(\beta)]$ and the ring $\C[\mf{M}_{\bc}(\alpha)]$ are generated by traces of oriented cycles $\omega$ in $\mathsf{G}$. If $\omega$ is an oriented cycle starting and ending at vertex $j$ and $X \in \mu^{-1}(\bc)$ then we write $\Tr \omega (X)$ for the trace of the endomorphism $\omega(X)$ of $V_j$ ($\Tr \omega (X)$ does not depend on the choice of $j$). Let $M_i$ denote the simple $\Pi^{\bc}$-module of dimension $\alpha^{(i)}$ and choose some $X_i \in \mu_{\alpha^{(i)}}^{-1}(\bc)$ above $M_i$. Then, as described above we have a commutative diagram 
	$$
	\begin{tikzcd}
		\mu^{-1}_{\beta}(\bc) \ar[r,"\iota_0"] \ar[d] & \mu^{-1}_{\alpha}(\bc) \ar[d] \\
		\mf{M}_{\bc}(\beta) \ar[r,"\iota"] & \mf{M}_{\bc}(\alpha)
	\end{tikzcd}
	$$
	where $\iota_0(X) = (X \o \Id_{n_1} + X_2 \o \Id_{n_2} + \cdots + X_k \o \Id_{n_k})$. If $[M] \in \mf{M}_{\bc}(\beta)$ and $X \in \mu^{-1}_{\beta}(\bc)$ any (not necessarily semi-simple) lift of $[M]$, then 
	\begin{align*}
		(\iota^* \Tr \omega)([M]) & = (\iota_0^* \Tr \omega)(X) \\
		& = \Tr \omega(X \o \Id_{n_1} +X_2 \o \Id_{n_2} + \cdots + X_k \o \Id_{n_k}) \\
		& = n_1 \Tr \omega(X) + n_2 \Tr \omega(X_2) + \cdots + n_k \Tr \omega(X_k) \\
		& = n_1 (\Tr \omega)([M]) + n_2 \Tr \omega(X_2) + \cdots + n_k \Tr \omega(X_k).
	\end{align*}
	Hence, $\Tr \omega = \frac{1}{n_1} (\iota^* \Tr \omega - (n_2 \Tr \omega(X_2) + \cdots + n_k \Tr \omega(X_k)))$ in $\C[\mf{M}_{\bc}(\beta)]$ and the result follows.
\end{proof}

In terms of framed quiver varieties $\mf{M}_{\bc}(\bw,\alpha)$,  Proposition~\ref{prop:leafclosurenormal} says that if 
\[
\tau = (\beta^{(0)},(n_1, \beta^{(1)}; \ds; n_k, \beta^{(k)}))
\]
is a decomposition where every $\beta^{(i)} \in \Sigma_{\bc}$, for $i = 1,\ds, k$, is a real root then $\overline{\LG}_{\tau} \cong \mf{M}_{\bc}(\bw,\beta^{(0)})$ is normal.

\section{Calogero--Moser varieties}\label{sec:symplc-refl-alg}

The canonical symplectic form on $\C^2$ induces a symplectic form on $V := (\C^2)^n$. Recall that $\Gamma \subset \SL_2(\C)$ is a finite group. Let $\Gamma_n = \Gamma \wr \s_n = \Gamma^n \rtimes \s_n \subset \Sp(V)$. This is a symplectic reflection group; that is, $\Gamma_n$ is generated by its symplectic reflections, which are the elements $s \in \Gamma_n$ with $\mr{rk}(1 - s)|_V = 2$. To each symplectic reflection $s$, we associate the degenerate $2$-form $\omega_s$ on $V$ which equals $\omega$ when restricted to $\mr{Im} (1 - s)$ and is zero on $\Ker (1 - s)$. If $\gamma \in \Gamma \setminus \{ 1 \}$ then we write $\gamma_i$ for the element of $\Gamma^n$ which is $\gamma$ in the $i$th position and one elsewhere. The transpositions in $\s_n$ are $s_{i,j}$. Assuming $n > 1$, the symplectic reflections in $\Gamma_n$ are $\gamma_i$, for $1 \le i \le n$ and $\gamma \in \Gamma \setminus \{ 1 \}$, together with $s_{i,j} \gamma_i \gamma_j^{-1}$, for $1 \le i \neq j \le n$ and $\gamma \in \Gamma$. When $n = 1$, $\Gamma_n = \Gamma$ and the symplectic reflections in $\Gamma_n$ are just elements of $\Gamma \setminus \{ 1 \}$. When $n = 0$, we set $V = \{ 0 \}$ and $\s_n = \Gamma_n = \{ 1 \}$.   

Given a conjugate invariant function $c \colon \mc{S} \to \C$, we define the \textit{symplectic reflection algebra} $\Hb_{c}(\Gamma_n)$ (at $t = 0$) to be the quotient of $T V \rtimes \Gamma_n$ by the relations
\[
u \o v - v \o u = \sum_{s \in \mc{S}} c(s) \omega_{s}(u,v) s, \quad \forall \, u,v \in V. 
\]
The centre of the symplectic reflection algebra $\Hb_{c}(\Gamma_n)$ is denoted $\Zb_{c}(\Gamma_n)$ and the Calegero--Moser variety is $\ZCB_{c}(\Gamma_n) := \Spec \Zb_{c}(\Gamma_n)$.
The ring $\Hb_{c}(\Gamma_n)$ is prime and thus $\Zb_{c}(\Gamma_n)$ is a domain. Moreover, it is known to be integrally closed \cite[Lemma~3.5]{EG} and hence $\ZCB_{c}(\Gamma_n)$ is a normal variety. Since $\Zb_{c}(\Gamma_n)$ has a quantization given by the spherical subalgebra of $\Hb_{c}(\Gamma_n)$ at $t = 1$, it has a Poisson structure that is generically non-degenerate; see \cite{EG}. When $n = 0$, we define $\ZCB_{c}(\Gamma_n) = \{ \mr{pt} \}$. 

Let us agree that when we write $\Gamma=\ZM/\ell\ZM$, we mean that $\Gamma$ is the subgroup of $\SL_2(\C)$ whose elements are of the form $\left(\begin{matrix}
        \zeta & 0\\
        0 &\zeta^{-1}\\
    \end{matrix}\right)$, where $\zeta$ is an $\ell$th root of unity.

\begin{rem}
\label{rem:C*-on-CM}
    Assume that $\Gamma=\ZM/\ell\ZM$. 
    The inclusion $\CM^\times\subset \SL_2(\C)$ given by 
    \[
    \xi\mapsto
    \left(\begin{matrix}
        \xi & 0\\
        0 &\xi^{-1}\\
    \end{matrix}\right)
    \]
    yields a $\CM^\times$-action on $V$. This action 
    induces a $\CM^\times$-action on $\Hb_{c}(\Gamma_n)$ by automorphisms (this action is trivial on $\CM\Gamma_n$).
    This induces a $\CM^\times$-action on $\ZCB_{c}(\Gamma_n)$.
    
\end{rem}

\subsection{The McKay correspondence}\label{rootsystemnotation}

Let $\Gamma\subset \SL(2,\C)$ be a finite group and $\mathsf{G}$ the simply laced affine Dynkin graph associated to $\Gamma$ via the McKay correspondence. Thus, $I$ is identified with the set of isomorphism classes of irreducible $\Gamma$-modules in such a way that $0 \in I$ corresponds to the trivial representation. Then $\Delta = \{ e_0, e_1, \ds, e_r \}$ is a set of simple roots for the affine root system $R \subset Q$. We let $\delta = \sum_{i \in I} \delta_i e_i$ denote the positive root that spans the radical of $( - , -)$ on $Q$. It satisfies $\delta_0 = 1$. For the remainder of the article, $W^{\aff} := W_{\mathsf{G}}$ is the affine Weyl group associated to the affine Dynkin diagram defined by $\Gamma$. Let $R_{\bc}$ denote the set of roots $\alpha$ in $R$ such that $\bc(\alpha) = 0$. 

If $\Gamma$ is nontrivial, then  roots $e_1, \ds, e_r$ are a set of simple roots for a finite root system $\Phi \subset R$. This root system is irreducible and $\theta \in \Phi^+$ denotes the longest positive root. Then $\delta = e_0 + \theta$. We denote by $\boldsymbol{\Lambda} \subset Q$ the root lattice for $\Phi$.    

If $\Gamma$ is trivial then we define $\delta := e_0$ and $\boldsymbol{\Lambda}=0$.

Let $L(\Lambda_0)$ be the integrable highest weight module with highest weight $\Lambda_0$ (the basic representation) for the affine Lie algebra $\mf{g}$ with root system $R$ (when $\Gamma$ is trivial, we take the infinite-dimensional Heisenberg Lie algebra as in \cite{NakHeisenberg}). We recall the properties of $L(\Lambda_0)$. First, recall that a \textit{weight} of $L(\Lambda_0)$ is an element $v \in \mathfrak{h}^*$ of the dual of the Cartan subalgebra such that $L(\Lambda_0)_v \neq 0$, when $\Gamma$ is nontrivial. For $\Gamma$ trivial, we replace $\mathfrak{h}^*$ with the two-dimensional vector space spanned by $\delta=e_0$ and $\Lambda_0$. Note that $P,Q \subset \mathfrak{h}^*$ are free abelian groups of rank $r+1$, with sum $P+Q$ a full lattice in  $\mathfrak{h}^*$, which has dimension $r+2$. A weight $v$ of $L(\Lambda_0)$ is called \textit{maximal} if $v + \delta$ is not a weight. Finally, for $v = \sum_i v_i \Lambda_i \in P$, then for every $\lambda \in \C$, we say that $v+\lambda \delta$ is \textit{dominant}  if $v_i \geq 0$ for all $i$ (otherwise $v+\lambda \delta$ is not dominant). 
\begin{lem}\label{lem:basicrepCarter}
	\begin{enumerate}
 \item[(i)] The weights of $L(\Lambda_0)$ are $\Lambda_0 - m \delta - \frac{1}{2}(\nu,\nu) \delta + \nu$ for some (unique) $m \ge 0, \nu \in \Lambda$. 
		\item[(ii)] The weights $\Lambda_0 - \frac{1}{2}(\nu,\nu) \delta + \nu$, $\nu \in \Lambda$, are precisely the maximal weights of $L(\Lambda_0)$.
		\item[(iii)] $\Lambda_0$ is the unique dominant maximal weight. 
		\item[(iv)] The maximal weights of $L(\Lambda_0)$ form a single $W^{\aff}$-orbit. 
		\item[(v)] There exists $w \in W^{\aff}$ such that 
		$$
		w\left(\Lambda_0 - m \delta - \frac{1}{2}(\nu,\nu) \delta + \nu\right) = \Lambda_0 - m \delta. 
		$$
	\end{enumerate}
\end{lem}

\begin{proof}
	These all follow from \cite[Theorem~20.23]{CarterBook}. Specifically, part (i) is \cite[Theorem~20.23(c)]{CarterBook}, (ii) is \cite[Theorem~20.23(b)]{CarterBook} and part (iii) is \cite[Theorem~20.23(a)]{CarterBook}. Part (iv) then follows from \cite[Corollary~20.15]{CarterBook} and \cite[Theorem~20.23(a)]{CarterBook}. This implies that there exists $w \in W^{\aff}$ such that $w(\Lambda_0 - \frac{1}{2}(\nu,\nu) \delta + \nu) = \Lambda_0$. Since $W^{\aff}$ acts trivially on $\delta$, part (v) follows from part (iv). 
\end{proof}

\subsection{Calogero--Moser varieties as quiver varieties}
\label{CM-as-quiv}

In the case of $\Gamma_n$ with $n > 1$, there are two types of conjugacy class of symplectic reflections. Namely, the set 
\[
\mc{S}_0 = \{ s_{i,j} \gamma_i \gamma_j^{-1} \, | \, \textrm{$1 \le i \neq j \le n$ and $\gamma \in \Gamma$} \}
\]
forms a single conjugacy class and each $\{ \gamma_i \, | \, \gamma \in \mc{C}, 1 \le i \le n \}$ is another conjugacy class, as $\mc{C}$ runs over all non-trivial conjugacy classes in $\Gamma$.  Therefore, $c = (c_1, \underline{c})$. Here $\underline{c} \colon \Gamma \smallsetminus \{ 1 \} \to \C$ is a conjugate invariant function and $c_1$ is the value of $c$ on $\mc{S}_0$. From $\underline{c}$ we get the element $z(\underline{c}) \in Z (\C \Gamma)$ by setting $z(\underline{c}) := \sum_{\gamma \in \Gamma \setminus \{ 1 \} } \underline{c}(\gamma)\gamma$. When $n = 1$, the conjugacy classes of symplectic reflections are just the $\mc{C}$ and $c = (\underline{c})$; there is no $c_1$. 

From $c$, we define a parameter $\bc \in \Hom_{\Z}(Q,\C)$ by  
\begin{equation}\label{eq:boldc}
   \bc (e_i) := -\frac{1}{2} c_1 |\Gamma| \Tr_{S_i} \mathbf{e}_{\Gamma} + \Tr_{S_i} z(\underline{c}), \quad \forall \, 0 \le i \le r.  
\end{equation}
Here $S_i$ is the irreducible representation (of dimension $\delta_i$) of $\Gamma$ corresponding to vertex $i$ and $\mathbf{e}_{\Gamma}$ is the trivial idempotent in $\C \Gamma$. When $n = 1$, one should take an arbitrary $c_1$ in \eqref{eq:boldc}. The variety $\XC_\bc(n\delta)  = \mf{M}_{{\bc}}(\Lambda_0, n \delta)$ does not depend on $c_1$ when $n=1$. The following is \cite[Theorem~1.4]{MarsdenWeinsteinStratification}, generalizing \cite[Theorem~1.13]{EG}.

\begin{thm}\label{thm:CMquiveriso}
    For all $n \ge 0$, there is an isomorphism $\ZCB_{c}(\Gamma_n) \iso \XC_{{\bc}}(n \delta)$ of Possion varieties. 
\end{thm}

Note that, with our conventions, the theorem still makes sense when $n = 0$ since $\XC_{{\bc}}(0)$ is a point. 

The number $\bc (\delta)$ is usually called the \textit{level}. Since we are at $t = 0$, there are isomorphisms $\ZCB_{c}(\Gamma_n) \iso \ZCB_{q c}(\Gamma_n)$ and $\XC_{{\bc}}(n \delta) \iso \XC_{q{\bc}}(n \delta)$ for any $q \in \Cs$, compatible with the identification of Theorem~\ref{thm:CMquiveriso}. Therefore, we only need to consider the cases $\bc(\delta) = 1$ and $\bc(\delta) = 0$.

\subsection{Parabolic subgroups}\label{sec:paraboliclabel}

A subgroup $\mr{P} \subset \Gamma_n$ is said to be a \textit{parabolic subgroup} if it is the stabiliser of some vector $v \in V$. Parabolic subgroups play an important role in the classification of symplectic leaves in the Calogero--Moser variety. First, we explicitly describe the parabolic subgroups of $\Gamma_n$. 

\begin{lem}\label{lem:parabolicwreathproduct}
	The parabolic subgroups of $\Gamma_n$ are, up to conjugacy, of the form $\Gamma_m \times \s_{\lambda}$, where $\lambda$ is a partition with $n = |\lambda| + m$. The normalizer of this parabolic can be described as
	$$
	N_{\Gamma_n}(\Gamma_m \times \s_{\lambda}) = \Gamma_m \times ((\s_{\lambda} \times \Gamma^{\ell(\lambda)}) \rtimes \s(\lambda)) = \Gamma_m \times N_{\Gamma_{n-m}}(\s_{\lambda}),
	$$ 
	and hence $N_{\Gamma_n}(\Gamma_m \times \s_{\lambda}) / (\Gamma_m \times \s_{\lambda}) \cong \prod_{i \ge 1} \Gamma_{n_i}$. 
\end{lem}

\begin{proof}
The fact that the parabolic subgroups are all conjugate to a subgroup of the form $\Gamma_m \times \s_{\lambda}$ is standard; see e.g. the proof of \cite[Proposition~3.4]{BellSchmitThiel}. If we let $U \subset (\C^2)^n$ denote the set of points of the form $(0,\ds,0,v_1,\ds, v_1,v_2,\ds,v_2,\ds, v_{\ell(\lambda)})$ where $v_i \in \C^2$ with $(\Gamma \cdot v_i) \cap (\Gamma \cdot v_j) = \emptyset$ for $i \neq j$ and there are $\lambda_1$ copies of $v_1$, $\lambda_2$ copies of $v_2$ etc. then $\Gamma_m \times \s_{\lambda}$ is the stabilizer of any $u \in U$ and $N_{\Gamma_n}(\Gamma_m \times \s_{\lambda})$ is the set of elements of $\Gamma_n$ mapping $U$ into itself. Then it is straight-forward to check that 
\[
N_{\Gamma_n}(\Gamma_m \times \s_{\lambda}) = \Gamma_m \times ((\s_{\lambda} \times \Gamma^{\ell(\lambda)}) \rtimes \s(\lambda)).
\]
The final isomorphism follows from the fact that $\Gamma^{\ell(\lambda)} \rtimes \s(\lambda) = \prod_{i \ge 1} \Gamma_{n_i}$.  
\end{proof}

Let $\mr{P}$ be a parabolic subgroup of $\Gamma_n$. By definition, it is normal in its normalizer $N(\mr{P})$. Let $\Xi := N_G(\mr{P})/\mr{P}$ be the quotient. The conjugacy class of $\mr{P}$ in $\Gamma_n$ is denoted $(\mr{P})$. The algebra $\Hb_{c}(\Gamma_n)$ has a canonical filtration given by placing $\Gamma_n$ in degree zero and $V$ in degree one. Then $\Zb_{c}(\Gamma_n)$ inherits a filtration by restriction and $\gr(\Zb_{c}(\Gamma_n)) \cong \C[V^*]^{\Gamma_n}$ by \cite[Theorem~3.3]{EG}. If $\mf{p} \subset \Zb_{c}(\Gamma_n)$ is the prime ideal defining the closure of a symplectic leaf of $\ZCB_{c}(\Gamma_n)$, then \cite[Theorem~2.8]{MartinoAssociated} says that $\gr(\mf{p})$ is a prime ideal defining the closure of a symplectic leaf in $V^*/\Gamma_n$. Since the leaves of $V^*/\Gamma_n$ are in bijection with conjugacy classes of parabolic subgroups of $\Gamma_n$, the leaves in $\ZCB_{c}(\Gamma_n)$ can also be labeled by conjugacy classes of parabolic subgroups of $\Gamma_n$. However, the same conjugacy class can label several different leaves. 

This labeling of symplectic leaves by conjugacy classes of parabolic subgroups is important because Losev has shown that there is a notion of induction of leaves whose construction depends on this labeling. Let $\mr{PSpec}_{(\mr{P})} \ZCB_{c}(\Gamma_n)$ denote the set of all leaves in $\ZCB_{c}(\Gamma_n)$ that are labeled by the conjugacy class $(\mr{P})$. Here $\mr{PSpec}$ denotes the set of Poisson prime ideals. We fix a representative $\mr{P}$ in $(\mr{P})$. There is a unique zero-dimensional leaf $\{ 0 \}$ in $V^*/\Gamma_n$; it is labeled by $(\Gamma_n)$. Next, we consider the algebra $\Hb_{c|\mr{P}}(\mr{P},(V^{\mr{P}} )^{\perp})$, the symplectic reflection algebra defined by the subgroup $\mr{P}$, the restriction $c|\mr{P}$, and the subspace $(V^{\mr{P}} )^{\perp} \subset V$, where the orthogonal is with respect to the symplectic form on $V$. The group $\Xi$ acts on the Calogero--Moser variety $\ZCB_{c|\mr{P}}(\mr{P},(V^{\mr{P}} )^{\perp})$, permuting the Poisson prime ideals. The set $\mr{PSpec}_{(\mr{P})} \ZCB_{c|\mr{P}}(\mr{P},(V^{\mr{P}} )^{\perp})$ of zero-dimensional leaves in $\mr{PSpec} \, \ZCB_{c|\mr{P}}(\mr{P},(V^{\mr{P}} )^{\perp})$ is stable under this action. 
By \cite[Theorem~1.3.2(4)]{LosevSRAComplete}:

\begin{thm}\label{thm:losevinductionleaves}
There exists a bijection 
\[
(\mr{PSpec}_{(\mr{P})}\ZCB_{\bc|\mr{P}} (\mr{P},(V^{\mr{P}} )^{\perp})) / {\Xi} \stackrel{1:1}{\longleftrightarrow} \mr{PSpec}_{(\mr{P})}\ZCB_{c}(\Gamma_n).
\]
\end{thm}

We note an immediate corollary. 

\begin{cor}\label{cor:level1leaflabel}
    If $\bc(\delta) \neq 0$ then the leaf $\LG \subset \ZCB_{c}(\Gamma_n)$ is labeled by $(\Gamma_m)$, where $2(n - m) = \dim \LG$. 
\end{cor}

\begin{proof}
    By Theorem~\ref{thm:losevinductionleaves}, the leaf $\LG$ is labeled by some conjugacy class $(\mr{P})$ with $\dim \LG = \dim V^{\mr{P}}$. Lemma~\ref{lem:parabolicwreathproduct} says that $\mr{P} =\Gamma_m \times \s_{\lambda}$ with $n = m + |\lambda|$. Moreover, this leaf corresponds to a $\Xi$-orbit of zero-dimensional leaves in 
    \[
    \ZCB_{\bc|\mr{P}} (\mr{P},(V^{\mr{P}} )^{\perp}) = \ZCB_{\bc|\Gamma_m} (\Gamma_m,\C^{2m}) \times \prod_{i} \ZCB_{c_1} (\s_{\lambda_i},\h_i),
    \]
    where $\h_i$ is the reflection representation for $\s_{\lambda_i}$. We note that 
    \begin{align*}
         \bc (\delta) & = \sum_{i= 0}^r (\dim S_i) \left( -\frac{1}{2} c_1 |\Gamma| \Tr_{S_i} \mathbf{e}_{\Gamma} + \Tr_{S_i} z(\underline{c}) \right)   \\
         & = -\frac{1}{2} c_1 |\Gamma| + \Tr_{\C \Gamma} z(\underline{c}) = -\frac{1}{2} c_1 |\Gamma|
    \end{align*}
    where we have used the fact that $z(\underline{c})$ is a sum over the non-trivial conjugacy classes in $\Gamma$. Thus, $\bc (\delta) \neq 0$ implies that $c_1 \neq 0$. In this case, each $\ZCB_{c_1} (\s_{\lambda_i},\h_i)$ is smooth by \cite[Theorem~1.24]{EG}. In particular, it has no zero-dimensional leaves unless $\h_i = \{ 0 \}$ i.e. $\lambda_i = 1$. Thus, $\lambda = (1^{n-m})$ and $\s_{\lambda} = \{ 1 \}$. 
\end{proof}

\section{Non-zero level}\label{sec:nonzero}

Throughout this section, we assume that $\bc(\delta) \neq 0$. Without loss of generality, $\bc(\delta) = -1$. 

\subsection{} Recall that $R_{\bc}^+$ is the set of roots in $R^+$ that dot to zero with $\bc$. We say that $\alpha \in R_{\bc}^+$ is \textit{minimal} if $\alpha$ cannot be written as a sum of two vectors in $R_{\bc}^+$. Let $\Delta(\bc) \subset R_{\bc}^+$ be the set of minimal vectors. 

\begin{lem}
	The set $\Delta(\bc)$ is a set of simple roots for the root system $R_{\bc}$ such that the corresponding positive roots are precisely  $R_{\bc}^+$. 
\end{lem}

\begin{proof}
	We note that $R_{\bc} = R_{\bc}^+ \cup - R_{\bc}^+$ since $R = R^+ \cup - R^+$. Moreover, every root in $R_{\bc}^+$ can be written as a positive (integer) sum of minimal roots. The lemma follows. 
\end{proof}

We say that a root subsystem $S \subset R$ is a \textit{parabolic root subsystem} if there exists $\Delta_S \subset \Delta$ such that $\Delta_S$ is a set of simple roots for $S$. 

\begin{lem}\label{lem:wparabolicrootsystem}
	There exists $w \in W^{\aff}$ such that $w(R_{\bc}) \cap \Delta$ is a set of simple roots for $w(R_{\bc})$. In particular, $w(R_{\bc})$ is a (proper) parabolic root subsystem of $R$. 
\end{lem}

\begin{proof}
We write $\bc = \bc_{\mr{re}} + \sqrt{-1} \bc_{\mr{im}}$, where $\bc_{\mr{re}}, \bc_{\mr{im}} \in \R^{r+1}$. Then $\bc_{\mr{re}}(\delta) = -1$ and $\bc_{\mr{im}}(\delta) = 0$. Since $\bc_{\mr{re}} (\delta) = -1$, \cite[Proposition~3.12(c)]{KacBook} implies that $\bc_{\mr{re}}$ is in the Tits cone. Therefore, there exists $w \in  W^{\aff}$ such that the real part $\bc_{\mr{re}}'$ of $\bc' := w^*(\bc)$ belongs to (dual) fundamental domain. 

Then $w(R_{\bc}) = R_{\bc'}$. Since $\bc'(\delta) = \bc(\delta) = -1$, the set $R_{\bc'}$ consists of real roots. The subgroup $W'$ of $W^{\aff}$ generated by these real roots is contained in the stabilizer subgroup $W^{\aff}_{\bc'}$ of $\bc'$. However, by  \cite[Proposition~3.12(a)]{KacBook}, $W^{\aff}_{\bc'}$ is generated by the $s_i$, $i = 0, \ds, r$ such that $s_i^*(\bc')= \bc'$. That is, by the reflections along the roots in $\Delta \cap R_{\bc'}$. Hence $W' = W^{\aff}_{\bc'}$ and $\Delta \cap R_{\bc'}$ is a set of simple roots for $R_{\bc'}$. In particular, $R_{\bc'}$ is a parabolic subsystem of $R$. Since it does not contain any imaginary roots, it must be a proper subsystem. 
\end{proof}

This means that, up to the action of $W^{\aff}$, the root system $R_{\bc}$ is given by deleting a certain number of nodes in the affine Dynkin diagram. Since $R_{\bc}$ is a finite root system, the element $w$ in Lemma~\ref{lem:wparabolicrootsystem} can be chosen so that $w(\Delta(\bc)) = w(R_{\bc}) \cap \Delta$. 

Recall that we have defined in \eqref{eq:twistedaction} a twisted action of $W^{\aff}$ on $Q$ by 
\[
w \star \alpha = w(\alpha - \Lambda_0) + \Lambda_0.
\]
Here we have taken $\bw = \Lambda_0$ and omitted it from the notation. For $\nu \in \boldsymbol{\Lambda}$ and $m$ a non-negative integer, define
$$
\gamma(m,\nu) := m \delta + (1/2)(\nu,\nu) \delta - \nu. 
$$

\begin{lem}\label{lem:basicrepCarternonLie}
Let $\alpha \in Q^+$. Then $\Lambda_0 - \alpha$ is a weight of $L(\Lambda_0)$ if and only if there exists $\nu \in \boldsymbol{\Lambda}$ and $m \ge 0$ such that 
\[
\alpha = \gamma(m,\nu).
\]
Moreover, if $\Lambda_0 - \alpha$ is a weight of $L(\Lambda_0)$ then there exists $w \in W^{\aff}$ such that $w \star \alpha = m \delta$. 
\end{lem}

\begin{proof}
This is a reformulation of Lemma~\ref{lem:basicrepCarter}(i) and (v). 
\end{proof}

Recall from Sections~\ref{root-sys} and \ref{sec:CBtrick} that we have defined the sets $\Sigma_{\bc} \subset R^+_{\bc}$ and $E_{\bc} \subset Q^+$ (for $\bw = \Lambda_0$). 

\begin{lem}\label{lem:thingsinSigmac}
    Assume $\bc(\delta) = -1$. 
    \begin{enumerate}
        \item[(i)] $\Sigma_{\bc} = \Delta(\bc)$; and 
        \item[(ii)] $n \delta \in E_{\bc}$. 
    \end{enumerate}
\end{lem}

\begin{proof}
 Part (i). If $\beta \in \Sigma_{\bc}$ then $\beta$ is a positive real root since $\bc(\delta) = -1$. In particular, $p(\beta) = 0$. Therefore, being in $\Sigma_{\bc}$ means that it does not admit any proper decomposition into a sum of positive roots in $R_{\bc}^+$. But this is precisely the definition of being in $\Delta(\bc)$. 

 Part (ii) is \cite[Proposition~4.2(i)]{BellamyCrawQuotient} since $n \delta \in E_{\bc}$ precisely when $e_{\infty} + n \delta \in \Sigma_{\bc_{n \delta}}$. 
\end{proof}

\subsection{Another presentation of the affine Weyl group}\label{sec:anotherpresentation}

Recall that the affine Weyl group has another presentation. We have $W^\aff=W\ltimes \boldsymbol{\Lambda}$, where $W$ is the non-affine Weyl group.	For each $\beta\in \boldsymbol{\Lambda}$, denote by $t_\beta$ the image of $\beta$ in $W^\aff$. Each element of $W^\aff$ can be written in a unique way in the form $w\cdot t_\beta$, where $w\in W$ and $\beta\in \boldsymbol{\Lambda}$. We can also extend the notation $t_\alpha$ to $\alpha\in Q$ by setting $t_\alpha:=t_{\pi(\alpha)}$ for each $\alpha\in Q$, where $\pi$ is the following map
	$$
	\pi\colon Q\to Q/\ZM\delta\cong \boldsymbol{\Lambda}.
	$$
	
	\begin{lem}	\label{lem:affWeyl_translation-roots}
		Assume $\alpha,\beta\in Q$. Then we have $t_\beta \star \alpha \equiv \alpha-\beta  \smod  \ZM\delta$.
	\end{lem}
	
	\begin{proof}
		This statement is a special case of \cite[(6.5.2)]{KacBook}.
	\end{proof}
	
Each vector $\alpha \in Q$ defines a linear functional $\overline{\alpha} \in \Hom_{\Z}(Q,\C)$ by $\overline{\alpha}(v) := (\alpha,v)$. Note that the kernel of the map $Q \to \Hom_{\Z}(Q,\C)$ is spanned by $\delta$.	The kernel of this map is $\ZM\delta$. 
	\begin{lem}
		\label{lem:affWeyl_translation-theta}
		For each $\beta\in Q$ and $\bc\in \Hom_{\Z}(Q,\C)$, we have $t_\beta^*(\bc)=\bc+ \bc(\delta)\overline\beta$.
	\end{lem}
	\begin{proof}
  The statement follows from \cite[(6.5.2)]{KacBook}. 
	\end{proof}

\subsection{The symplectic leaves}\label{sec:sympleaf1}

We define on $Q$ a quadratic function 
\begin{equation}\label{eq:qquadraticLambda0}
    \qu(\beta) := \beta_0 + (1/2)(\beta,\beta).
\end{equation}

\begin{remark}
    If $\Gamma \neq \{ 1 \}$ then the affine Dynkin graph $\mathsf{G}$ contains no loops. As in Remark~\ref{rem:noloopaffineLie}, the bilinear form $( - , - )$ extends to $P$ and satisfies $( \Lambda_0,\Lambda_0) = 0$; see \cite[(6.2.2)]{KacBook}. Then the quadratic function $\qu$ can be rewritten as $\qu(\beta) = (1/2)(\Lambda_0 + \beta, \Lambda_0 + \beta)$.   
\end{remark}

Write $\SS{\bc} := \N \Sigma_{\bc}$, which by Lemma~\ref{lem:thingsinSigmac} equals $\N \Delta(\bc)$. If $s = |\Delta(\bc)|$ then Lemma~\ref{lem:wparabolicrootsystem} implies that $\SS{\bc}  \cong \N^s$. Define a partial ordering $\succ$ on $\SS{\bc}$ by $\eta \succ \beta$ if $\eta \neq \beta$ and $\eta - \beta \in \SS{\bc}$. Let $\Xi(\bc)$ denote the set of all vectors $\beta \in \SS{\bc}$ such that $\eta \in \SS{\bc}$ and $\eta \succ \beta$ implies that $\qu(\eta) > \qu(\beta)$. In other words, if $\beta \in \SS{\bc}$ then $\beta$ belongs to $\Xi(\bc)$ if and only if it is maximal, under the partial ordering $\succ$, in the set $\{ \eta \in \SS{\bc} \, | \, \qu(\eta) \le \qu(\beta) \}$. Our goal is to prove the following. 

\begin{thm}\label{thm:nonzerolevelleaves}
There is a bijection between $\{ \beta \in \Xi(\bc) \, | \, \qu (\beta) \le n \}$ and the symplectic leaves of $\ZCB_{c}(\Gamma_n)$, $\beta \mapsto \LG(\beta)$, such that 
	\begin{enumerate}
		\item[(i)] $\dim \LG(\beta) = 2m$, and
		\item[(ii)] the leaf $\LG(\beta)$ is labeled by the parabolic conjugacy class $(\Gamma_m)$,
	\end{enumerate}    
where $m := n - \qu(\beta)$. 
\end{thm}

\begin{proof}
	Let $\LG$ be a leaf of $\XC_{\bc}(n \delta)$. Then, as explained in Section~\ref{sec:sympleavesdefn}, $\LG$ is labeled by some decomposition type $\tau = (\beta^{(0)},(n_1,\beta^{(1)};\ds; n_k, \beta^{(k)}))$. By Lemma~\ref{lem:basicrepCarternonLie}, $\beta^{(0)} = \gamma(m,\nu)$ for some $\nu \in \boldsymbol{\Lambda}$ and $m \ge 0$. Moreover, since $n \delta \in E_{\bc}$ by Lemma~\ref{lem:thingsinSigmac}(ii), we must have $0 \le m \le n$ with $m = n$ if and only if $\tau = (n \delta, (\emptyset))$. The roots $\beta^{(i)}$ for $i \ge 1$ belong to $R^+_{\bc}$. In particular, they are real roots since $\bc(\delta) = 1$; this forces $\beta^{(i)} \neq \beta^{(j)}$ for $i \neq j$ since $\LG \neq \emptyset$. 
We may write
	$$
	\beta := n_1 \beta^{(1)} + \cdots + n_k \beta^{(k)} \in \SS{\bc}.
	$$
	Since $\gamma(m,\nu) + \beta = n \delta$, we have $m \delta + (1/2)(\nu,\nu) \delta - \nu + \beta = n \delta$, or $m \delta + (1/2)(\nu,\nu) \delta - \nu = n \delta - \beta$. Taking the coefficient of $e_0$ gives $m + (1/2)(\nu,\nu) = n - \beta_0$ and $\nu = \beta - \beta_0 \delta$. Then
	\begin{align*}
		n - m & = \beta_0 + (1/2)(\nu,\nu) \\
		& = \beta_0 + (1/2)(\beta - \beta_0 \delta,\beta - \beta_0 \delta) \\
		& = \beta_0 + (1/2)(\beta,\beta) = \qu(\beta).
	\end{align*}
	Hence $0 \le \qu (\beta) \le n$ and $m = n - \qu(\beta)$. This implies that $\beta$ belongs to the set $\{ \beta \in \SS{\bc} \, | \, \qu(\beta) \le n \}$ and $\dim \LG = 2m = 2 (n - \qu(\beta))$.

Next we show that $\gamma(m,\nu) \in E_{\bc}$ implies that $\beta \in \Xi(\bc)$. If $\beta \notin \Xi(\bc)$ then there exists $\eta \succ \beta$ in $\SS{\bc}$ with $\varrho(\eta) \le \varrho(\beta)$. Let $m' = n - \varrho(\eta) \ge 0$ and $\nu' = \eta - \eta_0 \delta$ so that $\gamma(m',\nu') + \eta = n \delta$. Then  
\[
\gamma(m,\nu) - \gamma(m',\nu') = (n\delta - \beta) - (n \delta - \eta) = \eta - \beta \in\SS{\bc}. 
\]
But this implies that 
\begin{equation}\label{eq:decomposegammamnu}
    \gamma(m,\nu) = \gamma(m',\nu') + \zeta^{(1)} + \cdots + \zeta^{(s)}
\end{equation}
for some $\zeta^{(i)} \in \Delta(\bc)$. Since $p_{\Lambda_0}(\gamma(m',\nu')) = m' \ge m = p_{\Lambda_0}(\gamma(m,\nu))$ and $p(\zeta^{(i)}) = 0$, this contradicts the fact that $\gamma(m,\nu) \in E_{\bc}$. Thus, $\beta \in \Xi(\bc)$. 

Conversely, assume we have chosen $\beta \in \Xi(\bc)$ with $\varrho(\beta) \le n$. If $\nu := \beta - \beta_0 \delta \in \boldsymbol{\Lambda}$ then 
\[
n \delta - \beta = (n - \varrho(\beta)) \delta + \frac{1}{2} (\nu,\nu) \delta - \nu = \gamma(m,\nu) \ge 0,
\]
where $m = n - \varrho(\beta)$. Repeating the argument of the previous paragraphs shows that if $\gamma(m,\nu) \notin E_{\bc}$ then there must exist $m',\nu'$ and $\zeta^{(i)}$, with $m' \ge m$, such that decomposition \eqref{eq:decomposegammamnu} holds. But, again, if we write $\eta = n \delta - \gamma(m',\nu')$ then $\eta - \beta = \sum_i \zeta^{(i)} \in \SS{\bc}$ and hence $\eta \succ \beta$. Then $\varrho(\eta) = n - m' \le n - m = \varrho(\beta)$ implies that $\beta \notin \Xi(\bc)$. 

Part (ii). As shown in Corollary~\ref{cor:level1leaflabel}, the leaves are labeled by $(\Gamma_j)$ for some $j$ when $\bc(\delta) = -1$. The leaves of dimension $2m$ are all labeled by $(\Gamma_m)$, so (ii) follows from (i). 
\end{proof}

By Lemma~\ref{lem:parabolicwreathproduct}, the normalizer of $\Gamma_m$ in $\Gamma_n$ is just $\Gamma_m$ itself. Therefore, in this case, Losev's induction result says that there is bijection between the zero-dimensional leaves in $\ZCB_{c}(\Gamma_m)$ and the $2(n-m)$-dimensional leaves in $\ZCB_{c}(\Gamma_n)$. It is natural to ask what is the relation between Theorem~\ref{thm:nonzerolevelleaves} and Losev's induction result, Theorem~\ref{thm:losevinductionleaves}. This is clarified by the following corollary. For $\beta \in \Xi(\bc)$ with $\varrho(\beta) \le n$, write $\LG_n(\beta)$ for the corresponding leaf in $\ZCB_c(\Gamma_n)$.

\begin{cor}
The zero-dimensional leaves in $\ZCB_{c}(\Gamma_n)$ are labeled by $\{ \beta \in \Xi(\bc) \, | \, \varrho(\beta) = n \}$. The rule $\LG_{n}(\beta) \mapsto \LG_{n + k}(\beta)$, for $\beta \in \Xi(\bc)$ with $\varrho(\beta) = n$, defines a bijection  between zero-dimensional leaves in $\ZCB_{c}(\Gamma_n)$ and $2k$-dimensional leaves in $\ZCB_{c}(\Gamma_{n+k})$.  
\end{cor}

\begin{proof}
Theorem~\ref{thm:nonzerolevelleaves} makes it clear that for $0 \le m \le n$, the codimension $2m$ leaves in $\ZCB_{c}(\Gamma_n)$ are in bijection with the codimension $2m$ leaves in $\ZCB_{c}(\Gamma_{n+k})$ since they are both in bijection with the elements $\beta$ of $\Xi(\bc)$ satisfying $\qu(\beta) = m$. 
\end{proof}

The codimension two leaves are particularly easy to enumerate. Note that even though the root system $\Phi$ is irreducible, the system $\Phi_{\bc}$ is not irreducible in general. 

\begin{lem}\label{lem:codim2inCM}
	The codimension two leaves in $\ZCB_{c}(\Gamma_n)$ are in bjiection with the irreducible factors of $\Phi_{\bc}$. 
\end{lem}

\begin{proof}
	By Theorem~\ref{thm:nonzerolevelleaves}, we wish to find the vectors in $\Xi(\bc)$ with $\qu(\beta) = 1$. Recall $\qu(\beta) = \beta_0 + (1/2)(\beta,\beta)$. Since $R_{\bc}$ is a finite root system, $\beta \neq 0$ implies that $0 \neq (\beta,\beta) \in 2 \Z_{\ge 0}$. Thus, $(\beta,\beta) = 2$ and $\beta_0 = 0$. In other words, $\beta \in \SS{\bc} \cap \boldsymbol{\Lambda}$. Since $(\beta,\beta) = 2$, this forces $\beta \in \Phi_{\bc}$. Then, belonging to $\Xi(\bc)$ and $\Phi_{\bc}^+$ means that $\beta$ must be the highest positive root in one of the irreducible factors of $\Phi_{\bc}$. 
\end{proof}

\begin{example}
	Take $\Gamma = \Z/2\Z$ so that $\Gamma_n$ is the Weyl group of type $B$. Then $Q = \Z e_0 + \Z e_1$ and the real positive roots are 
 \[
 \{ m e_0 + (m+1) e_1 \, | \, m \in \Z_{\ge 0} \}\cup \{ m e_0 + (m-1) e_1 \, | \, m \in \Z_{> 0} \}.
 \]
Note that we assume $\bc (\delta) = -1$, which becomes $\bc_0 + \bc_1 = -1$. Hence, 
	$$
	\Delta(\bc) = \left\{ \begin{array}{ll}
		\{ m e_0 + (m+1) e_1 \} & \mbox{if } \bc = (m+1,-m), \, m \in \Z_{\ge 0}, \\
		\{ m e_0 + (m-1) e_1 \} & \mbox{if } \bc = (1-m,m), \, m \in \Z_{> 0}, \\
		\emptyset & \textrm{otherwise}.
	\end{array} \right. 
	$$
	Therefore, in the interesting cases, $\SS{\bc} = \N \alpha$ and $\qu(k\alpha)=k (k+m)$. Moreover, $\Xi(\bc) = \SS{\bc}$ in this case. Therefore, the leaves of $\ZCB_{c}(\Gamma_n)$ are in bijection with the set $\{ k \ge 0 \, | \, k(k+m) \le n\}$. 
	
More generally, if $\Gamma$ is again arbitrary and $R_{\bc}$ is a product of rank one root systems, equivalently $(\alpha, \beta) = 0$ for all $\alpha \neq \beta$ in $\Delta(\bc)$, then the above argument shows that the leaves of $\ZCB_{c}(\Gamma_n)$ are in bijection with the set $\{ k \in \N^s \, | \, \sum_i k_i(k_i+m_i) \le n\}$. Here $m_i$ is the coefficient of $e_0$ in $i$th root $\alpha$ in $\Delta(\bc)$. 
\end{example}

Finally, we describe the closure of leaves in $\ZCB_{c}(\Gamma_n)$. Recall that $\mathsf{G}'$ is the deframed graph corresponding to $\mathsf{G}$ and $\bw = \Lambda_0$. 

\begin{lem}\label{lem:basicrootfundemental}
	The root vector $e_{\infty} + \gamma(m,\nu)$ belongs to the fundamental region $F'$ for $\mathsf{G}'$ if and only if $\nu = 0$ and $m > 1$. 
\end{lem}

\begin{proof}
	We have 
	$$
	(e_{\infty} + \gamma(m,\nu),e_i) = \left\{ \begin{array}{ll}
		-(\nu,e_i) & i \ge 1 \\
		-(\nu,e_0) - 1 & i = 0 \\
		2 - (m + (1/2)(\nu,\nu)) & i = \infty. 
	\end{array} \right.
	$$
	Therefore, if $\nu = 0$ then $e_{\infty} + \gamma(m,\nu)$ belongs to the fundamental region if and only if $m > 1$. Thus, we must show that $e_{\infty} + \gamma(m,\nu)$ is not in the fundamental region if $\nu \neq 0$. 

	Consider the case $m = 0$. Let $v = \Lambda_0 - \frac{1}{2}(\nu,\nu) \delta + \nu$ be the corresponding  weight for the basic representation $L(\Lambda_0)$. By Lemma~\ref{lem:basicrepCarter}, it is dominant if and only if $\nu = 0$. But being dominant means that $(v, e_i) \ge 0$ for all $i = 0, \ds, r$. Since $(v,e_i) = (\nu, e_i)$ for $i \ge 1$ and $(v,e_0) = (\nu,e_0) + 1$, we deduce that the equations $-(\nu,e_i) \le 0$ for $i \ge 1$ and $-(\nu,e_0) - 1 \le 0$ forces $\nu = 0$. Thus, $\nu \neq 0$ implies that $e_{\infty} + \gamma(m,\nu)$ is not in the fundamental region $F'$. 
 
 Finally, if $e_{\infty} + \gamma(m,\nu) = e_{\infty}$ or $e_{\infty} + \delta$ then it is clearly not in the fundamental region.  
\end{proof}

\begin{prop}\label{prop:Cadmissiblemove}
	Let $m \ge 0$ and $\nu \in \boldsymbol{\Lambda}$. If $\gamma(m,\nu)$ belongs to $E_{\bc}$ then there exists a $\bc$-admissible $w \in W^{\aff}$ such that $w \star \gamma(m,\nu) = m \delta$. 
\end{prop}

\begin{proof}
As in the proof of Lemma~\ref{lem:basicrootfundemental}, we work with dimension vectors on the deframed graph $\mathsf{G}'$ (associated to the framing vector $\bw = \Lambda_0$). Then $s_i (e_{\infty} + \gamma) = e_{\infty} + (s_i \star \gamma)$ for $i = 0, \ds, r$. Let $F'$ denote the fundamental region for $\mathsf{G}'$. We must show that there is a $\bc$-admissible $w \in W^{\aff} \subset W_{\mathsf{G}'}$ such that $w(e_{\infty} + \gamma(m,\nu)) = e_{\infty} + m \delta$. 

Note that $w(e_{\infty} + \gamma(m,\nu))$ is a positive root for all $w \in W^{\aff}$ that are $\bc$-admissible since the coefficient of $e_{\infty}$ in $w(e_{\infty} + \gamma(m,\nu))$ is always one. Let $\beta$ be minimal among all such $w(e_{\infty} + \gamma(m,\nu))$ under the dominance ordering. We claim that $(\beta,e_i) \le 0$ for all $i \ge 0$. Indeed, if $\bc' = w^*(\bc)$ and $\bc_i' = 0$ then $(\beta,e_i) \le 0$ by \cite[Lemma~7.2]{CBmomap} since $\beta \in \Sigma_{\bc_{\beta}'}(\mathsf{G}')$. Therefore, if $(\beta,e_i) >0$ then necessarily $\bc_i' \neq 0$. But then $s_i$ is $\bc'$-admissible and $s_i(\beta) < \beta$, which contradicts the minimality of $\beta$. If we write $\beta = e_{\infty} + m' \delta + \frac{1}{2}(\nu',\nu') \delta - \nu'$ then $m' = p(\beta) = p(e_{\infty} + \gamma(m,\nu)) = m$. This implies that  
\begin{equation}\label{eq:betaeinfty}
    	(\beta,e_{\infty}) = 2 - (m + (1/2)(\nu',\nu')).
\end{equation}
If $m > 1$ then we deduce that $\beta$ is in the fundamental region since $(\nu',\nu') \ge 0$. In particular, $\beta  \in F' \cap (W_{\mathsf{G}'} \cdot (e_{\infty} + \gamma(m,\nu)))$. By Lemma~\ref{lem:basicrepCarter} and Lemma~\ref{lem:basicrootfundemental}, $e_{\infty} + m \delta \in F' \cap (W_{\mathsf{G}'} \cdot (e_{\infty} + \gamma(m,\nu)))$ too. But $|F' \cap (W_{\mathsf{G}'} \cdot (e_{\infty} + \gamma(m,\nu)))| \le 1$ by \cite[Proposition~3.12(b)]{KacBook}. The claim follows. 

If $m = 1$ then equation \eqref{eq:betaeinfty} becomes $(\beta,e_{\infty}) = 1 - (1/2)(\nu',\nu')$. Therefore, $\beta$ is again in $F'$ when $(\nu',\nu') > 0$ i.e. when $\nu' \neq 0$. But this contradicts Lemma~\ref{lem:basicrootfundemental}. If $\nu' = 0$ then $\beta = e_{\infty} + \delta = e_{\infty} + m \delta$. 

Finally, we consider the case $m = 0$. Again, if $(\nu',\nu') > 2$ then $\beta \in F'$ by \eqref{eq:betaeinfty}, which contradicts Lemma~\ref{lem:basicrootfundemental}. If $(\nu',\nu') = 2$ then $\nu'$ is in the finite root system $\Phi$. But then $(\beta,e_i) = - (\nu',e_i) \ge 0$ for $i = 1, \ds, r$ implies that $\nu' = \theta$ is the highest positive root and hence $\beta = e_{\infty} + e_0$. But then $(\beta, e_0) = 2-1 = 1$ contradicting the assumption that $(\beta,e_i) \le 0$ for all $i \ge 0$. If $\nu' = 0$ then $\beta = e_{\infty} = e_{\infty} + m \delta$, as required.  \qedhere

\end{proof}    

Note that the equation \eqref{eq:boldc} defines an isomorphism between the vector spaces of parameters $c$ for the symplectic reflection algebra and the space $\Hom_{\Z}(Q,\C)$. Therefore, it implicitly defines an action of $W^{\aff}$ on the space of parameters, which we also write $c \mapsto w^*(c)$, for $w \in W^{\aff}$.  

\begin{thm}\label{thm:leafclosurenormal}
	Let $\beta \in \Xi(\bc)$ with $m := n - \qu(\beta) \le n$. There exists a $\bc$-admissible element $w \in W^{\aff}$ such that $\overline{\LG_n(\beta)} \cong \ZCB_{w^*(c)}(\Gamma_m)$. In particular, the leaf closure $\overline{\LG_n(\beta)}$ is normal. 
\end{thm}

\begin{proof}
	We let $\nu = \beta - (\Lambda_0,\beta) \delta$. Then, by Proposition~\ref{prop:Cadmissiblemove}, there exists a $\bc$-admissible element $w \in W^{\aff}$ such that $w \star \gamma(m,\nu) = m \delta$. Write $\beta = \sum_{i = 1}^s n_i \eta^{(i)}$ with $\eta^{(i)} \in \Delta(\bc)$. Then the proof of Theorem~\ref{thm:nonzerolevelleaves} shows that $\LG_n(\beta)$ is the leaf labeled by the representation type
 \[
 \tau = (\gamma(m,\nu),(n_1,\eta^{(1)};\ds; n_s,\eta^{(s)})). 
 \]
Since every $\eta^{(i)}$ is a real root, Proposition~\ref{prop:leafclosurenormal} implies that $\overline{\LG}(\beta) \cong \XC_{\bc}(\gamma(m,\nu))$ is normal. The fact that $w$ is $\bc$-admissible means that 
 \[
 \XC_{\bc}(\gamma(m,\nu)) \cong \XC_{w^*(\bc)}(m \delta) \cong \ZCB_{w^*(c)}(\Gamma_m). \qedhere
 \]
\end{proof}

\begin{rem}\label{rem:finiteWisoparameters}
For any $w$ in the finite Weyl group $W$, the Calogero--Moser varieties $\ZCB_{c}(\Gamma_n)$ and $\ZCB_{w^*(c)}(\Gamma_n)$ are isomorphic. This follows from the isomorphism $\XC_{\bc}(n \delta) \cong \XC_{w^*(\bc)}(n \delta)$, which can be deduced from the reflection isomorphisms of Section~\ref{subs:adm-refl} just as in     \cite[Corollary~3.6]{BM} noting that if $1 \le i \le r$ such that $\bc_i = 0$ then we still have $\XC_{\bc}(n \delta) \cong \XC_{s_i^*(\bc)}(n \delta)$ because $s_i^* \bc = \bc$ in this case. 
\end{rem}

Noting that we assume $\bc(\delta) = -1$, the parameter $w^*(\bc)$ can be computed explicitly in terms of the vector $\beta$ as follows.

\begin{prop}\label{prop:wbccomp}
    	Let $\beta \in \Xi(\bc)$ with $m := n - \qu(\beta) \le n$. Then $\overline{\LG_n(\beta)} \cong \XC_{\bc'}(m \delta)$, where $\bc' =\bc + \overline{\beta}$.
\end{prop}

\begin{proof}
We continue with the setup of the proof of Theorem~\ref{thm:leafclosurenormal}. In particular, $w \in W^{\aff}$ is a $\bc$-admissible element with $w \star \gamma(m,\nu) = m \delta$ and $\gamma(m,\nu) + \beta = n \delta$. Recall from Section~\ref{sec:anotherpresentation} that $w$ can be written $w = x t_{u}$, for some (unique) $x \in W$ and $u \in Q$. Then $w \star \gamma(m,\nu) = m \delta$, together with Lemma~\ref{lem:affWeyl_translation-roots}, means that there exists $r \in \Z$ such that 
    \begin{align*}
      e_{\infty} + m \delta & =  w(e_{\infty} + m \delta + (1/2)(\nu,\nu) \delta - \nu) \\
      & = x t_u(e_{\infty} + m \delta + (1/2)(\nu,\nu) \delta - \nu) \\
        & = x (e_{\infty} + m \delta + (1/2)(\nu,\nu) \delta - \nu - u + r \delta) \\
        & = e_{\infty} + (m + r + (1/2)(\nu,\nu)) \delta - x(\nu + u).
    \end{align*}
    Comparing the coefficient of $e_0$ in these equalities implies that $m + r + (1/2)(\nu,\nu) = m$ i.e. $r = - (1/2)(\nu,\nu)$. Then $x(\nu + u) = 0$ i.e. $u = - \nu$. 

    Therefore, Lemma~\ref{lem:affWeyl_translation-theta} says that $w^*(\bc) = x^* t^*_{-\nu} (\bc) = x^*(\bc - \bc(\delta) \overline{\nu})$. Recall that $\nu = \beta - (\Lambda_0,\beta) \delta$ and hence $\overline{\beta} = \overline{\nu}$. Since $\bc(\delta) = -1$, it follows that $\overline{\LG_n(\beta)} \cong \XC_{x^*(\bc')}(m \delta)$. But $\XC_{x^*(\bc')}(m \delta) \cong \XC_{\bc'}(m \delta)$ by Remark~\ref{rem:finiteWisoparameters}.   
\end{proof}

The partial order on leaves is described as follows. 

\begin{prop}\label{prop:leafclosure1}
Assume $\beta, \eta \in \Xi(\bc)$ with $\beta \neq \eta$ and $\qu(\beta), \qu(\eta) \le n$. Then $\overline{\LG(\beta)} \cap \LG(\eta) \neq \emptyset$ in $\ZCB_{c}(\Gamma_n)$ if and only if $\eta \succ \beta$. 
\end{prop}

\begin{proof}
For brevity, set $\alpha = e_{\infty} + n \delta$. 

First, we assume that $\eta \succ \beta$. Then there exist simple $\Pi^{\bc_{n \delta}}$-modules $M_0,N_0$ and semi-simple $\Pi^{\bc_{n \delta}}$-modules $M_1 = N_{11}$ and $N_{12}$ such that $\dim M_0 = \alpha - \beta, \dim M_1 = \beta, \dim N_{0} = \alpha - \eta$ and $\dim N_{12} = \eta - \beta$. Let $N_1 = N_{11} \oplus N_{12}$ so that $\dim N_1 = \eta$. Then $[M_0 \oplus M_1] \in \LG(\beta)$ and $[N_0 \oplus N_{1}] \in \LG(\eta)$. To show that $\overline{\LG(\beta)} \cap \LG(\eta) \neq \emptyset$ it suffices, by \cite[Proposition~3.6]{BellSchedQuiver} and \cite[Corollary~3.25]{BellSchedQuiver}, to argue that the $M_i$ and $N_i$ can be chosen such that $\Stab_G (M_0 \oplus M_1) \subset \Stab_G(N_0 \oplus N_1)$, where $G = G(\alpha)$. Since there are no homomorphisms between $M_0$ and the $N_i$ or between $N_0$ and $N_{11},N_{12}$, we have 
	\begin{align*}
		\Stab_G (M_0 \oplus M_1) & = \mr{Aut}_{\Pi^{\bc}}(M_0) \times \mr{Aut}_{\Pi^{\bc}}(M_1) \cong \Cs \times \mr{Aut}_{\Pi^{\bc}}(M_1), \\
		\Stab_G (N_0 \oplus N_1) & = \mr{Aut}_{\Pi^{\bc}}(N_0) \times \mr{Aut}_{\Pi^{\bc}}(N_{1}) \cong \Cs \times \mr{Aut}_{\Pi^{\bc}}(M_{1} \oplus N_{12}). 
	\end{align*}
	The claim follows. 

 Conversely, if $\overline{\LG(\beta)} \cap \LG(\eta) \neq \emptyset$ with $\eta \neq \beta$ then there exists $[M] \in \LG(\beta)$ and $[N] \in \LG(\eta)$ such that $\Stab_G(M) \subsetneq \Stab_G(N)$. We wish to argue that $\eta \succ \beta$. As in the previous paragraph, we decompose $M = M_0 \oplus M_1$ and $N = N_0 \oplus N_1$. If $\Delta(\bc) = \{ \eta^{(1)}, \ds, \eta^{(s)} \}$ and $L(\eta^{(i)})$ are the simple $\Pi^{\bc}$-modules with $\dim L(\eta^{(i)}) = \eta^{(i)}$ then $M_1 = \bigoplus_{i = 1}^s L(\eta^{(i)}) \o U_i$, where $\dim U_i = \beta_i$, and 
 \[
 \Stab_{G(\alpha)}(M) = \Cs \o \bigotimes_{i = 1}^s GL(U_i), 
 \]
with first factor $\Cs = \mr{Aut}_{\Pi^{\bc}}(M_0)$. Let $\Cs_{M_1} \subset  \Stab_{G(\alpha)}(M)$ be the torus that acts by weight one on $M_1$ and acts trivially on $M_0$. Let $V_j$ be the vector space at vertex $j$ as in \eqref{eq:repspace}, with $V_{\infty} = W_0 \cong \C$ since $\alpha_{\infty} = 1$. Then $\Stab_{G(\alpha)}(M)$ acts on each $V_j$ and $M_1 = \bigoplus_j (M_1 \cap V_j)$ as a $\Stab_{G(\alpha)}(M)$-module. Notice that 
\[
V_j = (M_0 \cap V_j) \oplus (M_1 \cap V_j) = V_j(0) \oplus V_j(1),
\]
where $V_j(m)$ is the subspace of $V_j$ of weight $m$ with respect to $\Cs_{M_1}$. We have $M_0 \cap V_{\infty} = V_{\infty} \cong \C$ because $M_1 \cap V_{\infty} = 0$. Now, by assumption, $\Cs_{M_1} \subset \Stab_{G(\alpha)}(N)$. Crucially, $\Cs_{M_1}$ acts trivially on $N_0$. Indeed, $N_0$ is a simple $\Pi^{\bc}$-module so $\Cs_{M_1}$ must act as a scalar. But $V_{\infty} = M_0 \cap V_{\infty} = N_0 \cap V_{\infty} \subset N_0$ so $\Cs_{M_1}$ must act trivially on $N_0$. Therefore, for every $j \ge 0$, $V_j(0) = M_0 \cap V_j \supset N_0 \cap V_j$ which means that 
 \[
 \alpha_j - \beta_j = \dim M_0 \cap V_j \ge \dim N_0 \cap V_j =  \alpha_j - \eta_j
 \]
 and hence $\eta \succeq \beta$. If $\eta = \beta$ then $\Stab_G(M)$ would be a proper subgroup of $\Stab_G(N)$ which is also conjugate to $\Stab_G(N)$, which is impossible. Thus, $\eta \succ \beta$. 
 \end{proof}

\section{Zero level}\label{sec:zerolevel}

In this section we consider the parameters $c$ for which $\bc(\delta) = 0$. 

\subsection{} Recall that $R$ is the affine root system associated to the finite group $\Gamma$ and $\Phi$ the standard finite subsystem. 

\begin{lem}\label{lem:degeneratelambdaroots}
	If $\bc(\delta) = 0$ then $R_{\bc}^+ = \Phi^+_{\bc} \cup \{ m \delta \pm \alpha \, | \, m > 0, \alpha \in \Phi^+_{\bc} \} \cup \Z_{> 0} \delta$. 
\end{lem}

\begin{proof}
	Recall that $R^+ = \Phi^+ \cup \{ m \delta \pm \alpha \, | \, m > 0, \alpha \in \Phi^+ \} \cup \Z_{> 0} \delta$. The result follows.  
\end{proof}

This means that $R_{\bc} = (\Phi_{\bc} + \Z \delta) \cup (\Z \smallsetminus \{ 0 \} )\delta$. 

\begin{lem}\label{lem:Sigmalambdadegenerate}
	Assume $\bc(\delta) = 0$. 
 \begin{enumerate}
     \item[(i)] $\Sigma_{\bc} = \{ \textrm{minimal elements in } \Phi_{\bc}^+ \} \cup \{ \delta \} \cup \{ \delta - \alpha \ | \ \alpha \textrm{ maximal elements in } \Phi_{\bc}^+ \}$.
 \item[(ii)] If $\alpha \in E_{\bc}$ then $p_{\Lambda_0}(\alpha) = 0$ i.e. $e_{\infty} + \alpha$ is a real root. 
 \end{enumerate}
\end{lem}

\begin{proof}
	Part (i). It follows directly from Lemma~\ref{lem:degeneratelambdaroots} that the set $\{ \textrm{minimal elements in } \Phi_{\bc}^+ \} \cup \{ \delta \}$ is contained in $\Sigma_{\bc}$. Consider now $m \delta \pm \alpha$, where $m > 0$ and $\alpha \in \Phi_{\bc}^+$. The element $m \delta + \alpha$ can be written as the sum $\delta + \cdots + \delta$ plus a sum of vectors in $\{ \textrm{minimal elements in } \Phi_{\bc}^+ \}$ (adding up to $\alpha$) which shows that the real root $m \delta + \alpha$ does not belong to $\Sigma_{\bc}$. Similarly, if $m > 1$ then $m \delta - \alpha$ can be written as $m-1$ copies of $\delta$ plus the real root $\delta - \alpha$, again implying that it does not belong to $\Sigma_{\bc}$.   
	
	Finally, we consider $\delta - \alpha$, where $\alpha \in \Phi_{\bc}^+$. If $\alpha$ is not maximal then it is clear that $\delta - \alpha$ does not belong to $\Sigma_{\bc}$. If $\alpha$ is maximal, then the only roots in $R^+_{\bc}$ that are less that $\delta - \alpha$ are all of the form $\beta \in \Phi^+_{\bc}$, where $\alpha$ and $\beta$ belong to different irreducible factors of $\Phi^+_{\bc}$. In this case, we can never express $\delta - \alpha$ as a sum of such $\beta$; in other words, $\delta - \alpha$ is minimal in $R^+_{\bc}$. It follows that this vector belongs to $\Sigma_{\bc}$.  

 Part (ii). If $\alpha \in E_{\bc}$ then, by Lemma~\ref{lem:basicrepCarternonLie}, $\alpha = m \delta + (1/2) (\nu,\nu) \delta - \nu$. Let $\beta = (1/2) (\nu,\nu) \delta - \nu$, so that $e_{\infty} + \beta$ is a positive real root for $\mathsf{G}'$, again by Lemma~\ref{lem:basicrepCarternonLie}. Notice that $\bc(\delta) = 0$ and $\bc_{\alpha}(e_{\infty} + \alpha) = 0$ implies that $\bc_{\alpha}(e_{\infty} + \beta) = 0$ too. 
 If $m>0$, then $\alpha \in E_{\bc}$ implies 
 \[
 m = p_{\Lambda_0}(\alpha) > p_{\Lambda_0}(\beta) + p(\delta) + \cdots + p(\delta) = m.
 \]
Then we have $m = 0$. 
\end{proof}

\subsection{The symplectic leaves}

We enumerate the irreducible factors of $\Phi_{\bc}$ as 
$$
\Phi_{\bc} = \Phi_1 \cup \cdots \cup \Phi_{\ell}. 
$$
We recall that a composition of length $\ell$ is an $\ell$-tuple $(\rho_1, \ds, \rho_{\ell})$ of non-negative integers. 

\begin{prop}\label{prop:labelleaveslevel0}
	The symplectic leaves $\LG(\lambda,\rho)$ in $\ZCB_{c}(\Gamma_n)$ are labeled by pairs $(\lambda,\rho)$, where $\lambda \in \mc{P}$ is a partition, $\rho$ is a composition of length $\ell$ and $|\lambda| + |\rho| = n$.       
\end{prop}

\begin{proof}
	For each $i$, we denote by $\alpha(i)_1, \ds, \alpha(i)_{r_i}$ the simple roots in $\Phi_i^+$. If $\theta(i)$ is the longest root in $\Phi_i$ then 
	$$
	\theta(i) = h(i)_1 \alpha(i)_1 + \cdots + h(i)_{r_i} \alpha(i)_{r_i}.
	$$
	Let us explain how to associate to the pair $(\lambda,\rho)$ a decomposition of $n \delta$. Lemma~\ref{lem:Sigmalambdadegenerate}(1) implies that the root $\delta$ admits a decomposition $(1,\delta - \theta(i); h(i)_1,\alpha(i)_1; \dots ; h(i)_{r_i},\alpha(i)_{r_i})$ in $\Sigma_{\bc}$. Multiplying through by $\rho_i$ gives a decomposition of $\rho_i \delta$. Then the decomposition of $n \delta$ corresponding to $(\lambda,\rho)$ is 
	\begin{equation}\label{eq:degeneratelambdadecomp}
		(0,(\lambda_1,\delta;\lambda_2,\delta;\ds;\rho_1,\delta - \theta(1);\rho_1 h(1)_1, \alpha(1)_1; \dots ;\rho_1 h(1)_{r_1}, \alpha(1)_{r_1};\rho_2,\delta - \theta(2);\ds )),    
	\end{equation}
see \eqref{eq:taureptype} for the notation. We check that these are the only decompositions of $n \delta$. 

First, it is shown in \cite[Theorem~1.1]{CBdecomp} that $e_{\infty} + n \delta$ admits a \textit{canonical decomposition} with respect to $\Sigma_{\bc_{n \delta}}(\mathsf{G}')$ such that any other decomposition of $e_{\infty} + n \delta$ is a refinement of this decomposition. In this particular case, the canonical decomposition is computed in \cite[Proposition~4.2]{BellamyCrawQuotient} and equals $(0,(1, \delta;\ds; 1, \delta))$. Therefore, we must have $\beta^{(0)} = 0$ in every decomposition; that is, every decomposition is of the form $(0,(n_1,\beta^{(1)};\ds))$ where $(n_1,\beta^{(1)};\ds)$ is a decomposition of $n \delta$ in $\Sigma_{\bc}$. The occurrences of $\delta$ (with multiplicity) in such a decomposition define a partition $\lambda$. Discarding these, we are left with decompositions of $m \delta$ using only the roots in $\{ \textrm{minimal elements in } \Phi_{\bc}^+ \}$ and $\{ \delta - \alpha \ | \ \alpha \textrm{ maximal elements in } \Phi_{\bc}^+ \}$. All these roots are real, so only occur once (with multiplicity) in any decomposition. If $\delta - \theta(i)$ occurs (with multiplicity $\rho_i$ say), then the $\alpha(i)_j$ must also occur with multiplicity $\rho_i h(i)_j$. This implies that our decomposition has the form \eqref{eq:degeneratelambdadecomp}. 
\end{proof}

\begin{example}
	If $\bc = 0$, then $\Sigma_0 = \{ \delta, e_0, \ds, e_r \}$ and $\ZCB_{c}(\Gamma_n) = V^* / \Gamma_n$. Write $\delta = \sum_{i = 0}^r \delta_i e_i$. The decompositions of $n \delta$ are all of the form 
 \[
 \tau = (0,(\lambda_1, \delta; \ds;\lambda_k, \delta; m \delta_0, e_0; \ds; m \delta_r,e_r))
\]
where $n = m + \sum_{i = 0}^r \lambda_i$. Here, $m$ is acting as a composition of length $\ell = 1$. The leaf $\LG_{\tau}$ is labeled by the parabolic $\Gamma_m \times \s_{\lambda}$ and we see that there is a bijection between leaves and parabolic subgroups, as expected.  
\end{example}

\subsection{Leaf closures}

The combinatorics somewhat obscure the geometry of the situation when $\bc (\delta) = 0$. Write $X_{\bc}$ for the Calogero--Moser variety $\ZCB_{c}(\Gamma)$. This is the two-dimensional Calogero--Moser variety obtained when $n = 1$. It is a deformation of the Kleinian singularity $\C^2 / \Gamma$. Then $\ZCB_{c}(\Gamma_n) \cong S^n X_{\bc}$; see e.g.  \cite[Section~6.4]{BellamyThesis}. The variety $X_{\bc}$ has an open leaf $X_{\bc}^{\reg}$ and $\ell$ zero dimensional leaves $\{ p_i \}$, in bijection with the irreducible components $\Phi_i$ of $\Phi_{\bc}$. Then the leaves of $\ZCB_{c}(\Gamma_n)$ are of the form $S^{\lambda} X_{\bc}^{\reg} \times \{ \rho_1 p_1 + \cdots + \rho_{\ell} p_{\ell} \}$. 

The ordering on leaves is rather unnatural from the point of view of partition combinatorics. First, we consider points in the regular locus of $X_{\bc}$. Let $\pi_{\bc} \colon X_{\bc}^n \to S^n X_{\bc}$ be the quotient map. Then $S^{\eta} X_{\bc}^{\reg} \subset \overline{S^{\lambda} X_{\bc}^{\reg}}$ if and only if the stabilizer $\s_{\lambda}$ of a point in $\pi^{-1}_{\bc}(S^{\lambda} X_{\bc}^{\reg})$ is conjugate to a subgroup of the stabilizer $\s_{\eta}$ of a point in $\pi^{-1}(S^{\eta} X_{\bc}^{\reg})$. This happens precisely when we can write each row $\eta_i$ of the partition $\eta$ as a sum $\lambda_{i,1} + \cdots + \lambda_{i,r_i}$ of rows from $\lambda$. This motivates the following definition. 

\begin{defn}\label{defn:consituentparitions}
Let $\eta,\lambda$ be partitions and $\zeta,\rho$ compositions of length $\ell$. 
\begin{enumerate}
    \item[(i)] $\lambda$ is a \textit{constituent} of $\eta$ if $\eta = (\eta_1, \ds, \eta_s)$ and there exist $1 \le j_1 < \ds < j_k \le s = \ell(\eta)$ such that, after reordering the parts of $\lambda$, 
    \[
    \lambda = (\lambda_{1,1},\ds, \lambda_{1,r_1},\lambda_{2,1}, \ds,\lambda_{2,r_2},\ds,\lambda_{k,1},\ds, \lambda_{k,r_k})
    \]
    with $\eta_{j_i} = \lambda_{i,1} + \cdots + \lambda_{i,r_i}$ for each $1 \le i \le k$. 
    \item[(ii)] We say that $(\eta,\zeta) < (\lambda,\rho)$ if and only if $(\eta,\zeta) \neq (\lambda,\rho)$ and the following two conditions hold. First, $\lambda$ is a constituent of $\eta$. Secondly, if $1 \le j_1 < \ds < j_k \le \ell(\eta)$ are as in part (i) then there is a partition of the complement $I_1 \sqcup \cdots \sqcup I_{\ell}$ of $\{ 1, \ds, s \} \setminus \{ j_1, \ds, j_k \}$ such that 
    \[
    \zeta_i = \rho_i + \sum_{j \in I_i} \eta_j
    \]
    for all $1 \le i \le \ell$. 
\end{enumerate}
\end{defn} 

\begin{prop}\label{prop:leafrelationlevelzero}
    The closure order on leaves is given by $\LG{(\eta,\zeta)} \subset \overline{\LG}{(\lambda,\rho)}$ with $\LG{(\eta,\zeta)} \neq {\LG}{(\lambda,\rho)}$ if and only if $(\eta,\zeta) < (\lambda,\rho)$.     
\end{prop}

\begin{proof}
Given $(\lambda,\rho)$ with $\lambda = (\lambda_1,\ds, \lambda_s)$, let us partition 
\begin{equation}\label{eq:partitionJsn}
    \{1, \ds, n\} = J_1 \sqcup \cdots \sqcup J_s \sqcup J_{s+1} \sqcup \cdots \sqcup J_{s + \ell}
\end{equation}
linearly (that is, if $a \in J_i, b \in J_j$ and $i < j$ then $a < b$) such that $|J_i| = \lambda_i$ for $i \le s$,  and $|J_{i + s}| = \rho_i$ for $1 \le i \le \ell$. Then we define $C^{(\lambda,\rho)} X_{\bc}$ to be the set of points $(x_1,\ds, x_n) \in X_{\bc}^n$ such that (a) $x_i = x_j$ if and only if $i,j \in J_k$ for some $k$, (b) if $i \in J_j$ for some $j \le s$ then $x_i \in X_{\bc}^{\reg}$, and (c) if $i \in J_{j+s}$ for some $1 \le j \le \ell$ then $x_i = p_j$. The locally closed set $C^{(\lambda,\rho)} X_{\bc}$ is a connected component of 
\[
\pi^{-1}(S^{\lambda} X_{\bc}^{\reg} \times \{ \rho_1 p_1 + \cdots + \rho_{\ell} p_{\ell} \}) = \bigcup_{\sigma \in \s_n} \sigma(C^{(\lambda,\rho)} X_{\bc}).
\]
It is straight-forward to check that $C^{(\eta,\zeta)} X_{\bc}$ is conjugate (under $\s_n$) to a subset of $\overline{C^{(\lambda,\rho)} X_{\bc}}$ if and only if $(\eta,\zeta) \le (\lambda,\rho)$. Then the claim of the proposition follows from the fact that $\pi_{\bc}$ is a finite surjective map and 
\[
\overline{\LG}{(\lambda,\rho)} = \pi_{\bc}\left(\overline{\pi^{-1}_{\bc}(S^{\lambda} X_{\bc}^{\reg} \times \{ \rho_1 p_1 + \cdots + \rho_{\ell} p_{\ell} \})}\right) = \pi_{\bc}\left(\overline{C^{(\lambda,\rho)} X_{\bc}}\right).  \qedhere
\]
\end{proof}

\begin{lem}\label{lem:parabolicleaflevel0}
    The conjugacy class of parabolics associated to $\LG{(\lambda,\rho)}$ is $(\Gamma_{|\rho|} \times \s_{\lambda})$.
\end{lem}

\begin{proof}
    Let $\mf{p} \lhd \Zb_{\bc}(\Gamma_n)$ be the prime Poisson ideal defining the closure of $\LG_{(\lambda,\rho)}$. We must show that $\gr (\mf{p}) \lhd \C[V^*]^{\Gamma_n}$ defines the locus of points whose stabilizer is conjugate to a parabolic containing $\Gamma_{|\rho|} \times \s_{\lambda}$. We use the notation in the proof of Proposition~\ref{prop:leafrelationlevelzero}. 

    Note that $X_{\bc}^n = \Spec \Zb_{c}(\Gamma^n) = \Spec \Zb_{c}(\Gamma)^{\o n}$ and hence $\Zb_{c}(\Gamma_n) = \Zb_{c}(\Gamma^n)^{\s_n}$. Let $\mf{q} \lhd \Zb_{c}(\Gamma^n)$ be the prime ideal defining the closure of $C^{(\lambda,\rho)} X_{\bc}$ in $X_{\bc}^n$. The variety $V^* / \Gamma_n = X_0^n / \s_n$ has leaves $\LG{(\lambda,m)}$ labeled by conjugacy classes $(\Gamma_m \times \s_{\lambda})$ of parabolic subgroups. Here we think of $m$ as a composition of length one. The preimage of $\LG_{(\lambda,m)}$ in $X_0^n$ has irreducible components the images of $C^{(\lambda,m)} X_0$ under the action of $\Gamma_n$.  Assume that we have shown that the zero set of the ideal $\gr(\mf{q})$ contains the closure of $C^{(\lambda,m)} X_{0}$. Then the fact that $\mf{p} = \mf{q}^{\s_n}$ and that $\mc{F}_i \Zb_{c}(\Gamma_n) = (\mc{F}_i \Zb_{c}(\Gamma^n))^{\s_n}$ imply that $\gr(\mf{p}) = (\gr(\mf{q}))^{\s_n}$. Thus, the zero set of $\gr(\mf{p})$ contains $\pi_0(\overline{C^{(\lambda,m)} X_0}) = \overline{\LG}{(\lambda,m)}$. But the main result of \cite{MartinoAssociated} says that the zero set of $\gr(\mf{p})$ is the closure of a leaf in $S^n X_0 = V^* /\Gamma_n$ of dimension 
    \[
    \dim \LG{(\lambda,\rho)} = 2 |\lambda| = \dim \LG{(\lambda,m)}. 
    \]
    We deduce that the zero set of $\gr(\mf{p})$ is $\overline{\LG}{(\lambda,m)}$ as required. 

    Therefore, we are reduced to showing that the zero set of the ideal $\gr(\mf{q})$ contains the closure of $C^{(\lambda,m)} X_{0}$. If $\mf{m}_i \lhd \Zb_{c}(\Gamma)$ is the maximal ideal defining the point $p_i \in X_{\bc}$ and $\mf{q}_{\lambda}$ the ideal defining the closure of $C^{\lambda} X_{\bc}$ in $X_{\bc}^{n-m}$ then 
    \[
    \mf{q} = \mf{q}_{\lambda} \o \Zb_{c}(\Gamma)^{\o m} + \sum_{j = 1}^m \Zb_{c}(\Gamma)^{\o (n-m + j-1)} \o \mf{m}_{i_j} \o \Zb_{c}(\Gamma)^{\o (m-j)}.
    \]
    Since the filtration on $\Zb_{c}(\Gamma)^{\o n}$ is the product filtration, 
    \[
    \gr (\mf{q}) =  \gr(\mf{q}_{\lambda}) \o \Zb_{0}(\Gamma)^{\o (n-m)} + \sum_{j = 1}^m \Zb_{0}(\Gamma)^{\o (n-m + j-1)} \o \mf{m}_{0} \o \Zb_{0}(\Gamma)^{\o (m-j)},
    \]
    where $\mf{m}_0 = \gr(\mf{m}_i)$ is the augmentation ideal in $\Zb_0(\Gamma) = \C[\C^2]^{\Gamma}$. Thus, we may assume that $m = 0$ and $\lambda$ a partition of $n$. 
    
Recall that for a filtered algebra $R = \bigcup_i \mc{F}_i R$, the symbol $\sigma(x)$ of a non-zero element $x \in R$ is the image of $x$ in $(\mc{F}_i R) / (\mc{F}_{i-1} R) \subset \gr R$, where $i$ is the smallest integer such that $x \in \mc{F}_i R$. Recall also that $\Zb_{0}(\Gamma) = \C[\mathbf{A},\mathbf{B},\mathbf{C}] / (F)$, for some polynomial $F$. If $A,B,C \in \Zb_{c}(\Gamma)$ are elements whose symbols $\mathbf{A} = \sigma(A),\mathbf{B} = \sigma(B),\mathbf{C} = \sigma(C)$ are generator for $\Zb_0(\Gamma)$ then $\Zb_{0}(\Gamma^n)$ is generated by the symbols of the $A_i,B_i,C_i$ for $i = 1,\ds, n$, where $A_1 = A \o 1 \o \cdots \o 1$ etc. Using the partition \eqref{eq:partitionJsn}, the ideal $\mf{q}$ is generated by all $A_i - A_j, B_i - B_j, C_i - C_j$, where $1 \le i < j \le n$ with $i,j$ both belonging to some $J_k$ for $k \le s$. But then the ideal defining the closure of $C^{(\lambda,0)} X_{0}$ is generated by the $\sigma(A_i) - \sigma(A_j),\sigma(B_i) - \sigma(B_j), \sigma(C_i) - \sigma(C_j)$. Since $A_i$ and $A_j$ have the same degree under the filtration we have $\sigma(A_i - A_j) = \sigma(A_i) - \sigma(A_j)$ etc. It follows that $C^{(\lambda,0)} X_{0} \subset V(\gr(\mf{q}))$, as claimed. 
\end{proof}

In general, leaf closures are not normal when $\bc(\delta) = 0$. For instance, when $\Gamma = \{ 1 \}$ and $\bc = 0$, we have $\Gamma_n = \s_n$ and $\ZCB_{0} (\s_n) = \C^{2n} / \s_n$. This is the situation considered in Example~\ref{ex:symsn}. 

Let $\widetilde{\LG}{(\lambda,\rho)}$ denote the normalization of $\overline{\LG}{(\lambda,\rho)}$. Recall that if $\lambda = (\lambda_1,\lambda_2,\ds)$ is a partition then we set $\s(\lambda) := \prod_{i \ge 1} \s_{n_i}$ where $n_i := | \{ j \, | \, \lambda_j = i \} |$.

\begin{prop}
    Assume $\bc(\delta) = 0$. Then 
    \begin{equation}\label{eq:leafnormalizationlevel0}
     \widetilde{\LG}{(\lambda,\rho)} \cong \prod_{i = 1}^k \ZCB_{c}(\Gamma_{n_i}) = \ZCB_{c}(N_{\Gamma_n}(\mr{P})/\mr{P}),        
    \end{equation}
    where $(\mr{P} = \Gamma_m \rtimes \s_{\lambda})$ is the parabolic conjugacy class associated to $\LG_{(\lambda,\rho)}$.  
\end{prop}

\begin{proof}
Note that the conjugacy class of parabolic associated to $\LG{(\lambda,\rho)}$ is specified in Lemma~\ref{lem:parabolicleaflevel0}. The composition $\rho$ corresponds to the representation type   
  \[
  (\rho_1,\delta - \theta(1);\rho_1 h(1)_1, \alpha(1)_1; \dots ;,\rho_1 h(1)_{r_1}, \alpha(1)_{r_1};\rho_2,\delta - \theta(2);\ds ),
  \]
  all of whose terms are real roots. Therefore, Theorem~\ref{thm:quiverleafclosureagain} says that 
    \[
    \widetilde{\LG}{(\lambda,\rho)} \cong \mf{M}_{\bc}(\mathsf{G}(\Gamma),\delta)^{s}/\s(\lambda) \cong \prod_{i = 1}^k  X_{\bc}^{n_i}/\s_{n_i}, 
    \] 
    where $s = \ell(\lambda) = n_1 + \cdots + n_k$ and the second isomorphism is due to the fact that $\mf{M}_{\bc}(\delta) = X_{\bc}$. Since $\ZCB_{c}(\Gamma_n) = S^n X_{\bc}$ when $\bc (\delta) = 0$, the first isomorphism of \eqref{eq:leafnormalizationlevel0} follows. The second isomorphism is just the description of $N_{\Gamma_n}(\mr{P})/\mr{P}$ given in Lemma~\ref{lem:parabolicwreathproduct}. 
\end{proof}

\section{Combinatorics}\label{sec:combinatoricsCM}

In this section we introduce the additional combinatorics required to treat in greater generality the case where $\Gamma$ is of type $\mathsf{A}$ i.e. $\Gamma$ is a cyclic group. 

Assume $\ell\in \ZM_{>0}$. We set $\Gamma=\ZM/\ell\ZM$. The associated set of simple roots is $\Delta=\{e_0,e_1,\ldots,e_{\ell-1}\}$. The minimal imaginary root is $\delta=\sum_{i\in \ZM/\ell\ZM}e_i$.
We consider the following quiver $\Qu^\ell$. The set of vertices of $\Qu^\ell$ is $\Delta$ (we also identify this set with $\ZM/\ell\ZM$) and the arrows of $\Qu^\ell$ are of the form $e_i\to e_{i+1}$ for each $i\in\ZM/\ell\ZM$. 

We also allow $\ell=\infty$. In this case, we mean that $\Qu^\infty$ is the infinite linear quiver with set of vertices $\Delta=\{e_i,\,i\in\ZM\}$ and arrows $e_i\to e_{i+1}$ for $i\in \ZM$. Let us also use the convention that for $\ell=\infty$ we have $\ZM/\ell\ZM=\ZM$. Then we still have $\Delta=\{e_i,\,i\in \ZM/\ell\ZM\}$ for $\ell=\infty$.

	\bigskip
	
\subsection{Residues}
 
 Assume $n\in\ZM_{\geqslant 0}$. We will identify partitions with Young diagrams. 	The partition $\lambda$ corresponds to a Young diagram with $r$ lines such that the $i$th line contains $\lambda_i$ boxes. For example, the partition $(4,2,1)$ corresponds to the Young diagram
	
	$$
	\yng(4,2,1)
	$$

 Assume $\ell\in\ZM_{>0}\cup\{\infty\}$.
	We say that a box $b$ of the Young diagram is \emph{at position $(r,s)$} if it is in row $r$ and column $s$.  
	The $\ell$-\emph{residue} of the box $b$ is the 
	number $s-r$ modulo $\ell$; we say that the integer $s-r$ is the $\infty$-residue of the box $b$. Then we obtain a map 
	$$
	\Res_\ell\colon \PC\to \ZM^{\ZM/\ell\ZM},\qquad \lambda\mapsto \Res_\ell(\lambda),
	$$ 
	such that for each $i\in \ZM/\ell\ZM$ the number of boxes with $\ell$-residue $i$ in $\lambda$ is $(\Res_\ell(\lambda))_i$. 
	(In particular, we obtain a map $\Res_\infty\colon \PC\to \ZM^{\ZM}$.) For $\ell=\infty$, we mean that $\ZM^{\ZM/\ell\ZM}=\ZM^\ZM$ is the direct sum (and not the direct product) of $\ZM$ copies of $\ZM$. In other words, our convention is that for an element $\alpha=(\alpha_i)_{i\in \ZM}\in\ZM^{\ZM}$, only a finite number of integers $\alpha_i$ are non-zero.
	
	\bigskip
	
	\begin{example}
		For the partition $\lambda=(4,2,1)$ and $\ell=3$ the $3$-residues of the boxes are
		$$
		\young(0120,20,1)
		$$
		In this case $\Res_\ell(\lambda)=(3,2,2)$ because there are three boxes with residue $0$, two boxes with residue 
		$1$ and two boxes with residue $2$.
	\end{example}
	
	\bigskip
	
	We say that a box of a Young diagram is \emph{removable} if it has no boxes to its right or below it. In other words, a box $b$ is removable for $\lambda$ if $\lambda\backslash b$ is still a Young diagram. We say that a box $b$ is \emph{addable} for $\lambda$ if $b$ is not a box of $\lambda$ and $\lambda\cup b$ is still a Young diagram.  For $i\in \ZM/\ell\ZM$, we say that a box is $i$-addable or respectively $i$-removable if it is an addable or respectively removable box with $\ell$-residue $i$.

	For $\lambda,\mu\in\PC$, we write $\mu\leqslant \lambda$ if the Young diagram of $\mu$ can be obtained from 
	the Young diagram of $\lambda$ by removing a sequence of removable boxes.

	\subsection{$\ell$-cores}\label{sec:lcores}
	Assume $\ell\in \ZM_{>0}$.
	\begin{defn}
The partition $\lambda$ is an $\ell$-core if there is no partition $\mu\leqslant\lambda$ such that 
		the Young diagram of $\mu$ differs from the Young diagram of $\lambda$ by $\ell$ boxes with $\ell$ different $\ell$-residues.
	\end{defn}
	
	\bigskip
	
	See \cite{BJV} for more details about the combinatorics of $\ell$-cores. Let $\CC_\ell\subset \PC$ be the set of $\ell$-cores. 
	Set $\CC_\ell(n)=\PC(n)\cap \CC_\ell$.
	
	If a partition $\lambda$ is not an $\ell$-core, then we can get a smaller partition whose Young diagram is obtained from the Young diagram of $\lambda$ by 
	removing $\ell$ boxes with different $\ell$-residues. We can repeat this operation again and again until we get an $\ell$-core. 
	It is well-known that the $\ell$-core that we get is independent of the choice of the boxes that we remove. Then we get a function 
	$$
	\Core_\ell\colon \PC\to \CC_\ell.
	$$
	If $\mu=\Core_\ell(\lambda)$, we will say that the partition $\mu$ is {\it the $\ell$-core} of the partition $\lambda$.
	
	\bigskip
	
	\begin{example}
		The partition $(4,2,1)$ from the previous example is not a $3$-core because it is possible to remove the three bottom boxes. We get 
		$$
		\young(0120)
		$$
		But this is still not a $3$-core because we can remove three more boxes and we get 
		$$
		\young(0)
		$$
		This shows that the partition $(1)$ is the $3$-core of the partition $(4,2,1)$. 
	\end{example}
	
		\medskip
		\begin{remark}
		\label{rem:same-lcores}
	Assume that $\mu=\Core_\ell(\lambda)$ and $\mu$ is obtained from $\lambda$ by removing $r\ell$ boxes. Then $\Res_\ell(\lambda)=\Res_\ell(\mu)+r\delta$. In particular, if we have two partitions $\lambda^{(1)}$ and $\lambda^{(2)}$ with the same $\ell$-cores and such that $|\lambda^{(1)}|=|\lambda^{(2)}|$, then they have the same $\ell$-residues. More generally, if two partitions $\lambda^{(1)}$ and $\lambda^{(2)}$ have the same $\ell$-cores then we have $\Res_\ell(\lambda^{(1)})=\Res_\ell(\lambda^{(2)})+r\delta$, where $r=(|\lambda^{(1)}|-|\lambda^{(2)}|)/\ell$. 
	
	\end{remark}

	For $\nu\in \CC_\ell$, set $\PC_\nu=\{\lambda\in \PC;~\Core_\ell(\lambda)=\nu\}$ and $\PC_\nu(n)=\PC_\nu\cap\PC(n)$.
	
	\bigskip

	\subsection{Action of the affine Weyl group}\label{sub:action-weyl}
	Assume $\ell\in \ZM_{>0}$. Then $W^\aff$ is the affine Weyl group of type $\widetilde{\mathsf{A}}_{\ell-1}$. For $\ell\geqslant 2$ it is the Coxeter group with associated Coxeter system $(W^\aff,S^\aff)$, where 
	$S^\aff=\{s_i~|~i \in \ZM/\ell\ZM\}$ and the Coxeter graph has vertices the elements of $\ZM/\ell\ZM$ and we have an edge between $i$ and $i+1$ for each $i\in \ZM/\ell\ZM$. We also extend this notion to the case $\ell=1$ by setting $W^\aff=1$ in this case.  We denote by $\ell$ the length function $\ell\colon W^\aff\to \ZM_{\geqslant 0}$. 
	
	The non-affine Weyl group $W$ (isomorphic to the symmetric group $\mathfrak S_\ell$) is a parabolic subgroup of $W^\aff$ generated by $s_1,\ldots, s_{\ell-1}$ (for $\ell=1$ we mean that $W=1$).
	
	\medskip

 We consider the $W^\aff$-action on $Q\cong \ZM^{\Delta}\cong \ZM^{\ZM/\ell\ZM}$ given by 
 $s_j\star\alpha=\alpha'$, where 
	$$
    \alpha_i'=
	\begin{cases}
	\alpha_i & \text{if $i \neq j$,}\\
	1_{j,0} + \alpha_{i+1}+\alpha_{i-1}-\alpha_i & \text{if $i=j$.}\\
	\end{cases}
	$$
 The $W^\aff$-action here is the action defined in \eqref{eq:twistedaction}, where $\star=\star_{\Lambda_0}$.

 We also consider the $W^\aff$-action on $\Hom_\ZM(Q,\CM)\cong\CM^{\ZM/\ell\ZM}$
  given by 
 $$
 s_j^*(\bc)_i=
 \begin{cases}
\bc_i & \mbox{ if } i\not\in\{j-1,j,j+1\},\\
\bc_j +\bc_i & \mbox{ if } i\in\{j-1,j+1\},\\
-\bc_i & \mbox{ if } i = j.
\end{cases}
 $$
 Similarly to \S\ref{subs:graphs}, we write $w^*(\bc)$ and not $w(\bc)$ to distinguish this action from the usual action on the root lattice. This definition agrees with the definition given in \S\ref{subs:graphs}.
 In type $\widetilde{\mathsf{A}}_{\ell-1}$, the $\ZM$-linear map $Q\to \CM^{\ZM/\ell\ZM}$, $\alpha\mapsto \overline \alpha$, introduced in Section~\ref{sec:anotherpresentation} is given explicitly by
	$$
	(\overline{e_r})_i=2\cdot 1_{i,r}-1_{i,r+1}-1_{i,r-1}.
	$$
Note that this map is different from the obvious inclusion $\ZM^{\ZM/\ell\ZM}\subset\CM^{\ZM/\ell\ZM}$.
 
	\bigskip

	\bigskip

	\begin{remark}
	\label{rem:W-act-lcores}
	\cite[Section~3]{BJV} defined an $W^\aff$-action on $\CC_\ell$. Let us recall this construction. Fix $i\in \ZM/\ell\ZM$ and $\nu\in\CC_\ell$.
	\begin{itemize}
		\item[(1)] Assume that $\nu$ has neither $i$-removable boxes nor $i$-addable boxes, then we have $s_i(\nu)=\nu$.
		\item[(2)] Assume that $\nu$ has no $i$-removable boxes and has at least one $i$-addable box. Then $s_i(\nu)$ is obtained from $\nu$ by addition of all $i$-addable boxes.
		\item[(3)] Assume that $\nu$ has no $i$-addable boxes and has at least one $i$-removable box. Then $s_i(\nu)$ is obtained from $\nu$ by removing all $i$-removable boxes.
		\item[(4)] The situation when the $\ell$-core $\nu$ has an $i$-addable box and an $i$-removable box at the same time is impossible.
	\end{itemize}
    \end{remark}
    
	By construction, the map $\Res_\ell\colon \CC_\ell\to \ZM^{\ZM/\ell\ZM}$ is $W^\aff$-invariant. Moreover, the $\ell$-residue of the empty partition is zero. 
	The stabilizer of the empty partition in $W^\aff$ is $W$ and the stabilizer of $0\in\ZM^{\ZM/\ell\ZM}$ in  $W^\aff$ is also $W$. This implies that we have $W^\aff$-invariant bijections
	$$
	\begin{array}{ccccl}
	W^\aff/W&\cong& \CC_\ell&\cong& W^\aff\star 0\subset \ZM^{\ZM/\ell\ZM},\\
	wW&\mapsto& w(\emptyset)& \mapsto & w\star 0.
	\end{array}
	$$
Since $\Res_\ell$ is a $W^\aff$-invariant map and $\Res_\ell(\emptyset)=0$, the bijection $\CC_\ell\cong W^\aff\star 0$ is given by the map $\Res_\ell$. In particular, we see that an element $\ZM^{\ZM/\ell\ZM}$ is a residue of an $\ell$-core if and only if it is in the $W^\aff$-orbit of $0$.

Moreover, since we have $w\star(\alpha+n\delta)=w\star\alpha+n\delta$ and since each $W^\aff$-orbit in $\ZM^{\ZM/\ell\ZM}$ contains exactly one element of the form $n\delta$ (see \cite[Lem.~2.8]{BM}), each element $\alpha\in\ZM^{\ZM/\ell\ZM}$ has a unique presentation in the form 
\begin{equation}
\label{eq:d-via-core-n}
\alpha=\Res_\ell(\nu)+n\delta,\qquad\nu\in\CC_\ell, n\in\ZM.
\end{equation}

The following lemma is a reformulation of \cite[Remark~3.2.3]{BJV}.
\begin{lem}
	\label{lem:nu-w-i}	
	Fix $\nu\in\CC_\ell$ and $i\in \ZM/\ell\ZM$. Let $w$ be the unique element of $W^\aff$ such that $w(\emptyset)=\nu$ and such that $w$ is the shortest element in the coset $wW\in W^\aff/W$.
	The situations $(1)$, $(2)$, $(3)$ in Remark~\ref{rem:W-act-lcores} are equivalent to the following situations $(1)$, $(2)$, $(3)$ respectively:
	\begin{itemize}
		\item[$(1)$] $s_iw\in wW$ and $\ell(s_iw)>\ell(w)$,
		\item[$(2)$] $s_iw\not\in wW$ and $\ell(s_iw)>\ell(w)$,
		\item[$(3)$] $s_iw\not\in wW$ and $\ell(s_iw)<\ell(w)$.
	\end{itemize}
	
\end{lem}

	\bigskip

\subsection{$J$-cores}\label{sec:Jcore} Fix a subset $J\subset \ZM/\ell\ZM$. 

\begin{defn}
	We say that a box of a Young diagram is \emph{$J$-removable} if it is removable and its residue is in $J$. We say that a Young diagram is a \emph{$J$-core} if it has no $J$-removable boxes. Denote by $\CC_J$ the set of all $J$-cores. 
	
	To each partition $\lambda\in \PC$ we can associate a partition $\Core_J(\lambda)\in \CC_J$ obtained from it by removing $J$-removable boxes (probably in several steps). The result $\Core_J(\lambda)$ does not depend on the order of operations.
\end{defn}

\medskip
\begin{lem}
\label{lem:lcore-of-Jcore}
	For each $\mu\in \CC_J$, we have $\Core_\ell(\mu)\in \CC_J$.
\end{lem}
\begin{proof}
This statement is quite obvious when we see the partition $\mu$ as an abacus, see for example \cite[\S{2}]{BJV} for then definition of an abacus.

However we can give another proof based on the representation theory of quivers and the results of Section~\ref{sec:quiver-cycl}. Fix some $J$-standard $\bc\in \CM^{\ZM/\ell\ZM}$. Since $\nu$ is a $J$-core, the representation $A_\mu$ constructed in Section~\ref{subs:C-fixed-points} is simple by Lemma~\ref{lem:structure-Amu}. Then the dimension vector $\Res_\ell(\mu)$ of this representation is in $E_\bc$. 

Now, let $\nu$ be the $\ell$-core of $\mu$. Assume that $\nu$ is obtained from $\mu$ by removing $r\ell$ boxes. Then $\Res_\ell(\mu)=\Res_\ell(\nu)+r\delta\in E_\bc$ and hence Lemma~\ref{lem:descr-Eth} implies that $\nu$ is a $J$-core.

\end{proof}

\section{Quiver varieties for the cyclic quiver}\label{sec:quiver-cycl}
\medskip

\subsection{Quiver varieties for the cyclic quiver}\label{sub:quiver-cycl}

Assume $\ell\in \ZM_{>0}\cup\{\infty\}$.  
Consider a dimension vector $\alpha=(\alpha_i)_{i \in \ZM/\ell\ZM}\in \ZM_{\geqslant 0}^{\ZM/\ell\ZM}\cong Q^+$ for the quiver $\Qu^\ell$. For $\ell=\infty$, we always assume additionally that $\alpha\in \ZM^\ZM$ has a finite number of non-zero components.

Let $\Qov^\ell$ be the double quiver of $\Qu^\ell$. That is, $\Qov^\ell$ is the quiver obtained from $\Qu^\ell$ by adding an opposite arrow to each arrow of $\Qu^\ell$. We would also like to have a framed version adding a $1$-dimensional framing only for the vertex $e_0$. Similar to the notation given in \S\ref{rootsystemnotation}, let $\Qulinf$ be the quiver obtained from $\Qu^\ell$ by adding an extra vertex $e_\infty$ and an extra arrow $e_\infty\to e_0$. Denote by $\oQulinf$ the double quiver of $\Qulinf$. For each dimension vector $\alpha$ as above for the quiver $\Qu^\ell$, consider the dimension vector $\widehat\alpha=\alpha+e_\infty$ for the quiver $\Qov^\ell$ (we just add the dimension $1$ component for the extra vertex). Let us consider the quiver variety $\XC_\bc(\alpha)$ as in \S\ref{sec:quivernotation}. Since we wish to define an action of $\Cs$ on $\XC_\bc(\alpha)$, we spell out in greater detail the definition of $\XC_\bc(\alpha)$ in this particular case.

An element of $\Rep(\oQulinf,\widehat\alpha)$ is a tuple $(X,Y,x,y)$, where 
$$
X=(X_i)_{i\in \ZM/\ell\ZM},\quad X_i\in \Hom(\CM^{\alpha_{i+1}}, \CM^{\alpha_{i}}) \qquad Y=(Y_i)_{i\in \ZM/\ell\ZM},\quad Y_i\in \Hom(\CM^{\alpha_i}, \CM^{\alpha_{i+1}}),
$$
$$ 
x\in \Hom(\CM,\CM^{\alpha_0}),\qquad y\in \Hom(\CM^{\alpha_0},\CM).
$$

The group $G(\alpha)$ acts on $\Rep(\oQulinf,\widehat\alpha)$. 
We consider
$$\fonction{\widehat\mu_\alpha}{\Rep(\oQulinf,\widehat\alpha)}{\bigoplus_{i \in \ZM/\ell\ZM} 
\End(\CM^{\alpha_i})}{(X_i,Y_i,x,y)_{i \in \ZM/\ell\ZM}}{(X_iY_i - Y_{i-1}X_{i-1}+1_{i,0}xy)_{i \in \ZM/\ell\ZM}}$$
the corresponding {\it moment map}. 
If $\bc=(\bc_i)_{i \in \ZM/\ell\ZM}\in\CM^{\ZM/\ell\ZM}$, we denote by $\Irm_\bc(\alpha)$ the family $(\bc_i \Id_{\CM^{\alpha_i}})_{i \in \ZM/\ell\ZM}$.
Finally, we set
$$\YC_\bc(\alpha)=\widehat\mu_\alpha^{-1}(\Irm_\bc(\alpha)).$$
We get the following description of the variety: 
$$\XC_\bc(\alpha)=\YC_\bc(\alpha)/\!\!/G(\alpha).$$ 

\begin{remark}\label{rem:dim-negative}
We extend the definition of $\XC_\bc(\alpha)$ to the case where $\alpha \in \ZM^{\ZM/\ell\ZM}$ 
by the convention that $\XC_\bc(\alpha)=\vide$ whenever at least one of the $\alpha_i$ is negative. 
\end{remark}

Let $\Rep(\Qu)$ be the category of (finite dimensional) representations of the quiver $\Qu$. We can view each element of $\Rep(\Qu,\alpha)$ as an object in $\Rep(\Qu)$ with dimension vector $\alpha$.

Now, assume $\ell\in \ZM_{>0}$. 

\begin{defn}
\label{def:map-i-inf}
Consider the following map $\iota\colon \Rep(\oQuii)\to  \Rep(\oQulinf)$.

For each finite dimensional representation $(X,Y,x,y)$ of $\oQuii$ with the underlying vector space $V=\bigoplus_{j\in \ZM} V_j$ we can associate a representation $(X',Y',x',y')$ of $\oQulinf$ with the underlying vector space $V'=\bigoplus_{i\in \ZM/\ell\ZM} V'_i$ where
$$
V'_i = \bigoplus_{\substack{j \in \ZM \\ j \equiv i \smod l}} V_j,
\qquad 
X'_i = \bigoplus_{\substack{j \in \ZM \\ j \equiv i \smod l}} X_j,
\quad  
Y'_i = \bigoplus_{\substack{j \in \ZM \\ j \equiv i \smod l}} Y_j,
$$
$x'$ is the composition of $x$ with the natural map $V_0\to V'_0$, $y'$ is the composition of $y$ with the natural map $V'_0\to V_0$.
\end{defn}

\subsection{Reflection isomorphism}
\label{subs:Lusztig-iso}

By Section~\ref{subs:adm-refl}, we have an isomorphism
\equat\label{eq:iso-lusztig}
\XC_{s_j^*(\bc)}(s_j\star\alpha) \cong \XC_\bc(\alpha)\qquad \mbox{ if }\bc_j\ne 0.
\endequat
Note that this isomorphism takes into account the convention of Remark~\ref{rem:dim-negative}.

The isomorphism above motivates one to consider the following equivalence relation on the set $\ZM^{\ZM/\ell\ZM} \times \CM^{\ZM/\ell\ZM}$. Let $\sim$ be the transitive closure of 
$$
(\alpha,\bc)\sim (s_i\star\alpha,s_i^*(\bc)) \qquad \mbox{ if } \bc_i\ne 0.
$$
The isomorphism \eqref{eq:iso-lusztig} implies that if $(\alpha,\bc)\sim(\alpha',\bc')$ then we have an isomorphism of algebraic varieties $\XC_\bc(\alpha)\cong \XC_{\bc'}(\alpha')$.

\begin{remark}
\label{rem:equiv=shortest}
Let $W^\aff_\bc$ denote the stabilizer of $\bc$ in $W^\aff$. Assume that $\bc$ is such that $W^\aff_\bc$ is a parabolic subgroup of $W^\aff$. Then we can describe the set of pairs that are equivalent to $(\alpha,\bc)$ in the following way. They are of the form $(w\star\alpha,w^*(\bc))$, where $w$ is the element of shortest length in the right coset $wW^\aff_\bc\in W^\aff/W^\aff_\bc$.
\end{remark}

\subsection{Quiver varieties vs Calogero--Moser varieties}\label{sub:quiver-cm}

Assume $\ell\in\ZM_{>0}$. Recall that we assumed $\Gamma=\ZM/\ell\ZM$. Let us review the isomorphism in Theorem~\ref{thm:CMquiveriso} in this case. We include the additional $\CM^\times$-action.

Assume that $n \ge 2$.
We set $\zeta=\exp(2\pi \sqrt{-1} /\ell)$. We denote by $s$ the permutation matrix corresponding to the 
transposition $(1,2)$ and we set
$$t=\diag(\zeta,1,\dots,1) \in \Gamma_n.$$
Then $s$, $t$, $t^2$,\dots, $t^{\ell-1}$ is a set of representatives of conjugacy 
classes of reflections in $\Gamma_n$. For simplicity, we set 
$$
a=c_s\qquad\text{and}\qquad k_j=\frac{1}{\ell} \sum_{i=1}^{\ell-1} \zeta^{-i(j-1)} c_{t^i}
$$
for $j \in \ZM/\ell\ZM$. Then
\equat\label{eq:k}
k_0+\cdots + k_{\ell-1} = 0\qquad \text{and}\qquad c_{t^i}=\sum_{j \in \ZM/\ell\ZM} \zeta^{i(j-1)} k_j
\endequat
for $1 \le i \le \ell-1$. Finally, if $i \in \ZM/\ell\ZM$, we set

\equat
\bc_i=
\begin{cases}
k_{-i}-k_{1-i} & \text{if $i \neq 0$,}\\
-a +k_0-k_1 & \text{if $i=0$.}\\
\end{cases}
\endequat\label{eq:theta-k}
and $\bc=(\bc_i)_{i \in \ZM/\ell\ZM}$.

\bigskip

There is a $\CM^\times$-action on $\XC_\bc(n\delta)$ given by $\xi\cdot (X,Y,x,y)=(\xi^{-1}X,\xi Y,x,y)$. 
The following proposition is a $\CM^\times$-equivariant version of Theorem~\ref{thm:CMquiveriso} in this situation, see also Remark~\ref{rem:C*-on-CM}. However, the choice of the parameter $\bc$ in terms of $c$ here is different from the choice made in Theorem~\ref{thm:CMquiveriso} by multiplication by a constant. This choice is made to be compatible with \cite{BM}.

The following result is proved in~\cite[Theorem~3.10]{GordonQuiver}. (Note that our $k_i$ is related to Gordon's $H_i$ via $H_i=k_{-i}-k_{1-i}$.) 

\begin{prop}\label{prop:quiver-cm}
There is a $\CM^\times$-equivariant 
isomorphism of varieties
$$\ZCB_{\! c} \longiso \XC_\bc(n\delta).$$
\end{prop}

In the above isomorphism, the parameter $a$ of the variety $\ZCB_{\! c}$ corresponds to $~-(\sum_{i\in \ZM/\ell\ZM}\bc_i)=- \bc(\delta)$ for $\XC_\bc(n\delta)$. 
Note that $a=-\bc(\delta)$ is invariant under the transformation of the parameter $\bc\mapsto s_j^*(\bc)$. From now on, we assume $a\ne 0$. 

\begin{remark}
	All statements in Section~\ref{sub:quiver-cm} also make sense for $n=1$ with the following modifications. We have no transposition $s$, so we have no parameter $a$. On the other hand, for $n=1$, the variety $\XC_\bc(n\delta)$ does not depend on $\bc_0$. 
 	
	We can also use the convention that for $n=0$ the Calogero--Moser variety is a point. Then Proposition~\ref{prop:quiver-cm} still holds.
\end{remark}

Recall also from~\cite[\S{11}]{EG} the following result, which follows from Proposition~\ref{prop:quiver-cm}.

\bigskip

\begin{lem}\label{lem:dim-CM}
	If $n \ge 0$, then $\XC_\bc(n\delta)$ is normal of dimension $2n$.
\end{lem}

\subsection{Semisimple representations in $\Rep_{\bc}(\Qov^\ell)$}

Denote by $\Rep_{\bc}(\Qov^\ell)$ the additive category of representations $(X,Y)$ of $\Qov^\ell$ satisfying the moment map relations $\widehat\mu_\alpha(X,Y)=\Irm_\bc(\alpha)$, where $\alpha$ is the dimension vector of the representation $(X,Y)$. 
Since we assume $a\ne 0$, we are in the setup of Section~\ref{sec:nonzero}. Then the set $\Sigma_\bc$ of dimension vectors of simple representations in $\Rep_{\bc}(\Qov^\ell)$ is the same as the set of simple roots $\Delta(\bc)$ in $R_\bc$. The set of dimension vectors of semisimple representations is the set $\SS{\bc}$ of (possibly zero) sums of the elements of $\Sigma_\bc$. Moreover, since each element of $\SS{\bc}$ has a unique decomposition as a sum of elements of $\Sigma_\bc=\Delta(\bc)$, for each $\alpha\in \SS{\bc}$ there exists a unique up to isomorphism semisimple representation in $\Rep_\bc(\Qov^\ell)$. Let us denote this representation by $L(\alpha)$. We see in particular that for each $\alpha\in \ZM_{\geqslant 0}^{\ZM/\ell\ZM}$, the variety $\mf{M}_\bc(\alpha)$ is either a singleton or empty. More precisely, we have 
$$
\mf{M}_\bc(\alpha)
=
\begin{cases}
\{L(\alpha)\}&\mbox{ if $\alpha\in \SS{\bc}$},\\
\emptyset &\mbox{ else}.
\end{cases}
$$

\subsection{Symplectic leaves}
\label{subs:symp-leaves}

Denote by $\Rep_{\bc}(\oQulinf)$ the category of representations $(X,Y,x,y)$ of $\oQulinf$ whose dimension vector is of the form $\widehat\alpha$ for some $\alpha\in \ZM^{\ZM/\ell\ZM}$ and satisfying the moment map relations $\widehat\mu_\alpha(X,Y)=\Irm_\bc(\alpha)$. This category is not additive because we have imposed that the representations have dimension $1$ at the vertex $\infty$. However, it does make sense to add an object of $\Rep_{\bc}(\Qov^\ell)$ and an object of $\Rep_{\bc}(\oQulinf)$, getting an object of $\Rep_{\bc}(\oQulinf)$.

An object $M$ of $\Rep_{\bc}(\oQulinf)$ is indecomposable as a representation of the quiver $\oQulinf$ if and only if the only possible decomposition $M=M_0\oplus M_1$ with $M_0\in \Rep_{\bc}(\oQulinf)$ and $M_1\in \Rep_{\bc}(\Qov^\ell)$ is $M=M\oplus 0$.

\begin{remark}
\label{rk:couple-equiv-nd}
Assume $\alpha\in E_\bc$. Then, by Proposition~\ref{prop:Cadmissiblemove}, the pair $(\alpha,\bc)$ is equivalent to a pair of the form $(n\delta,\bc')$ with $n\geqslant 0$. In particular, by Proposition~\ref{prop:quiver-cm}, the variety $\XC_\bc(\alpha)$ is isomorphic to the Calogero--Moser variety.
\end{remark}

Each object $M\in\Rep_{\bc}(\Qov^\ell)$ has a unique decomposition $M=M_0\oplus M_1$ such that $M_0\in \Rep_{\bc}(\oQulinf)$, $M_1\in \Rep_{\bc}(\Qov^\ell)$ and $M_0$ is indecomposable. Set $\dim^{\rm reg}M=\dim M_0\in \ZM^{\ZM/\ell\ZM}$.

Take a point $[M]\in\XC_\bc(\alpha)$ presented by a semisimple representation $M\in \Rep_{\bc}(\oQulinf)$.

\begin{lem}
\label{lem:symple-leaves-descr}
Two points of $[M], [M']\in \XC_\bc(\alpha)$ are in the same symplectic leaf if and only if we have $\rdim(M)=\rdim(M')$.
\end{lem}
\begin{proof}
Let us decompose $M$ in a direct sum of simple representations $M=\bigoplus_{r=0}^k M_r$, where $M_0\in \Rep_{\bc}(\oQulinf)$ and other summands are in $\Rep_{\bc}(\Qov^\ell)$.

Once we know the dimension vector $\alpha'$ of $M_0$, we know automatically $k$ and the dimension vectors of $M_1, M_2, \ldots,M_k$ (up to a permutation) because $L(\alpha-\alpha')$ is the unique semisimple representation in $\Rep_{\bc}(\Qov^\ell)$ of dimension vector $\alpha-\alpha'$. Then the statement follows from the description of symplectic leaves given in \cite[Theorem~1.9]{BellSchedQuiver}. 
\end{proof}

For two dimension vectors $\alpha$ and $\alpha'$ we set $\LG^\alpha_{\alpha'}=\{[M]\in \XC_\bc(\alpha);~\rdim(M)=\alpha'\}$. By Lemma~\ref{lem:symple-leaves-descr} $\LG^\alpha_{\alpha'}$ is either a symplectic leaf of $\XC_\bc(\alpha)$ or is empty. Note that the labeling used here is different from the labeling of symplectic leaves in Section~\ref{sec:nonzero}. The leaf $\LG^\alpha_{\alpha'}$ here corresponds to $\LG(\alpha-\alpha')$ in Section~\ref{sec:nonzero}.

\begin{lem}
\label{lem:order-sym-leaves}
The symplectic leaves $\LG^\alpha_{\alpha'}\subset \XC_\bc(\alpha)$ define a finite stratification of $\XC_\bc(\alpha)$ into locally closed subsets.
For two symplectic leaves $\LG^\alpha_{\alpha'}$ and $\LG^\alpha_{\alpha''}$ of $\XC_\bc(\alpha)$ we have $\LG^\alpha_{\alpha'}\subset \overline{\LG^\alpha_{\alpha''}}$ if and only if $\alpha''-\alpha'\in \SS{\bc}$.
\end{lem}
\begin{proof}
This statement is a special case of \cite[Proposition~3.6]{BellSchedQuiver}. See also the proof of Proposition~\ref{prop:leafclosure1} for more details.

\end{proof}

\begin{prop}
\label{prop:can-decomp-of-d}
 \mbox{ }
 $(i)$ For each dimension vector $\alpha$ such that $\XC_\bc(\alpha)\ne \emptyset$, there is a decomposition $\alpha=\alpha^0+\alpha^1$ such that $\alpha^0\in E_\bc$ and $\alpha^1\in\SS{\bc}$ such that for any other decomposition $\alpha=\alpha'^0+\alpha'^1$ with $\alpha'^0\in E_\bc$ and $\alpha'^1\in\SS{\bc}$ we have $\alpha^0-\alpha'^0\in \SS{\bc}$. 

$(ii)$ $\LG_{\alpha^0}$ is the unique open symplectic leaf in $\XC_\bc(\alpha)$.

$(iii)$ We have an isomorphism of varieties 
$$
\XC_\bc(\alpha^0)\cong \XC_\bc(\alpha),\qquad [M]\mapsto [M\oplus L(\alpha^1)].
$$
\end{prop}

\begin{proof}
	
By \cite[Corollary~1.45]{MartinoThesis}, the smooth locus of $\XC_\bc(\alpha)$ is a symplectic leaf. Then it should be of the form $\LG^\alpha_{\alpha^0}$ for some $\alpha^0$.
	 
Since $\XC_\bc(\alpha)$ is irreducible by \cite[Corollary~1.4]{CBdecomp}, we have $\overline{\LG^\alpha_{\alpha^0}}=\XC_\bc(\alpha)$.
Then, by Lemma~\ref{lem:order-sym-leaves}  for any other symplectic leaf $\LG^\alpha_{\alpha'_0}$ we have $\alpha^0-\alpha'^0\in\SS{\bc}$. This proves $(i)$ and $(ii)$.

Part $(iii)$ follows from \cite[Theorem~1.1]{CBdecomp}.
\end{proof}

Now, we set $\XC_\bc(\alpha)^{\rm reg}=\LG^\alpha_{\alpha^0}$. 
Assume that $\alpha$ and $\alpha'$ are such that $\LG^\alpha_{\alpha'}$ is non-empty.

\begin{lem}
The closure of $\LG^\alpha_{\alpha'}$ is isomorphic to $\XC_{\bc}(\alpha')$. 
\end{lem}

\begin{proof}
Lemma~\ref{lem:normalization-leaf-cl} and Proposition~\ref{prop:leafclosurenormal} imply that the map
$$
\phi\colon\XC_{\bc}(\alpha')\to \overline{\LG^\alpha_{\alpha'}},\qquad [M]\mapsto [M\oplus L(\alpha-\alpha')]
$$
is an isomorphism.

\end{proof}

\begin{cor}
\label{coro:norm-cl-leaf}
The closure of each symplectic leaf $\LG^\alpha_{\alpha'}$ of the variety $\XC_{\bc}(\alpha)$ is isomorphic to a variety of the form $\XC_{\bc'}(r\delta)$ for some $r\geqslant 0$ and some $\bc'\in \CM^{\ZM/\ell\ZM}$. 
\end{cor}
\begin{proof}

First of all, note that $\alpha'\in E_\bc$. By Remark~\ref{rk:couple-equiv-nd}, the pair $(\alpha',\bc)$ is equivalent to some pair of the form $(r\delta,\bc')$ where $r\geqslant 0$ and $\bc'\in \CM^{\ZM/\ell\ZM}$. Then the isomorphism \eqref{eq:iso-lusztig} yields $\XC_{\bc}(\alpha')\cong \XC_{\bc'}(r\delta)$.
\end{proof}

Combining the corollary above with Proposition~\ref{prop:quiver-cm} yields the following theorem.

\begin{thm}\label{thm:GM1nclosure}
The closure of each symplectic leaf of the Calogero--Moser variety of type $G(\ell,1,n)$ with $a\ne 0$ is isomorphic to a Calogero--Moser variety of type $G(\ell,1,r)$ for some $r\in [0;n]$. In particular, all leaf closures are normal when $a \neq0$. 
\end{thm}

\begin{remark}
\label{rem:param-c'-cyclic}
We explain the relationship between the parameters of the two Calogero--Moser varieties in Theorem~\ref{thm:GM1nclosure}. The Calogero--Moser variety $\ZCB_{c}(G(\ell,1,n))$ is isomorphic to the quiver variety $\XC_\bc(n\delta)$.  The closure of the symplectic leaf $\LG^{n\delta}_{\alpha'}$ is isomorphic to  $\XC_\bc(\alpha')$. Just as in the proof of Proposition~\ref{prop:wbccomp}, we can find $w\in W^\aff$ that realizes the equivalence between $(\alpha',\bc)$ and $(w\star\alpha' = r \delta,w^*(\bc))$. Set $\bc'=w^*(\bc)$. We have an isomorphism $\XC_\bc(\alpha')\cong \XC_{\bc'}(r\delta)$. Since $w\star\alpha'=r\delta$, Lemma~\ref{lem:affWeyl_translation-roots} implies that $w=x t_{\alpha'}$, where $x\in W$. Then, by Lemma~\ref{lem:affWeyl_translation-theta}, $\bc'=(xt_{\alpha'})^*(\bc)=x^*(\bc-a\overline{\alpha'})$. Moreover, the action the element $x\in W$ on $\CM^{\ZM/\ell\ZM}$ corresponds to some permutation of the parameters $k_0,k_1,\ldots,k_{\ell-1}$ (see \cite[Remark~3.5]{BM}) and a permutation of the parameters does not change the Calogero--Moser variety up to isomorphism, see Remark~\ref{rem:finiteWisoparameters}.

Therefore, the parameters $a, k_0,k_1,\ldots,k_{\ell-1}$ (corresponding to $\bc$) of the original the Calogero--Moser variety $\ZCB_{c}(G(\ell,1,n))$ are related to the parameters $a', k'_0,k'_1,\ldots,k'_{\ell-1}$ (corresponding to $\bc'$) of the new Calogero--Moser variety $\overline{\LG^{n\delta}_{\alpha'}} \cong \ZCB_{c'}(G(\ell,1,r))$ as follows:
$$
a'=a,\qquad k'_i=k_i+(\alpha'_{1-i}-\alpha'_{-i}).
$$
In the case where $n$, resp $r$, is equal to $1$, we can forget the parameter $a$, resp. $a'$. In the case $r=0$, the variety $\XC_{\bc'}(r\delta)$ is just a point.
\end{remark}

\subsection{$\CM^\times$-fixed points}
\label{subs:C-fixed-points}

For each $J\subset \ZM/\ell\ZM$ we denote by $W_J$ the parabolic subgroup of $W^\aff$ generated by $s_i$ for $i\in J$. Let us say that $\bc$ is \emph{$J$-standard} if the stabilizer $W^\aff_\bc$ of $\bc$ in $W^\aff$ is equal to $W_J$. We say that $\bc\in\CM^{\ZM/\ell\ZM}$ is \emph{standard} it is $J$-standard for some $J\subset \ZM/\ell\ZM$. For a standard $\bc$, the set $J$ is the set of indices $i\in \ZM/\ell\ZM$ such that $\bc_i=0$.

Now, let us describe the $\CM^\times$-fixed points of $\XC_\bc(\alpha)$. First of all, each pair $(\alpha,\bc)$ is equivalent to a pair whose $\bc$ is standard by Lemma~\ref{lem:wparabolicrootsystem}.

The following lemma is obvious.
	\begin{lem}
			\label{lem:Sigma-th-simple}
		Assume that $\bc$ is $J$-standard. Then we have $\Sigma_\bc=\{e_i;~i\in J\}$.
	\end{lem}
	Let us now assume that $\bc$ is $J$-standard. 
For each partition $\mu$, we construct a $\CM^\times$-fixed point in $\XC_{\bc}(\Res_\ell(\mu))$. This construction is essentially the same as \cite[Section~5]{Przez}. However \cite{Przez} assumes that the variety $\XC_{\bc}(\Res_\ell(\mu))$ is smooth and we don't need this assumption. 

We are going to use the Frobenius forms of partitions: each partition $\mu$ can be described by some $k\in\ZM_{\geqslant 0}$ and $a_1,\ldots, a_k, b_1,\ldots,b_k\in \ZM_{\geqslant 0}$, where $k$ is maximal such that the Young diagram of $\mu$ contains a box in position $(k,k)$ and for each $r\in [1;k]$ there are $a_r$ boxes on the right of $(r,r)$ and $b_r$ boxes below $(r,r)$.
In other words, we see the Young diagram of the partition $\mu$ as a union of $k$ hooks. The box at position $(i,j)$ is in the $r$th hook if $\min(i,j)=r$. The numbers $a_r$ and $b_r$ are the lengths of the arm and of the leg of $r$th hook respectively. 

For $i\in \ZM$, we use the convention that $\bc_i$ means $\bc_{(i\smod l)}$. Set $\beta_r=\sum_{i=-b_r}^{a_r}\bc_i$. 

Let $V$ be a complex vector space with basis $\{v_{r,j};~r\in[1;k];~j\in[-b_r,a_r]\}$. It has a $\ZM$-grading $V=\bigoplus_{j\in \ZM} V_j$ such that $v_{r,j}\in V_j$. Consider two endomorphisms $X$ and $Y$ of this vector space given by
$$
X(v_{r,j})=
\begin{cases}
v_{r,j-1} &\mbox{ if } j>-b_r,\\
0 &\mbox{ if } j=-b_r,
\end{cases}
$$
and
$$
Y(v_{r,j})=
\begin{cases}
(\sum_{i=-b_r}^j \bc_{i})v_{r,j+1}+\sum_{t>r} \beta_tv_{t,j+1} &\mbox{ if } j\in[-b_r,-1]\\
-(\sum_{i=j+1}^{a_r}\bc_{i})v_{r,j+1}-\sum_{t<r}\beta_tv_{t,j-1} &\mbox{ if } j\in [0;a_r-1],\\
-\sum_{t<r}\beta_tv_{t,j-1} &\mbox{ if } j=a_r,\\
\end{cases}
$$

Consider also the linear maps $x\colon \CM\to V_0$ and $y\colon V_0\to \CM$ given by 
$$
x(1)=-\sum_{r=1}^k \beta_r v_{r,0}\qquad  \mbox{ and } \qquad y(v_{r,0})=1.
$$

Then $(X,Y,x,y)$ yields a representation $A^\infty_\mu$ of the quiver $\oQuii$. Applying the map $\iota$ as in Definition~\ref{def:map-i-inf}, we get a representation $A_\mu$ of the quiver $\oQulinf$. It satisfies the moment map relation $\widehat\mu_\alpha(A_\mu)=I_\bc(\alpha)$.

\begin{lem}
\label{lem:structure-Amu}

Assume that $\bc$ is $J$-standard.

\begin{enumerate}
\item[(i)] If $\mu$ is a $J$-core, then $A_\mu$ is simple.

\item[(ii)] Assume that  $b$ is a removable box of $\mu$ with $\ell$-residue $i\in J$. Then we have either a short exact sequence
$$
0\to  L(e_i)\to A_{\mu}\to A_{\mu\backslash b}\to 0
$$
or a short exact sequence
$$
0\to A_{\mu\backslash b}\to A_{\mu}\to  L(e_i)\to 0.
$$
\end{enumerate}

\end{lem}
\begin{proof}
First, we prove (ii). Assume that $b$ is the box as in the statement.  Assume that it is in the $r$th hook. Let $j$ be the $\infty$-residue of $b$.

Assume first $j<0$. We have $X(v_{r,j})=Y(v_{r,j})=0$.  Then the vector $v_{r,j}$ spans a subrepresentation isomorphic to $L(e_i)$. We get a short exact sequence
$$
0\to  L(e_j)\to A_{\mu}\to A_{\mu\backslash b}\to 0.
$$

Now, assume $j\geqslant 0$. Then we see that $A_{\mu\backslash b}$ is a subrepresentation of $A_{\mu}$. It is spanned by all basis vectors except $v_{r,j}$.
Then we have a short exact sequence
$$
0\to A_{\mu\backslash b}\to A_{\mu}\to  L(e_j)\to 0.
$$

\medskip
Now, let us prove (i). First of all, we note that the assumption that $\bc$ is $J$-standard implies that if for some $a,b\in \ZM$, $a\leqslant b$ we have $\bc_a+\bc_{a+1}+\ldots+\bc_{b-1}+\bc_b=0$, then we have $\bc_a=\bc_{a+1}=\ldots=\bc_{b-1}=\bc_b=0$. 
If $\mu$ is a $J$-core, then the numbers $\beta_1, \beta_2,\ldots,\beta_k$ are non-zero. Indeed, if some $\beta_r$ is zero, then $\beta_k$ is also zero. Then the $\ell$-residues of all boxes of the $k$th hook are in $J$. In particular, the $k$th hook contains a removable box whose residue is in $J$. This contradicts the fact that $\mu$ is a $J$-core.

In view of Lemma~\ref{lem:Sigma-th-simple}, if the representation $A_\mu$ is not simple, then it must either contain a subrepresentation of the form $L(e_i)$, or it must have a quotient of the form $L(e_i)$. Let us show that both situations are impossible when $\mu$ is a $J$-core.

Assume that $A_\mu$ has a subrepresentation isomorphic to $L(e_i)$. Let $v$ be a vector that spans this subrepresentation. We can write $v=\sum\limits_{j\in \ZM, j\equiv i \smod l}v_j$, where $v_j\in V_j$. Take $j$ in this decomposition such that $v_j\ne 0$. Then the vector $v_j$ also spans a subrepresentation of $A_\mu$ isomorphic to $L(e_i)$.

Let $t$ be the number of boxes of $\mu$ with the $\infty$-residue $j$. Write $v_j=\sum_{r=1}^t\lambda_r v_{r,j}$. Then $X(v)=0$ is only possible when $\lambda_1=\ldots=\lambda_{t-1}=0$, so the vector $v_{t,j}$ spans $L(e_i)$.

Assume $j<0$.  Since the box $b$ corresponding to the vector $v_{t,j}$ cannot be removable, the diagram of $\mu$ either contains the box below $b$ or the box on the right of $b$. In the first case we must have $X(v_{t,j})\ne 0$ and in the second case we must have $Y(v_{t,j})\ne 0$. This is a contradiction.

Assume $j>0$. Then $X(v_{t,j})\ne 0$. This is a contradiction.

Assume $j=0$. Then, since $\beta_1\ne 0$, $Y(v_{t,0})\ne 0$ is only possible for $t=1$. However, this implies that $\mu$ contains only one hook (i.e., we have $k=1$). Since the box $b$ corresponding to the vector $v_{1,0}$ cannot be removable, the diagram of $\mu$ either contains the box below $b$ or the box on the right of $b$. The first case is not possible because it implies $X(v_{1,0})\ne 0$. In the second case we must have $\bc_1+\bc_2+\ldots+\bc_{a_1}=0$. However, this implies ${a_1}=0$ and then the unique box with $\infty$-residue $a_1$ is removable. This is a contradiction.

Now, assume that $A_\mu$ has a quotient isomorphic to $L(e_i)$. Then the dual representation $A_\mu^*$ contains a submodule isomorphic to $L(e_i)$. An argument as above show that this is impossible if $A_\mu$ is a $J$-core.

\end{proof}

Denote by $A'_\mu$ the semisimplification of $A_\mu$, i.e., $A'_\mu$ is the direct sum of the Jordan-H\"older subquotients of $A_\mu$.

\begin{cor} 
\label{coro:structure-Amu}
	Assume $\mu\in \PC$ and set $\lambda=\Core_J(\mu)$. Then the representation $A'_\mu$ has the following decomposition in a direct sum of simple representations
	$$
	A'_\mu=A_\lambda\oplus \bigoplus_j L(e_j),
	$$
where the sum is taken by the multiset of $\ell$-residues of $\mu\smallsetminus\lambda$.
\end{cor}

\begin{defn}
	We say that the representation $(X,Y,x,y)$ of $\oQulinf$ is $\ZM$-\emph{gradable} if it is isomorphic to the image by $\iota$ (see Definition~\ref{def:map-i-inf}) of some representation $L$ of $\oQuii$. In this case we say that $L$ is a \emph{graded lift} of $(X,Y,x,y)$.
\end{defn}

A  $\ZM$-gradable representation yields a $\CM^\times$-fixed point in $\XC_\bc(\alpha)$.

\begin{lem}
Assume that $(X,Y,x,y)$ is simple and $\ZM$-gradable. Then its $\ZM$-grading is unique. 
\end{lem} 
\begin{proof}
Since we assume $a\ne 0$, the vector $v=x(1)$ must be non-zero (here $1$ is a vector spanning the $\infty$-component of the representation, which is isomorphic to $\CM$). Then $v$ should be in $\ZM$-degree $0$. Since the representation is simple, the vectors of the form $X^{a_1}Y^{b_1}\ldots X^{a_k}Y^{b_k}(v)$ and the vector $1$ span the representation. But then vector $X^{a_1}Y^{b_1}\ldots X^{a_k}Y^{b_k}(v)$ must be in $\ZM$-degree $b_1-a_1+\ldots+b_k-a_k$. This shows that the $\ZM$-grading is unique. 
\end{proof}

\begin{example}
\label{ex:Res-Amu}
If $\mu$ is a $J$-core, then the representation $A_\mu$ is simple. It is $\ZM$-gradable by construction. Its graded lift $A^\infty_\mu$ is unique. The $\ZM$-graded dimension of the graded lift $A^\infty_\mu$ is $\Res_\infty(\mu)$.
\end{example}

\begin{cor} 
\label{coro:isom-if-same-Jcore}
For $\mu,\eta \in \PC_\nu(n\ell+|\nu|)$, the representations $A'_{\mu}$ and $A'_{\eta}$ are isomorphic if and only if $\mu$ and $\eta$ have the same $J$-cores.
\end{cor}

\begin{proof}
Let $\lambda$ and $\rho$ be the $J$-cores of $\mu$ and $\eta$ respectively.

Assume that $A'_{\mu}$ and $A'_{\eta}$ are isomorphic. We see from Corollary~\ref{coro:structure-Amu} that the representations $A_{\lambda}$ and $A_{\rho}$ are also isomorphic. Now, Example~\ref{ex:Res-Amu} implies $\Res_\infty(\lambda)=\Res_\infty(\rho)$, this yields $\lambda=\rho$. 

Now, assume that we have $\lambda=\rho$. Since we have $\mu,\eta\in \PC_\nu(n\ell+|\nu|)$, the partitions $\mu$ and $\eta$ have the same residues equal to $\Res_\ell(\nu)+n\delta$. Then $\mu\smallsetminus\lambda$ and $\eta \smallsetminus\rho$ have the same residues. Then Corollary~\ref{coro:structure-Amu} implies that $A'_{\mu}$ and $A'_{\eta}$ are isomorphic.
\end{proof}

\begin{remark}
For each partition $\mu$, we have a $\CM^\times$-fixed point $[A'_\mu]\in \XC_\bc(\Res_\ell(\mu))$ presented by the representation $A'_\mu$. 
 Assume $\alpha\in E_\bc$, see Remark~\ref{rk:couple-equiv-nd}. Write $\alpha=\Res_\ell(\nu)+n\delta$, where $\nu\in\CC_\ell$. By \cite[Proposition~8.3~(i)]{GordonQuiver}, the $\CM^\times$-fixed points in $\XC_\bc(\alpha)$ are parameterized by $J$-cores of elements of $\PC_\nu(n\ell+|\nu|)$. On the other hand, we have already constructed the same number of $\CM^\times$-fixed points $[A'_\mu]$ for $\mu\in\PC_\nu(n\ell+|\nu|)$, see Corollary~\ref{coro:isom-if-same-Jcore}.
This implies that each $\CM^\times$-fixed point in $\XC_\bc(\alpha)$ is of the form $[A'_\mu]$.
\end{remark}

\subsection{Combinatorial parameterization of symplectic leaves}

\begin{lem}
\label{lem:equiv-simple-nd}
The following conditions are equivalent.
\begin{enumerate}
    \item[(i)]  The pair $(\alpha,\bc)$ is equivalent to a pair of the form $(n\delta,\bc')$ with $n\geqslant 0$.
    \item[(ii)] We have $\alpha\in E_\bc$.
\end{enumerate}
\end{lem}
\begin{proof}
$(ii)$ implies $(i)$ by Remark~\ref{rk:couple-equiv-nd}.

Now, let us prove that $(i)$ implies $(ii)$. Assume that $(\alpha,\bc)$ satisfies $(i)$. Since the isomorphism \eqref{eq:iso-lusztig} sends simple representations to simple representations by construction, it is enough to assume $\alpha=n\delta$. Let $\alpha^0$ be associated to $\alpha=n\delta$ and $\bc$ as in Proposition~\ref{prop:can-decomp-of-d}. Then $(ii)$ is equivalent to $\alpha^0=\alpha$. 

Assume that $\alpha^0\ne \alpha$. Since the pair $(\alpha^0,\bc)$ satisfies $(ii)$, it also satisfies $(i)$. So, it must be equivalent to a pair of the form $(n'\delta,\bc')$. Since we have $\alpha^0-n'\delta\in \ZM_{\geqslant 0}^{\ZM/\ell\ZM}$ and $0\ne n\delta-\alpha^0\in \ZM_{\geqslant 0}^{\ZM/\ell\ZM}$,  we get $n>n'$.

Now, we get $\XC_{\bc}(n\delta)\cong \XC_{\bc}(\alpha^0)$ by Proposition~\ref{prop:can-decomp-of-d} $(iii)$ and we have $\XC_{\bc}(\alpha^0)\cong \XC_{\bc'}(n'\delta)$ by \eqref{eq:iso-lusztig}. This is impossible because by Lemma~\ref{lem:dim-CM} we have $\dim \XC_{\bc}(n\delta)=2n$, $\dim \XC_{\bc'}(n'\delta)=2n'$ and $n'<n$.
\end{proof}

The following lemma is a combinatorial version of Proposition~\ref{prop:Cadmissiblemove}. 

\begin{lem}\label{lem:descr-Eth}
Assume that $\bc$ is $J$-standard. Then $\alpha\in E_\bc$ if and only if we have 
$$
\alpha=\Res_\ell(\nu)+r\delta
$$
with $r\geqslant 0$ and $\nu\in \CC_\ell\cap\CC_J$.
\end{lem}
\begin{proof}
The parabolic subgroup $W_J$ of $W^{\rm aff}$ is the stabilizer of $\bc$ in $W^{\rm aff}$.
Writing $\alpha=\Res_\ell(\nu)+r\delta$ as in \eqref{eq:d-via-core-n}, we have $r\in\ZM$ and $\nu\in\CC_\ell$.

\medskip
Assume $\alpha\in E_\bc$. Then Lemma~\ref{lem:equiv-simple-nd} implies that $r\geqslant 0$ and we can find $x\in W^\aff$ (see  Remark~\ref{rem:equiv=shortest}) such that $x(\alpha)=r\delta$ and $x$ is the shortest element in the coset $xW_J\in W^\aff/W_J$.  

Let $w$ be the shortest element in $x^{-1}W$. We have $\nu=x^{-1}(\emptyset)=w(\emptyset)$. Assume that $\nu$ is not a $J$-core. Then we have $|s_i(\nu)|<|\nu|$ for some $i\in J$; this corresponds to the case $(3)$ in  Remark~\ref{rem:W-act-lcores}. Then Lemma~\ref{lem:nu-w-i} implies $\ell(s_iw)<\ell(w)$. Then we also have $\ell(s_ix^{-1})<\ell(x^{-1})$ or equivalently $\ell(xs_i)<\ell(x)$. This contradicts the fact that $x$ is the shortest element in $xW_J$. Thus, $\nu$ must be a $J$-core.

\medskip
Now, assume that we have $\alpha=\Res_\ell(\nu)+r\delta$ for $r\geqslant 0$ and $\nu\in\CC_\ell\cap\CC_J$. Let $w$ be the element of $W^\aff$ such that $w(\emptyset)=\nu$ and such that $w$ is the shortest element in $wW$. It is enough to prove that $w$ in the shortest element in $W_Jw$. Indeed, if we prove this, then by Remark~\ref{rem:equiv=shortest} we have $(\alpha,\bc)\sim (w^{-1}\star\alpha,{(w^{-1}})^*(\bc))=(r\delta,(w^{-1})^*(\bc))$ and then by Lemma~\ref{lem:equiv-simple-nd} we have $\alpha\in E_\bc$.

 Since $\nu$ is a $J$-core, for each $i\in J$ we have $|s_i(\nu)|\geqslant |\nu|$. This means that for each $i\in J$, we are either in situation $(1)$ or in situation $(2)$ of Remark~\ref{rem:W-act-lcores}. In both cases Lemma~\ref{lem:nu-w-i} yields $\ell(s_iw)>\ell(w)$. 
\end{proof}

\begin{remark}
\label{rem:dim-nd+nu}
Assume that $\bc$ is $J$-standard and fix $\alpha\in E_\bc$. By Lemma~\ref{lem:descr-Eth} above, we can write $\alpha$ in the form $\alpha=\Res_\ell(\nu)+n\delta$ with $n\geqslant 0$ and $\nu\in \CC_\ell\cap\CC_J$. Then by Lemma~\ref{lem:equiv-simple-nd}, the pair $(\alpha,\bc)$ is equivalent to $(n\delta,\bc')$ for some $\bc'\in \CM^{\ZM/\ell\ZM}$. 
Then Lemma~\ref{lem:dim-CM} implies that the variety $\XC_\bc(\Res_\ell(\nu)+n\delta)$ is normal of dimension $2n$.
\end{remark}

We see that the elements of $E_\bc$ are in bijection with the pairs $(\nu,r)$ where $\nu$ is an $\ell$-core that is a $J$-core and $r\in \ZM_{\geqslant 0}$.

Assume that $\bc$ is $J$-standard. Then we have a partial order $\succcurlyeq$ on $E_\bc$ given by $\alpha\succcurlyeq\alpha'$ if $\alpha-\alpha'\in \sum_{j\in J} \ZM_{\geqslant 0} e_j$. In other words, we have $\alpha\succcurlyeq\alpha'$ if and only if $\LG^\alpha_{\alpha'}\ne \emptyset$. Using the bijection above, we may consider the order $\succcurlyeq$ as an order on the set $(\CC_\ell\cap\CC_J)\times \ZM_{\geqslant 0}$. 

\begin{lem}
We have $(\nu_1,r_1)\succcurlyeq (\nu_2,r_2)$ if and only if $r_1\geqslant r_2$ and there exists a partition $\lambda\in \PC_{\nu_1}(|\nu_1|+\ell(r_1-r_2))$ such that $\Core_J(\lambda)=\nu_2$.
\end{lem}
\begin{proof}
Assume $(\nu_1,r_1)\succcurlyeq (\nu_2,r_2)$. Then $\dim \XC_\bc(\Res_\ell(\nu_1)+r_1\delta)=2r_1$ and  $\dim \XC_\bc(\Res_\ell(\nu_2)+r_2\delta)=2r_2$ by Remark~\ref{rem:dim-nd+nu}. By Corollary~\ref{coro:norm-cl-leaf} and its proof, the normalization of the closure of the symplectic leaf $\LG^{\Res_\ell(\nu_1)+r_1\delta}_{\Res_\ell(\nu_2)+r_2\delta}$ is isomorphic to $\XC_\bc(\Res_\ell(\nu_2)+r_2\delta)$. In particular, 
$$
\dim\XC_\bc(\Res_\ell(\nu_1)+r_1\delta)\geqslant\dim\LG^{\Res_\ell(\nu_1)+r_1\delta}_{\Res_\ell(\nu_2)+r_2\delta}
$$ 
implies $r_1\geqslant r_2$.

Now, $(\nu_1,r_1)\succcurlyeq (\nu_2,r_2)$ implies $\Res_\ell(\nu_1)+r_1\delta \succcurlyeq \Res_\ell(\nu_2)+r_2\delta$ and then $(\nu_1,r_1-r_2)\succcurlyeq (\nu_2,0)$. This means that the variety $\XC_\bc(\Res_\ell(\nu_1)+(r_1-r_2)\delta)$ contains the symplectic leaf $\LG^{\Res_\ell(\nu_1)+(r_1-r_2)\delta}_{\Res_\ell(\nu_2)}$. This symplectic leaf is $0$-dimensional, so it is a $\CM^\times$-fixed point. Then by Section~\ref{subs:C-fixed-points}, this should be a point of the form $[A'_\lambda]$ for some $\lambda\in \PC_{\nu_1}(|\nu_1|+l(r_1-r_2))$. By Corollary~\ref{coro:structure-Amu}, we have $\dim^\reg(A'_\lambda)=\Res_\ell(\Core_J(\lambda))$. Then $[A'_\lambda]\in \LG^{\Res_\ell(\nu_1)+(r_1-r_2)\delta}_{\Res_\ell(\nu_2)}$ implies $\Core_J(\lambda)=\nu_2$.

\medskip
Conversely, if $r_1\geqslant r_2$ and there exists such a partition $\lambda$, then the $\CM^\times$-fixed point $[A'_\lambda]$ of $\XC_\bc(\Res_\ell(\nu_1)+(r_1-r_2)\delta)$ is a symplectic leaf. Since $\dim^\reg(A'_\lambda)=\Res_\ell(\Core_J(\lambda))=\Res_\ell(\nu_2)$, this is the symplectic leaf $\LG^{\Res_\ell(\nu_1)+(r_1-r_2)\delta}_{\Res_\ell(\nu_2)}$. Then we have $(\nu_1,r_1-r_2)\succcurlyeq (\nu_2,0)$. This implies $\Res_\ell(\nu_1)+(r_1-r_2)\delta \succcurlyeq \Res_\ell(\nu_2)$ and hence $(\nu_1,r_1)\succcurlyeq (\nu_2,r_2)$.
\end{proof}

Assume that $\bc$ is $J$-standard and $\alpha\in E_\bc$. Write $\alpha=\Res_\ell(\nu)+n\delta$, $\nu\in\CC_\ell\cap\CC_J,n\geqslant 0$.
\begin{cor}
For $\alpha'\in \ZM/\ell\ZM$, the following conditions are equivalent.
\begin{enumerate}
    \item[(i)] We have $\LG^\alpha_{\alpha'}\ne \emptyset$.
    \item[(ii)] There exists a partition $\lambda\in\PC_\nu(n'\ell+|\nu|)$ for some $n'\in[0;n]$ such that 
    \[
    \alpha'=\Res_\ell(\Core_J(\lambda))+(n-n')\delta.
    \]
\end{enumerate}
\end{cor}

\begin{proof}
Write $\alpha'=\Res_\ell(\nu')+r'\delta$. Then $\LG^\alpha_{\alpha'}\ne \emptyset$ is equivalent to $(\nu,n)\succcurlyeq (\nu',r')$. By the lemma above, this is equivalent to $n\geqslant r'$ and the existence of a partition $\lambda\in \PC_{\nu}(\ell(n-r')+|\nu|)$ such that $\Core_J(\lambda)=\nu'$. Moreover, the condition $\Core_J(\lambda)=\nu'$ is equivalent to $\Res_\ell(\Core_J(\lambda))=\Res_\ell(\nu')=\alpha'-r'\delta$. Now we see that $(i)$ is equivalent to $(ii)$ with $n'=n-r'$.
\end{proof}

In particular, we see that the symplectic leaves of $\XC_\bc(\alpha)$ are parameterized by $\ell$-cores of $J$-cores of elements of $\PC_\nu(n'\ell+|\nu|)$ for $n'\in [0;n]$. Note that by Lemma~\ref{lem:lcore-of-Jcore}, the $\ell$-cores of $J$-cores are also $J$-cores. We sum up this in the following proposition.

\begin{prop}
\label{prop:comb-param-leaves}
The symplectic leaves of $\XC_\bc(\alpha)$ are paremeterized by a subset of the set $\CC_\ell\cap\CC_J$. This subset is the image of the set $\coprod_{n'=0}^n \PC_\nu(n'\ell+|\nu|)$ by the map $\Core_\ell\circ \Core_J$.    
\end{prop}

Since each pair $(n\delta,\bc)$ is equivalent to a pair of the form $(\alpha,\bc')$ such that $\bc'$ is $J$-standard for some $J$ and $\alpha\in E_{\bc'}$ (see Lemma~\ref{lem:equiv-simple-nd}), the description above gives a parameterization of the symplectic leaves of an arbitrary Calogero--Moser variety of type $G(\ell,1,n)$ with $a\ne 0$.

\begin{example}
\label{ex:comb-l=2}

Assume $\ell=2$. In this case the set $\CC_2$ of $2$-cores is labeled by nonnegative integers. We have $\CC_2=\{\nu_t,t\in \ZM_{\geqslant 0}\}$ where $\nu_t$ is the partition staircase $\nu_t=(t,t-1,t-2,\ldots,2,1)$ of $t(t+1)/2$ ($\nu_0$ is the empty partition). The two possible non-trivial examples of $J$ are $J_0=\{0\}$ and $J_1=\{1\}$. Then the $2$-cores $\nu_2, \nu_4, \nu_6,\ldots$ are $J_0$-cores and not $J_1$-cores, the $2$-cores $\nu_1, \nu_3, \nu_5,\ldots$ are $J_1$-cores and not $J_0$-cores, the $2$-core $\nu_0=\emptyset$ is a $J_0$-core and a $J_1$-core.  
Let us write $c=(c_1,c_\gamma)$ for the parameters of the Calogero--Moser variety to make this example easily comparable to Section~\ref{sec:hyperoctahedral}. We always assume $c_1\ne 0$ here.
Let us understand how the combinatorial parametrization of symplectic leaves of the Calogero--Moser variety looks from the point of view of Proposition \ref{prop:comb-param-leaves}.
Assume that we have $c_\gamma=\pm mc_1$ with integer $m> 0$. Then the Calogero--Moser variety is isomorphic to the quiver variety $\XC_{\bc'}(n\delta)$ with $\bc'=(\pm m-1,\mp m)$. Then there is a $\bc'$-admissible $w$ such that 
$$
w(\bc')=
\begin{cases}
(-1,0) & \mbox{ if $m$ is even},\\
(0,-1) & \mbox{ if $m$ is odd}.
\end{cases}
$$
Set $\bc=w(\bc')$. Set $J=J_1$ if $m$ is even and $J=J_0$ if $m$ is odd.  The parameter $\bc$ is $J$-standand. Note that we have $w(\nu_0)=\nu_{m-1}$, in particular the Calogero--Moser variety is isomorphic to $\XC_{\bc'}(n\delta)\cong \XC_\bc(n\delta+\Res_\ell(\nu_{m-1}))$. 

Now, we would like to find all possible $\nu'\in (\Core_\ell\circ \Core_J)(\coprod_{n'=0}^n \PC_{\nu_{m-1}}(n'\ell+|\nu_{m-1}|))$. 
If $n\leqslant m$, then each $\lambda\in \coprod_{n'=0}^n \PC_{\nu_{m-1}}(n'\ell+|\nu_{m-1}|))$ is a $J$-core. (To see this, we need to use the presentation of cores by abaci, see \cite[Section 2]{BJV}.) In this case, we get $(\Core_\ell\circ \Core_J)(\lambda)=\Core_\ell(\lambda)=\nu_{m-1}$. So, the only possible $\nu'$ that we may get is $\nu'=\nu_{m-1}$. In this case the Calogero--Moser variety is smooth. Now, if $n\geqslant m+1$ then it is also possible to get $\nu'=\nu_{m+1}$. If $n\geqslant 2(m+2)$ then it is also possible to get $\nu'=\nu_{m+3}$, etc. If $n\geqslant k(m+k)$ then it is also possible to get $\nu'=\nu_{m-1+2k}$. Finally, we see that the symplectic leaves are labeled by the following subset of $\CC_2\cap\CC_J$: $\{\nu_{m-1},\nu_{m+1},\nu_{m+3},\ldots,\nu_{m-1+2k}\}$ where $k$ is the maximal non-negative integer such that $n\geqslant k(m+k)$.

The case $m=0$ is a bit special; $\nu_{m-1}$ does not make sense here. Assume $c_\gamma=0$. Then the Calogero--Moser variety is isomorphic to $\XC_\bc(n\delta)$ with $\bc=(-1,0)$. Set $J=J_1$. In this case, the Calogero--Moser variety is never smooth for $n>0$, the symplectic leaves are labeled by $\{\nu_{0},\nu_{1},\nu_{3},\ldots,\nu_{2k-1}\}$ where $k$ is the maximal positive integer such that $n\geqslant k^2$.
\end{example}

\section{Slices to symplectic leaves.}\label{sec:slices}

In this section we describe the transverse slices to the symplectic leaves in the Calogero--Moser variety.

\subsection{Transverse Singularities}

In this section, we explain how Crawley-Boevey's \'etale local normal form \cite{CBnormal} can be used to explicitly compute a transverse slice to each symplectic leaf $\LG$ in a quiver variety. A transverse slice $(S,o)$ to a $2k$-dimensional leaf $\LG$ is a pointed Poisson variety together with a local isomorphism $(\mf{M}_{\bc}(\alpha),p) \cong (\mathbb{A}^{2k} \times S, (0,o))$ of Poisson varieties for each $p \in \LG$. Here $\mathbb{A}^{2k}$ is equipped with the usual symplectic structure. "Local isomorphism'' usually means formally local, however it is consequence of Crawley-Boevey's construction that "local isomorphism'' will mean \'etale locally in this article. 
Since the local isomorphism is Poisson, $\{ o \}$ will be a symplectic leaf in $S$ so that the leaf in $\mathbb{A}^{2k} \times S$ containing $(0,o)$ is $\mathbb{A}^{2k} \times \{ o \}$.

Before we state the result, we introduce one piece of notation. Let $L \subset H$ denote the set of loops in the set $H$ of oriented edges of our graph $\mathsf{G}$. For each dimension vector $\alpha$, taking trace at loops defines a $G(\alpha)$-equivariant map $\Tr_L \colon \Rep(\alpha) \to \C^{2\ell}$, where $2 \ell = |L|$. The restriction of this map to $\mu^{-1}(\bc)$ descends to $\mf{M}_{\bc}(\alpha)$. We define $\mf{M}_{\bc}^0(\alpha) = \Tr_L^{-1}(0)$. Then, provided $\alpha_i \neq 0$ for all $i$, we have $\mf{M}_{\bc}(\alpha) \cong \mf{M}_{\bc}^0(\alpha) \times \C^{2 \ell}$.

Given a leaf $\LG \subset \mf{M}_{\bc}(\alpha)$ labeled by a representation type $\tau = (n_1,\beta^{(1)}; \ds;n_k,\beta^{(k)})$ as in \eqref{eq:taureptype}, we define a new graph $\mathsf{G}(\tau)$ (the \textit{ext-graph}) as follows. The vertices of $\mathsf{G}(\tau)$ are $f_1,\ds, f_k$, thought of as being labeled by the roots $\beta^{(i)}$. In $\mathsf{G}(\tau)$, there are $-(\beta^{(i)},\beta^{(j)})$ edges between $i$ and $j$ if $i \neq j$ and $p(\beta^{(i)})$ loops at $i$. The $k$-tuple $\bn = (n_1, \ds, n_k)$ forms a dimension vector for the ext-graph $\mathsf{G}(\tau)$. Note that when $\beta^{(i)} \neq \beta^{(j)}$, \cite[Proposition~2.6]{CBKimura} says that $-(\beta^{(i)},\beta^{(j)}) = \dim \Ext^1_{\Pi^{\bc}}(M_i,M_j) \ge 0$, where $M_i,M_j$ are simple $\Pi^{\bc}$-modules of dimension $\beta^{(i)}$ and $\beta^{(j)}$ respectively. If $\beta^{(i)}= \beta^{(j)}$ then $p(\beta^{(i)}) \ge 0$ implies that $-(\beta^{(i)},\beta^{(j)}) < 0$ only when $\beta^{(i)}$ is real, in which case $i = j$ because real roots can only appear once in any given representation type.

\begin{thm}\label{thm:CBLunaslice}
	For each $p \in \LG$ there is a local isomorphism 
	$$
	(\mf{M}_{\bc}(\mathsf{G},\alpha), p) \cong (\mf{M}_{0}(\mathsf{G}(\tau),\bn),0) \cong (\mf{M}_{0}^0(\mathsf{G}(\tau),\bn) \times \mathbb{A}^{2 \ell},(0,0)).
	$$
	Moreover, $(\mf{M}_{0}^0(\mathsf{G}(\tau),\bn) ,0)$ is a transverse slice to $\LG$ at $p$. 
\end{thm}

\begin{proof}
The local isomorphism follows directly from \cite[Corollary~4.10]{CBnormal}. The only thing to check is the last statement. For $(\mf{M}_{0}^0(\mathsf{G}(\tau),\bn) ,0)$ to be a transverse slice, we require (a) the isomorphism $\mf{M}_{0}(\mathsf{G}(\tau),\bn) \cong \mf{M}_{0}^0(\mathsf{G}(\tau),\bn) \times \mathbb{A}^{2 \ell}$ is Poisson, with the non-degenerate Poisson structure on $\mathbb{A}^{2 \ell}$, and (b) $\dim \mc{L} = 2 \ell$.

Recall that $\mc{O}(\mf{M}_{0}(\mathsf{G}(\tau),\bn))$ is generated by traces of oriented cycles. If $\mf{l}$ is the necklace Lie algebra of $\mathsf{G}(\tau)$ then the morphism $\mf{l} \to \mc{O}(\mf{M}_{0}(\mathsf{G}(\tau),\bn))$ given by $\omega \mapsto \Tr \omega$ is a surjective Lie algebra homomorphism, where the codomain is a Lie algebra via the Poisson bracket; see \cite[Theorem~1.8]{VdenBDouble}. As a closed subvariety, $\mathbb{A}^{2 \ell} \subset \mf{M}_{0}(\mathsf{G}(\tau),\bn)$ is the zero set of all functions $\Tr \omega$ as $\omega$ runs over all oriented cycles that contain at least one arrow between distinct vertices (equivalently, whose support is more than one vertex). If $I \subset \mf{l}$ is the span of all such cycles then one can check from the necklace formula that $I$ is an ideal of $\mf{l}$. Moreover, the formula implies that $\mc{O}(\mf{M}_{0}(\mathsf{G}(\tau),\bn)) \to \mc{O}( \mathbb{A}^{2 \ell} )$ is a morphism of Poisson algebras confirming (a).

The leaf in $\mf{M}_{0}(\mathsf{G}(\tau),\bn)$ containing the point $0$ is labeled by representation type $(n_1,f_1;n_2,f_2;\ds; n_k,f_k)$ and hence has dimension 
\[
2 \sum_{i=1}^k p(f_i) = 2 \sum_{i=1}^k p(\beta^{(i)}) = 2 \ell,
\]
confirming (b). 
\end{proof}

We end with basic observations that allow us to explicitly identify the transverse slices in examples; these observations will be useful later. The first is immediate. 

\begin{lem}\label{lem:quiverbraking1}
	Assume that $0$ is a vertex for the graph $\mathsf{G}$ and $\alpha$ a dimension vector with $\alpha_0 = 1$. If $\ell$ is the number of loops at vertex $0$ then $\mf{M}_{\bc}(\mathsf{G},\alpha) \cong \mf{M}_{\bc}(\mathsf{G}^{\circ},\alpha) \times \C^{2 \ell}$, where $\mathsf{G}^{\circ}$ is the graph obtained by removing the $\ell$ loops at $0$. 
\end{lem}

\begin{remark}\label{rem:loopsatone}
In the setting of Theorem~\ref{thm:CBLunaslice}, assume that $n_i = 1$ for all imaginary $\beta^{(i)}$. Since $\beta$ is imaginary if and only if $p(\beta) > 0$, this means that $\mbf{n}_i = 1$ for all vertices $i$ in $\mathsf{G}(\tau)$ with at least one loop. Then combining Theorem~\ref{thm:CBLunaslice} with Lemma~\ref{lem:quiverbraking1} shows that $\mf{M}_{0}^0(\mathsf{G}(\tau),\bn) \cong \mf{M}_{0}(\mathsf{G}(\tau)^{\circ},\bn)$, where $\mathsf{G}(\tau)^{\circ}$ now means "remove all loops in $\mathsf{G}(\tau)$".
\end{remark}

Let $\mathsf{G}(1), \ds, \mathsf{G}(k)$ be a collection of graphs with dimension vectors $\alpha(1), \ds, \alpha(k)$ and distinguished vertex $0_j \in \mathsf{G}(j)$ such that $\alpha(j)_{0_j} = 1$ and $\mathsf{G}(j)$ has no loops at $0_j$ for all $j$. We form the new graph $\mathsf{G}$ by gluing these graphs at the distinguished vertex (labeled $0$ in $\mathsf{G}$). Note that each $\mathsf{G}(j)$ is a full subgraph of $\mathsf{G}$. There is a dimension vector $\alpha$ for $\mathsf{G}$ with $\alpha_0 = 1$ such that $\alpha |_{\mathsf{G}(j)} = \alpha(j)$ for all $j$. Let $\bc(j)$ be parameters for the $\mathsf{G}(j)$ such that $\bc(j)(e_{0_j}) = 0$. They extend to a parameter $\bc$ for $\mathsf{G}$ with $\bc(e_0) = 0$.

\begin{lem}\label{lem:quiverbraking2}
Assume that the graphs $\mathsf{G}(1), \ds, \mathsf{G}(k)$ etc. are given as above. If $\dim \mf{M}_{\bc}(\mathsf{G},\alpha) = 2 p(\alpha)$ and $\dim \mf{M}_{\bc(j)}(\mathsf{G}(j),\alpha(j)) = 2 p (\alpha(j))$ for all $i = 1, \ds, k$ then there is an isomorphism  
		$$
	\prod_{j = 1}^k \mf{M}_{\bc(j)}(\alpha(j)) \iso \mf{M}_{\bc}(\alpha).
		$$

\end{lem}

\begin{proof}
Let $\Rep(\alpha)$ be the representation space, with $\bw = 0$, as in \eqref{eq:repspace}. Then $\Rep(\alpha) = \prod_{j = 1}^k \Rep(\alpha(j))$ and there is a diagonal embedding $G(\alpha) \hookrightarrow \prod_{j = 1}^k G(\alpha(j))$ such that the identification is $G(\alpha)$-equivariant. We claim that the resulting embedding 
\[
\bigotimes_{j = 1}^k \C[\Rep(\alpha(j))]^{G(\alpha(j))} \hookrightarrow    \C[\Rep(\alpha)]^{G(\alpha)}
\]
is an isomorphism. As noted in the proof of Proposition~\ref{prop:leafclosurenormal}, $\C[\Rep(\alpha)]^{G(\alpha)}$ is generated by all traces $\Tr \omega$ of oriented cycles $\omega$ in $\mathsf{G}$. If such a cycle does not pass through $0$ then it is entirely contained in one $\mathsf{G}(j)$ and is clearly invariant under the larger group. If $\omega$ does pass through $0$, then 
we may write
    \[
       \Tr \omega = \prod_{x} \Tr \omega{(x)}
    \]
    where $\omega{(x)}$ is a cycle beginning and ending at $0$, passing through that vertex just once. Then $\omega{(x)}$ again lies entirely in some $\mathsf{G}(j)$ and is invariant under the larger group. The claim follows.

Next, the relations defining $\mu^{-1}(\bc)$ at vertex $0$ can be written $\sum_{j = 1}^k \mu_j^{-1}(\bc(j)) = 0$, where $\mu_j$ is the moment map for $G(\alpha(j))$ acting on $\Rep(\alpha(j))$. Thus, $I(\mu^{-1}(\bc)) \subset I(\mu_1^{-1}(\bc(1)),\ds, \mu_k^{-1}(\bc(k)))$ in $\C[\Rep(\alpha)]$ and we get a surjection
\begin{align*}
\C[\mu^{-1}(\bc)]^{G(\alpha)} & = \left( \C[\Rep(\alpha)] / I(\mu^{-1}(\bc))\right)^{G(\alpha)} \to \left( \C[\Rep(\alpha)] / I(\mu_1^{-1}(\bc(1)),\ds, \mu_k^{-1}(\bc(k)))\right)^{G(\alpha)} \\
& = \left( \C[\Rep(\alpha)] / I(\mu_1^{-1}(\bc(1)),\ds, \mu_k^{-1}(\bc(k)))\right)^{\prod_j G(\alpha(j))} = \bigotimes_{j = 1}^k \C[\mu_j^{-1}(\bc(j))]^{G(\alpha(j))}. 
\end{align*}
Thus, there is a closed embedding 
\begin{equation}\label{eq:productclosedembed}
	\prod_{j = 1}^k \mf{M}_{\bc}(\alpha(j)) \hookrightarrow \mf{M}_{\bc}(\alpha).    
\end{equation}
Finally, write $\alpha(j) = \alpha(j)' + e_{0}$ so that $(\alpha(j)',\alpha(s)') = 0$ for $j \neq s$ and $\alpha = \alpha(1)' + \cdots + \alpha(k)' + e_{0}$. Then, using that $(e_0,e_0) = 2$,  
\begin{align*}
    p(\alpha) & = 1 - (1/2)( \alpha(1)' + \cdots + \alpha(k)' + e_{0},\alpha(1)' + \cdots ) \\
    & = - \sum_{i = 1}^k  ((1/2)(\alpha(i)',\alpha(i)') + (e_{0},\alpha(i)'))  \\
    & = k - (1/2) \sum_{i = 1}^k \left( (\alpha(i)',\alpha(i)') + 2 (e_{0},\alpha(i)') + (e_{0},e_0) \right) \\
    & = \sum_{i = 1}^k p(\alpha(i)).
\end{align*}
Therefore, the assumption $\dim \mf{M}_{\bc}(\mathsf{G},\alpha) = 2 p(\alpha)$ and $\dim \mf{M}_{\bc(j)}(\mathsf{G}(j),\alpha(j)) = 2 p (\alpha(j))$ for all $i = 1, \ds, k$ implies that \eqref{eq:productclosedembed} is a closed embedding of irreducible varieties of the same dimension. Thus, it is an isomorphism. 
\end{proof}

Lemma~\ref{lem:quiverbraking1} says that we can forget any loops at a vertex $i$ with $\alpha_i = 1$ and Lemma~\ref{lem:quiverbraking2} says that we can "break'' the graph into its connected parts at vertex $i$ if it has no loops and $\alpha_i = 1$. 

\subsection{Transverse Singularities in Calogero--Moser varieties}\label{sec:nonzerotransverse}

Assume first that $\bc(\delta) = -1$. We fix a leaf $\LG(\beta)$ in $\ZCB_{c}(\Gamma_n)$. Recall that the set $\Delta(\bc) = \{ \eta^{(1)}, \ds, \eta^{(s)} \}$ of minimal roots in $R^+_{\bc}$ is a set of simple roots for the root system $R_{\bc}$.

As explained previously, we may assume that $w(\Delta(\bc)) = w(R^+_{\bc}) \cap \Delta$ and hence $R^+_{\bc}$ is conjugate to a parabolic subsystem of $R$. Let $\mathsf{G}(\bc)$ be the subgraph of the affine Dynkin graph $\mathsf{G}(\Gamma)$ obtained by deleting the vertices not in $w(R^+_{\bc}) \cap \Delta$. Then $\mathsf{G}(\bc)$ is a disjoint union of finite Dynkin diagrams. Since $\beta \in \Xi(\bc)$, we may write $\beta = v_1 \eta^{(1)} + \cdots + v_s \eta^{(s)}$ and think of $\bv = (v_1, \ds, v_s)$ as a dimension vector for $\mathsf{G}(\bc)$. We define a framing vector $\bw = (w_1, \ds, w_s)$ for $\mathsf{G}(\bc)$ by 
$$
w_i := (\beta + \Lambda_0, \eta^{(i)}).
$$

\begin{thm}\label{thm:transversenonzerolevel}
	Let $\beta \in \Xi(\bc)$ with $\varrho(\beta) \le n$.  The transverse slice to $\LG(\beta)$ in $\ZCB_{c}(\Gamma_n)$ is isomorphic to the framed quiver variety $\mf{M}_0(\mathsf{G}(\bc),\bw,\bv)$. 
\end{thm}

\begin{proof}
Recall from the proof of Theorem~\ref{thm:nonzerolevelleaves} that $\beta$ corresponds to the representation type 
\[
\tau = (\gamma(m,\nu),(v_1,\eta^{(1)};\ds; v_s,\eta^{(s)})),
\]
where $\nu = \beta - \beta_0 \delta$ and $m = n - \qu(\beta)$. The ext-graph $\mathsf{G}(\tau)$ in this case is obtained from $\mathsf{G}(\bc)$ by adding one additional vertex $\infty$ corresponding to the vector $\gamma(m,\nu)$ and $m = p(\gamma(m,\nu))$ loops at the vertex $\infty$. The dimension at vertex $\infty$ is $1$ and the number of arrows from $\infty$ to $\eta^{(i)}$ is 
\[
 - (e_{\infty} + \gamma(m,\nu),\eta^{(i)}) = (\beta - n \delta - e_{\infty},\eta^{(i)}) = (\beta - e_{\infty}, \eta^{(i)}) = (\beta + \Lambda_0, \eta^{(i)})=w_i.
\]
Therefore, Lemma~\ref{lem:quiverbraking1} says that $\mf{M}_{0}(\mathsf{G}(\tau),e_{\infty} + \bv) \cong \mf{M}_{0}(\mathsf{G}(\tau)^{\circ},e_{\infty} + \bv) \times \C^{2m}$; see also Remark~\ref{rem:loopsatone}. But $\mathsf{G}(\tau)^{\circ}$ equals the graph $\mathsf{G}(\bc)'$ we obtain by deframing $\mathsf{G}(\bc)$ with respect to framing vector $\bw$; see Section~\ref{sec:CBtrick}. In other words, 
\[
\mf{M}_{0}(\mathsf{G}(\tau)^{\circ},e_{\infty} + \bv) \cong \mf{M}_0(\mathsf{G}(\bc),\bw,\bv).
\]
Therefore, Theorem~\ref{thm:CBLunaslice} implies that a transverse slice to $\LG(\beta)$ is given by $\mf{M}_0(\mathsf{G}(\bc),\bw,\bv)$. 
\end{proof}

Our next goal is to show that if $\mathsf{G}$ is any finite type (ADE) graph and $\Gamma\subset \SL(2,\C)$ the group corresponding to the affine Dynkin graph $\widetilde{\mathsf{G}}$ then every framed Nakajima quiver variety associated to $\mathsf{G}$ can be realized as a transverse slice to $\LG$ in $\ZCB_{c}(\Gamma_n)$ for some leaf $\LG$ and $n,c$. 

\begin{lem}\label{lem:rootssimplesub}
    Let $\Phi \subset R$ be the finite and affine root system respectively associated to $\Gamma \subset \SL(2,\C)$. Assume $\bc(\delta) = 1$ and we are given $\{ \eta^{(1)}, \ds, \eta^{(r)} \} \subset R^+_{\bc}$ such that the Cartan matrix $C$, with $C_{i,j} := (\eta^{(i)},\eta^{(j)})$ for $1 \le i, j \le r$, is equivalent to the Cartan matrix of $\Phi$. Then $\{ \eta^{(1)}, \ds, \eta^{(r)} \}$ is a set of simple roots for $R_{\bc}^+$ and $R_{\bc}^+ \cong \Phi$. 
\end{lem}

\begin{proof}
First we note that $C_{i,j} \le 0$ for $i \neq j$ since $C$ is assumed equivalent to the Cartan matrix for $\Phi$. This implies that $\{ \eta^{(1)}, \ds, \eta^{(r)} \}$ is a set of simple roots for a root subsystem $\Psi$ of $R_{\bc}$ by \cite[Lemma~1]{DyerLehrer}. The rank of $C$ is $r$. By Lemma~\ref{lem:wparabolicrootsystem},  $R_{\bc}$ is a equivalent to a parabolic root system in $R$. Therefore, the rank of $R_{\bc}^+$ is at most $r$ too. This means that $R_{\bc}$ must also have rank exactly $r$. Since $R_{\bc}$ is a equivalent to a parabolic root system in $R$, its Cartan matrix $C(\bc)$ with respect to $\Delta(\bc)$ is equivalent to a (proper) submatrix of the Cartan matrix for the affine root system $R$. 

Therefore, we are in the following situation: we have a finite root system $R_{\bc}$ obtained by deleting one vertex of the simply-laced affine Dynkin graph $\widetilde{\mathsf{G}}$, which contains a subsystem $\Psi$ of type $\mathsf{G}$. We claim that $\Psi = R_{\bc}$. First note that $R_{\bc}$ must be irreducible. Indeed, $\Psi$ is an irreducible subsystem of the same rank as $R_{\bc}$, which means that the Weyl group $W(\Phi)$ acts irreducibly on its reflection representation $\mathfrak{h} = \Span_{\R}(\Psi) = \Span_{\R}(R_{\bc})$, so $W(R_{\bc})$ must also act irreducibly since $W(\Psi) \subset W(R_{\bc})$. In type $\mathsf{A},\mathsf{D},\mathsf{E}_6$ the equality is clear since all irreducible parabolic subsystems of $\widetilde{\mathsf{A}}_{r}$, resp. $\widetilde{\mathsf{D}}_{r}$ or $\widetilde{\mathsf{E}}_{6}$, obtained by deleting one vertex are of type $\mathsf{A}_{r}$, resp. of type $\mathsf{D}_r$ or type $\mathsf{E}_6$. In $\widetilde{\mathsf{E}}_{7}$, the irreducible parabolic subsystems of rank $7$ are of type $\mathsf{E}_{7}$ or $\mathsf{A}_{7}$. But $|\Psi| = |\mathsf{E}_{7}| = 126 > 56 = |\mathsf{A}_{7}|$ so we must have $\Psi = R_{\bc}$ of type $\mathsf{E}_{7}$. Finally, in $\widetilde{\mathsf{E}}_{8}$, the irreducible parabolic subsystems of rank $7$ are of type $\mathsf{E}_{8}, \mathsf{D}_8,\mathsf{A}_{8}$. But $|\Psi| = |\mathsf{E}_{8}| = 240 > 112 = |\mathsf{D}_{8}|, 72 = |\mathsf{A}_{8}|$, so we must have $\Psi = R_{\bc}$. 
    \end{proof}

\begin{prop}\label{prop:transversearbitrary}
    Let $\mathsf{G}$ be a finite type Dynkin graph and $\Gamma\subset \SL(2,\C)$ the group corresponding to the affine Dynkin graph $\widetilde{\mathsf{G}}$. For any pair $(\bw,\bv)$ of dimension vectors for $\mathsf{G}$, there exists $n, c$ and a leaf $\LG$ such that $\mf{M}_0(\mathsf{G},\bw,\bv)$ is isomorphic to the transverse slice to $\LG$ in $\ZCB_{c}(\Gamma_n)$.   
\end{prop}

\begin{proof}
Let $\Phi$ denote the finite root system associated to $\mathsf{G}$. First note that \cite[Proposition~3.9]{McGertyNevinsGalois} says that there exists $\bv' \le \bv$ such that $\mf{M}_0(\mathsf{G},\bw,\bv) \cong \mf{M}_0(\mathsf{G},\bw,\bv')$ with $\sum_{i=1}^r w_i \Lambda_i - \sum_i v_i' e_i$ a dominant weight; here $\bv' = (v_1', \ds, v_r')$. This means that 
    \[
    k_i := w_i - \sum_j v_j' (e_i,e_j) \in \Z_{\ge 0}, \quad \forall \, 1 \le i \le r. 
    \]
    Let $\bc_i = -k_i$ and $\bc_0 = 1 + \sum_{i = 1}^r \delta_i k_i$ so that $\bc(\delta) = 1$. Let $\eta^{(i)} = e_i + k_i \delta$. Then $\eta^{(i)} \in R_{\bc}^+$ for all $i$. Notice that $C_{i,j} := (\eta^{(i)},\eta^{(j)}) = (e_i,e_j)$, implying that the Cartan matrix $C$ equals the Cartan matrix for $\Phi$. Lemma~\ref{lem:rootssimplesub} then implies that $\Delta(\bc) = \{ \eta^{(1)}, \ds, \eta^{(r)} \}$.  If $\beta := \sum_{i = 1}^r v_i' \eta^{(i)}$ then we claim that $\beta \in \Xi(\bc)$. Indeed, assume that $\eta = \sum_{i = 1}^r u_i \eta^{(i)} \in \SS{\bc}$ with $\eta \succ \beta$; that is, $u_i \ge v_i'$ for all $i$ and $\eta \neq \beta$. Then, 
    \begin{align*}
        2(\qu(\eta) - \qu(\beta)) & = 2 \eta_0 + (\eta,\eta) - 2 \beta_0 - (\beta,\beta) \\
        & = 2 \sum_{i = 1}^r u_i k_i + (\eta,\eta) - 2 \sum_{i = 1}^r v_i' k_i - \sum_{i,j} v_i' v_j' (e_i,e_j) \\
        & = 2 \sum_{i = 1}^r u_i \left(w_i - \sum_{j = 1}^r v_j' (e_i,e_j) \right) + (\eta,\eta) - 2 \sum_{i = 1}^r v_i' \left(w_i - \sum_{j = 1}^r v_j' (e_i,e_j) \right)- \sum_{i,j} v_i' v_j' (e_i,e_j) \\
        & = 2 \sum_{i = 1}^r (u_i - v_i') w_i - 2 \sum_{i,j} u_i v_j' (e_i,e_j) + \sum_{i,j} u_i u_j(e_i,e_j) + \sum_{i,j} v_i'v_j' (e_i,e_j)\\
        & = 2 \sum_{i = 1}^r (u_i - v_i') w_i + (\mbf{u} - \bv',\mbf{u} - \bv'),
    \end{align*}
    where $\mbf{u} = \sum_{i = 1}^r u_i e_i$. Now, $( - , - )$ is positive definite on $\boldsymbol{\Lambda}$ and $0 \neq \mbf{u} - \bv' \in \boldsymbol{\Lambda}$. Hence, $(\mbf{u} - \bv',\mbf{u} - \bv') > 0$. This implies that $\qu(\eta) > \qu(\beta)$ and thus $\beta \in \Xi(\bc)$. 

    Finally, note that 
    \[
    (\beta+ \Lambda_0, \eta^{(i)}) = \sum_{j = 1}^r v_j'(\eta^{(j)},\eta^{(i)}) + (\Lambda_0, \eta^{(i)}) = \sum_{j = 1}^r v_j'(e_j,e_i) + k_i = w_i. 
    \]
    Then the claim of the proposition follows from Theorem~\ref{thm:transversenonzerolevel} if we take any $n \ge \varrho(\beta)$. 
\end{proof}

Note that taking $n = \varrho(\beta)$ in the proof of Proposition~\ref{prop:transversearbitrary} shows that the finite type quiver variety $\mf{M}_0(\mathsf{G},\bw,\bv)$ can be realized as a transverse slice to a zero-dimensional leaf in some Calogero--Moser variety. 

\begin{rem}\label{rem:rangeofslices}
If $\mathsf{G} = \mathsf{A}_{\ell-1}$, for some $\ell \ge 2$, then the quiver varieties $\mf{M}_0(\mathsf{G},\bw,\bv)$ are isomorphic to type $\mathsf{A}$ Slodowy slices $S(\mu,\nu)$ (and any type $\mathsf{A}$ Slodowy slice can be realized this way); see \cite[Section~8]{Nak1994} and \cite{MaffeiA}. By \cite{MVyb}, we may equivalently identify $\mf{M}_0(\mathsf{G},\bw,\bv)$ with a slice in the affine Grassmaniann of type $\mathsf{A}$. Therefore, given any Slodowy slice of type $\mathsf{A}$ or any slice in the affine Grassmaniann of type $\mathsf{A}$, one can always find $n,c$ and a leaf $\LG \subset \ZCB_{c}(\Gamma_n)$, with $\Gamma = \Z / \ell \Z$, such that the singularity transverse to $\LG$ is isomorphic to this slice. 
\end{rem}

Now we assume that $\bc(\delta) = 0$. For each irreducible factor $\Phi_i$ of $\Phi_{\bc}$ we have a finite subgroup $\Gamma(i)$ of $\SL(2,\C)$ whose affine Dynkin diagram (via the McKay correspondence) is the affinization of the Dynkin diagram of $\Phi_i$. 

\begin{lem}\label{lem:transversezerolevel}
	If $\rho = (0,(n,\delta - \theta(i); n h(i)_1,\alpha(i)_1; \dots ))$ and $p \in \LG{(\emptyset,\rho)}$, then $(\ZCB_{c}(\Gamma_n),p) \cong (\C^{2n}/\Gamma(i)_n, 0)$.  
\end{lem}

\begin{proof}
	If we consider first the representation type $(1, \delta - \theta(i);h(i)_1,\alpha(i)_1; \dots )$ then the corresponding ext-graph is just $\mathsf{G}(\Gamma(i))$, with dimension vector $\delta(i)$, the minimal imaginary root for $\mathsf{G}(\Gamma(i))$. This implies that the ext-graph for $(\emptyset,\rho)$ is the affine Dynkin diagram associated to $\Gamma(i)$ and the dimension vectors are $(\Lambda_0, n \delta(i))$. Then the isomorphism follows from Theorem~\ref{thm:CMquiveriso}. 
\end{proof}

\begin{prop}
	Assume $\bc(\delta) = 0$. If $p \in \LG{(\lambda,\rho)}$ then 
	$$
	(\ZCB_{c}(\Gamma_n),p) \cong \left(\C^{2 |\lambda|}/ \s_{\lambda} \times \prod_{i = 1}^{\ell} \C^{2 \rho_i} / \Gamma(i)_{\rho_i}, 0\right).
	$$
\end{prop}

\begin{proof}
Recall from \eqref{eq:degeneratelambdadecomp} that the leaf $\LG_{(\lambda,\rho)}$ is labeled by representation type 
\[
		(0,(\lambda_1,\delta;\lambda_2,\delta;\ds;\rho_1,\delta - \theta(1);\rho_1 h(1)_1, \alpha(1)_1; \dots ;\rho_1 h(1)_{r_1}, \alpha(1)_{r_1};\rho_2,\delta - \theta(2);\ds )).    
\]
	The ext-graph associated to $(\lambda,\rho)$ has a central vertex corresponding to the factor $e_{\infty}$. The dimension at this vertex is one. One can check that the hypothesis of Lemma~\ref{lem:quiverbraking2} hold because $p(e_{\infty} + n \delta) = p_{\Lambda_0}(n \delta) = n$ and 
 \[
 \sum_{i = 1}^{\ell(\lambda)} p_{\Lambda_0}(\lambda_i \delta) + \sum_{j = 1}^{\ell} p_{\Lambda_0}(\rho_j \delta(j)) =  \sum_{i = 1}^{\ell(\lambda)} \lambda_i + \sum_{j = 1}^{\ell} \rho_j = |\lambda | + |\rho | = n.
 \]
Therefore, the associated quiver variety can be expressed as a product, one factor for each connected component of the graph we get by breaking the graph at this central vertex; see \eqref{eq:diagramquiverbreaking} for a visualization. The representation type $(0,(\lambda_i,\delta))$ has ext-graph the one vertex and one loop (or Jordan) graph and dimension vectors $(\Lambda_0,\lambda_i e_0)$. The ext-graph of the representation type $(0,(\rho_i,\delta - \theta(i);\rho_i h(i)_1, \alpha(i)_1; \dots ;\rho_i h(i)_{r_i}, \alpha(i)_{r_i}))$ is described in Lemma~\ref{lem:transversezerolevel}. Thus,
 \[
 \mf{M}_0(\mathsf{G}(\tau),\mathbf{n}) \cong \prod_{j = 1}^{\ell(\lambda)} \XC_0(\mathsf{G}_J, \lambda_j \delta) \times \prod_{i = 1}^{\ell} \XC_0(\mathsf{G}(\Gamma),\rho_i \delta(i)). 
 \]
It has been shown in \cite[Lemma~2.11]{GanGinzburg} that each $\XC_0(\mathsf{G}_J,\lambda_j \delta)$ is isomorphic to $\C^{2 \lambda_j} / \s_{\lambda_j}$. Theorem~\ref{thm:CMquiveriso} says that $\XC_0(\mathsf{G}(\Gamma),\rho_i \delta(i)) \cong \C^{2 \rho_i} / \Gamma(i)_{\rho_i}$. The claim follows. 
\end{proof}

\section{The hyperoctahedral group}\label{sec:hyperoctahedral}

As an extended example, we consider the case where $\Gamma_n$ is the hyperoctahedral group. That is, $\Gamma = \Z / 2 \Z$ and $\Gamma_n$ is the Weyl group of type $\mathsf{B}_n$. 

\subsection{The proof of Theorem~\ref{thm:mainsing}}\label{sec:mainsingproof}

The rational Cherednik algebra $\Hb_{c}(\Gamma_n)$ (at $t = 0$) associated to the Weyl group $\Gamma_n = W(\B_n)$ depends on the choice of a pair of parameters $c = (c_1,c_{\gamma}) \in \C^2$. We recall the following theorem by Martino \cite[Theorem~8.2]{MarsdenWeinsteinStratification}. 

\begin{thm}\label{thm:BSympleaves} 
	Let $c = (c_1,c_{\gamma})$.
	\begin{enumerate}
		\item[(i)] \label{thm:BSympleaves:singular}  $\ZCB_{c}(\Gamma_n)$ is singular if and only if $c_1 = 0$ or $c_{\gamma} = m c_1$ for some integer $-n < m < n$.  
		\item[(ii)] If $c_1 = 0$ but $c_{\gamma} \neq 0$, then the symplectic leaves of $\ZCB_{c}(W)$ are parameterized by the set $\mc{P}(n)$ of partitions of $n$. For $\lambda \in \mc{P}(n)$, the corresponding leaf $\LG_{\lambda}$ has dimension $2 \ell(\lambda)$, where $\ell(\lambda)$ is the length of $\lambda$.
		\item[(iii)] If $c_{\gamma} = m c_1$, with $c_1 \neq 0$, then there is a bijection $k \mapsto \LG_k$,
		$$
		\{ \textrm{symplectic leaves $\LG$ of $\ZCB_{c}(W)$} \ \} \stackrel{1 : 1}{\longleftrightarrow} \{ k \in \Z_{\ge 0} \ | \ k(k+m) \le n \} \;.
		$$
		Moreover, $\dim \LG_{k} = 2(n - k(k+m))$.
	\end{enumerate}
\end{thm}

In this case, the graph $\mathsf{G}$ is the affine Dynkin quiver of type $\widetilde{\mathsf{A}}_1$ (with $2$ vertices). 
The minimal imaginary root is $\delta = e_0 + e_1$. By Theorem~\ref{thm:CMquiveriso}, there is an isomorphism $\ZCB_{c}(\Gamma_n) \cong \XC_{{\bc}}(n \delta)$, where 
 \begin{equation}\label{eq:lambdakappc1}
\bc = (- c_1 + c_{\gamma}, - c_{\gamma}).
 \end{equation}
Part (ii) of Theorem~\ref{thm:BSympleaves} is a special case of Proposition~\ref{prop:labelleaveslevel0}, where the condition $c_{\gamma} \neq 0$ means that the surface $X_{\bc}$ is smooth and hence the leaves are labeled by pairs $(\lambda,\rho)$, where $\rho$ is an empty composition. The leaf $\LG({\lambda})$ for $\lambda \in \mc{P}(n)$ is labeled by the representation type
\begin{equation}\label{eq:leafLrho}
	\tau = (0,(\lambda_1,\delta; \ds ;\lambda_{\ell},\delta)).
\end{equation}
Part (iii) of Theorem~\ref{thm:BSympleaves} is a special case of Theorem~\ref{thm:nonzerolevelleaves}. It was already shown in \cite[Proposition~5.7]{MartinoThesis} that, when $c_{\gamma} = \pm m c_1$, the leaf $\LG_k$ is labeled by the representation type
\begin{equation}\label{eq:leafLk}
	\tau = (n \delta - k \eta^{(1)},(k,\eta^{(1)})),
\end{equation}
where $\eta^{(1)}$ is unique (real) root in $R^+_{\bc}$. Thus, $\beta = k \eta^{(1)} \in \Xi({\bc})$. Specifically, $\eta^{(1)} = m e_0 + (m-1) e_1$ when $c_{\gamma} = m c_1$ with $1 \le m < n$ and $\eta^{(1)} = m e_0 + (m+1) e_1$ when $c_{\gamma} = - m c_1$ with $0 \le m < n$. See also Example \ref{ex:comb-l=2} for the combinatorial parameterization of symplectic leaves in this case.

Recall from the introduction that $\mc{O}(k,N)$ denotes the $GL(N)$-orbit of all matrices $X \in \mf{gl}(N)$ of rank $k$ with $X^2 = 0$. 

\begin{lem}\label{lem:singLk}
	Assume $c_{\gamma} = m c_1$ for $1 \le m \le n-1$ and $c_1 \neq 0$. Then there is an isomorphism of \'etale germs 
    \[
    (\ZCB_{c}(\Gamma_n),p) \cong (\overline{\mc{O}(k,2k+m)} \times T^* (\C^{n - k(k+ m)}),0)
    \]
    for any $p \in \LG_k$ and $k \in \Z_{\ge 0}$ such that $k(k+m) \le n$. 
\end{lem}

\begin{proof}
We note, for the computations below, that 
$$
(e_{\infty},\delta) = -1, \quad (e_{\infty}, \eta^{(1)}) = -m, \quad (\delta,\eta^{(1)}) = 0. 
$$
The lemma is a special case of Theorem~\ref{thm:transversenonzerolevel}. In this case, there exists simple representations $M,N$ of the deformed preprojective algebra such that $\dim M = e_{\infty} + n \delta - k \eta^{(1)}$, $\dim N = \eta^{(1)}$ and $p$ corresponds to the point $M \oplus N^{\oplus k}$. We have
\begin{align}
	p(e_{\infty} + n \delta - k \eta^{(1)}) & = p_{\Lambda_0}(n \delta - k \eta^{(1)}) \\
 & = 1 - (1/2) ((e_{\infty},e_{\infty}) + 2n (e_{\infty}, \delta) + k^2 (\eta^{(1)},\eta^{(1)}) - 2 k (e_{\infty},\eta^{(1)})) \\
  & = n - k(m + k).  \label{eq:peinftycomputation}  
\end{align}
    Theorem~\ref{thm:CBLunaslice} says that the germ $(\XC_{{\bc}}(n \delta),p)$ is isomorphic to the quiver variety $(\mathfrak{M}_0(\mathsf{G}(\bc),(k,1)),0)$, where $\mathsf{G}(\bc)$ is the graph with $2$ vertices, $p(\eta^{(1)}) = 0$ loops at the first vertex, $p(e_{\infty} + n \delta - k \eta^{(1)}) = n - k(m + k)$ loops at the second vertex and 
	$$
	- (e_{\infty} + n \delta - k \eta^{(1)},\eta^{(1)}) = 2k + m
	$$
	edges between the first and second vertex. This means that $\mathsf{G}(\bc)$ is the $(2k+m)$-Kronecker graph. Notice that $(k,1)$ is an indivisible imaginary root with $p((k,1)) = k(k+m)$. It is straight-forward to check that $\mathfrak{M}_0(\mathsf{G}(\bc),(k,1)) = \overline{\mc{O}(k,2k + m)}$.
\end{proof}

Identical to Lemma~\ref{lem:singLk}, we have:

\begin{lem}\label{lem:singLk2}
	Assume $c_{\gamma} = -m c_1$ for $0 \le m \le n-1$ and $c_1 \neq 0$. Then there is an isomorphism of \'etale germs 
    \[
    (\ZCB_{c}(\Gamma_n),p) \cong (\overline{\mc{O}(k,2k+m)} \times T^*( \C^{n - k(k+ m)}),0)
    \]
    for $p \in \LG_k$ and $k \in \Z_{\ge 0}$ such that $k(k+m) \le n$.
\end{lem}

\begin{lem}
	Assume $c_1 = 0$ and $c_\gamma\ne 0$. If $\lambda \in \mc{P}(n)$ and $p \in \LG_{\lambda}$ then we have an isomorphism of \'etale germs 
	$$
	(\ZCB_{c}(\Gamma_n),p) \cong ((\h \times \h^*)/\s_{\lambda} \times T^* \C^{\ell(\lambda)},0), 
	$$ 
 where $\h$ denoted the reflection representation for $\s_{\lambda}$. 
\end{lem}

\begin{proof}
 In this case, there exist pairwise non-isomorphic simple representations $M_{\infty},M_1, \ds, M_{\ell}$ of the deformed preprojective algebra such that $\dim M_{\infty} = e_{\infty}$ and $\dim M_i = \delta$ for $i > 0$ so that the point $p$ corresponds to the semi-simple representation $M_0 \oplus M_1^{\oplus \lambda_1} \oplus \cdots \oplus M_{\ell}^{\oplus \lambda_{\ell}}$. Theorem~\ref{thm:CBLunaslice} says that $(\ZCB_{c}(\Gamma_n),p)$ is isomorphic to the quiver variety $\mathfrak{M}_0(\mathsf{G}(\bc),(1,\lambda_1, \ds, \lambda_{\ell}))$, where $\mathsf{G}(\bc)$ is a graph with one central vertex $\infty$, $\ell$ outer vertices, and an edge between $\infty$ and each outer vertex and a single loop at each outer vertex. This is illustrated as
 
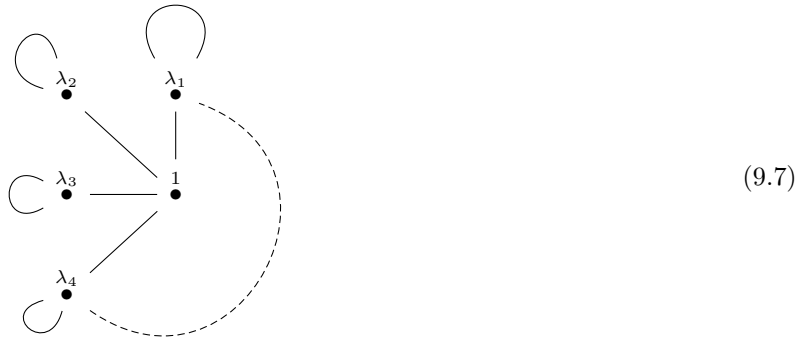
\begin{equation}\label{eq:diagramquiverbreaking}
      \begin{tikzcd}
 	\stackrel{\lambda_2}{\bullet} \arrow[out=165,in=105,loop,looseness=5,no head] & \stackrel{\lambda_1}{\bullet} \arrow[out=120,in=60,loop,looseness=5,no head] & \\
 	\stackrel{\lambda_3}{\bullet} \arrow[out=150,in=210,loop,looseness=5,no head] & \ar[dl,no head]  \ar[l,no head] \ar[lu,no head] \ar[u,no head]    \stackrel{1}{\bullet}  & \\
 	\stackrel{\lambda_4}{\bullet}\arrow[out=195,in=255,loop,looseness=6,no head]  \ar[out=325,in=340,uur,dashed,no head,looseness=2,no head] & & 
 \end{tikzcd}
 \end{equation}
 Let $\mathsf{G}_J$ be the one loop (and one vertex) graph. Since the dimension at the central vertex is $1$ and one can check that the hypothesis of Lemma~\ref{lem:quiverbraking2} hold, we can "break'' the quiver variety at this vertex to get an isomorphism 
 $$
 \mf{M}_{0}(\mathsf{G}(\bc),(1,\lambda_1, \ds, \lambda_{\ell})) \cong \XC_0(\mathsf{G}_J,(\lambda_1)) \times \cdots \times \XC_0(\mathsf{G}_J,(\lambda_{\ell})). 
 $$
 It is well-known \cite[Lemma~2.11]{GanGinzburg} that if $\h$ is the reflection representation for $\s_n$ then 
 $$
 \XC_0(\mathsf{G}_J,(n)) \cong \C^{2n} / \s_n \cong (\h \times \h^*)/ \s_n \times T^* \C.
 $$
The result follows. 
\end{proof}

\begin{remark}
	In the case $c_1 = 1$ and $c_{\gamma} = m c_1$, and $n = k(k+m)$, there is a bijection between the symplectic leaves $\mc{O}(r, 2k + m)$ in $\overline{\mc{O}(k,2k + m)}$ and the leaves $\LG_r$ in $\ZCB_{c}(\Gamma_n)$. This bijection preserves dimension. Moreover, in both cases the order one gets is a total ordering. 
\end{remark}

\subsection{Leaf closures} In this section, we give a proof of Theorem~\ref{thm:mainclosureB}.

\begin{thm}\label{thm:quiverleafclosurempos}
	If $c_1 = 1$ and $c_{\gamma} = m$ with $1 \le m < n$ then 
	$$
	\overline{\LG}_k \cong \ZCB_{c'}(\Gamma_{n - k(k+m)})
	$$
	where $c_1' = 1$ and $c_{\gamma}' = m+2k$. 
\end{thm}

\begin{proof}
In order to apply Theorem~\ref{thm:quiverleafclosure}, we work in the quiver variety $\XC_{{\bc}}(n \delta)$. By \eqref{eq:lambdakappc1}, the parameter ${\bc}$ equals $(m-1,-m)$. As noted in the proof of Lemma~\ref{lem:singLk}, the semisimple representations $M$ belonging to $\LG_k$ are of the form $M = M_0 \oplus M_1^{\oplus k}$, where $\dim M_0 = e_{\infty} + n \delta - k \eta^{(1)}$ and $\dim M_1 = \eta^{(1)}$. Since $m \ge 0$, $\eta^{(1)} = m e_0 + (m-1) e_1$ as in \cite[Lemma~5.4]{MartinoThesis}. Then the decomposition $\tau$ labeling $\LG_k$ is $(n \delta - k \eta^{(1)},(k,\eta^{(1)}))$. Proposition~\ref{prop:leafclosurenormal} says that the leaf closure is normal and this closure can be identified with  
	$$
	\XC_{\bc}(n \delta - k \eta^{(1)}) \times \mf{M}_{{\bc}}(\eta^{(1)}) \cong \XC_{{\bc}}(n \delta - k \eta^{(1)})
	$$
	because $\eta^{(1)}$ is a real root so $\mf{M}_{{\bc}}(\eta^{(1)}) = \{ \mr{pt} \}$. Therefore, we need to identify $\XC_{{\bc}}(n \delta - k \eta^{(1)})$ with a Calogero--Moser variety. Note that $p_{\Lambda_0}(n \delta - k \eta^{(1)}) = n - k(m+k)$ by \eqref{eq:peinftycomputation}. We need to find a sequence of admissible reflections taking $e_{\infty} + n \delta - k \eta^{(1)}$ to $e_{\infty} + (n - k(m+k)) \delta$. It is shown in the proof of \cite[Proposition 5.7(2)]{MartinoThesis} that  
	\begin{equation}\label{eq:reflectionsbvminusbeta}
	e_{\infty} + n \delta - k \eta^{(1)} = \left\{ \begin{array}{ll}
		(s_1 s_0)^k(e_{\infty} + (n - k(m+k)) \delta) & \textrm{if $c_{\gamma} = m c_1$ for $1 \le m < n$,} \\
		s_0(s_1 s_0)^{k-1}(e_{\infty} + (n - k(m+k)) \delta) & \textrm{if $c_{\gamma} = -m c_1$ for $0 \le m < n$.} \\	
	\end{array}\right. 	
	\end{equation}
	In other words,   
	$$
	(s_0 s_1)^k (e_{\infty} + n \delta - k \eta^{(1)}) = e_{\infty} + (n - k(m+k)) \delta =: \alpha'.
	$$ 
	Since $(e_{\infty},e_0) = -1, (e_{\infty}, e_1) = 0$ and $(e_0, e_1) = -2$, we have 
	$$
	s_0(a e_{\infty} + b e_0 + c e_1) = a e_{\infty} + (a - b + 2c) e_0 + c e_1 , \quad s_1(a e_{\infty} + b e_0 + c e_1) = a e_{\infty} + b e_0 + (2b - c) e_1.
	$$
	The dual action on parameters is given by 
 $$
	s_0^* ({\bc}) = (-{\bc}_0,{\bc}_1 + 2 {\bc}_0), \quad s_1^* ({\bc}) = ({\bc}_0 + 2 {\bc}_1,- {\bc}_1).
	$$ 
	Hence $s_0^* s_1^*({\bc}) = (- {\bc}_0-2 {\bc}_1,2{\bc}_0+3 {\bc}_1)$. Applying this to ${\bc} = (m-1,-m)$, an induction shows that 
	$$
	(s_0^* s_1^*)^k(m-1,-m) = (m + 2k-1,-m-2k) =: {\bc}'. 
	$$
	Moreover, a quick induction shows that the sequence of reflections $s_0 s_1 \cdots s_0 s_1$ is admissible for this ${\bc}$. Therefore, $\XC_{{\bc}}(n \delta - k \eta^{(1)}) \cong \XC_{{\bc}'}(\alpha')$. If we write $\bc' = (-c_1' + c_{\gamma}', -c_{\gamma}')$, as in \eqref{eq:lambdakappc1}, then $c_1' = 1$ and $c_{\gamma}' = m + 2k$. We deduce that $\XC_{{\bc}}(n \delta - k \eta^{(1)}) \cong \ZCB_{c'}(\Gamma_{n - k(k+m)})$ where $c' = (1,m+2k)$.
\end{proof}

\begin{remark}
    We could also apply Proposition~\ref{prop:wbccomp} to compute the parameter $\bc'$ in the proof of Theorem~\ref{thm:quiverleafclosurempos}, thus avoiding the computations with $(s_0 s_1)^k$. Since $\bc(\delta) = -1$ in our setting, Proposition~\ref{prop:wbccomp} shows that 
    \[
    \overline{\LG}_k \cong \XC_{\bc}(n \delta - \beta) \cong \XC_{\bc + \overline{\beta}}((n- k(m+k))\delta)),
    \]
    where $\beta = k \eta^{(1)} = km e_0 + k(m-1) e_1$. Then $\overline{\beta} = (2k,-2k)$ and hence $\bc' = (m+2k-1,-m-2k)$.  
\end{remark}

\begin{thm}\label{thm:Bnminusm}
	If $c_1 = 1$ and $c_{\gamma} = -m$ with $0 \le m < n$ then $\overline{\LG}_k \cong \ZCB_{c'}(\Gamma_{n - k(k+m)})$ where $c_1' = 1$ and $c_{\gamma}' = m+2k$. 
\end{thm}

\begin{proof}
	This is identical to the proof of Theorem~\ref{thm:quiverleafclosurempos}, except that now ${\bc} = (-(m+1),m)$ and 
 \[
 n \delta - k \eta^{(1)} = (n-km) e_0 + (n-k(m+1)) e_1.
 \]
This means that $(n - k(m+k)) \delta = (s_0 s_1)^{k-1} s_0 \star (n \delta - k \eta^{(1)})$ and 
	$$
	{\bc}' := (s_0^* s_1^*)^{k-1} s_0^* (-(m+1),m) = (m + 2k -1, - m - 2k).  
	$$
	since $s_0^*({\bc}) = (m+1, -m - 2)$. Once again, induction shows that the sequence $s_0, s_1, \ds, s_0, s_1, s_0$ is admissible for ${\bc}$.  Therefore, $\XC_{{\bc}}(n \delta - k \eta^{(1)}) \cong \XC_{{\bc}'}((n - k(m+k)) \delta)$. If we write $\bc' = (-c_1' + c_{\gamma}', -c_{\gamma}')$, as in \eqref{eq:lambdakappc1}, then $c_1' = 1$ and $c_{\gamma} = m + 2k$ and thus $c_1 = 1$. We deduce that $\XC_{{\bc}}(n \delta - k \eta^{(1)}) \cong \ZCB_{c'}(\Gamma_{n - k(k+m)})$ where $c' = (1,m+2k)$.
\end{proof}

\begin{thm}
	If $c_1 = 0, c_\gamma = 1$ and $\rho \in \mc{P}(n)$ then $\widetilde{\LG}_{\rho} \cong \ZCB_{c'}(\Z_2 \wr \s(\rho))$, where $c'(t) = 1$ for all reflections $t \in \Z_2^{\ell}$ and $c'(s) = 0$ for all reflections (= transpositions) in $\s(\rho)$.
\end{thm}

\begin{proof}
By \eqref{eq:lambdakappc1}, the parameter ${\bc}$ equals $(-1,1)$. As noted in the proof of Lemma~\ref{lem:singLk}, the semisimple representations $M$ belonging to $\LG_{\rho}$ are of the form $M = M_0 \oplus M_1^{\oplus \rho_1} \oplus \cdots \oplus M_{\ell}^{\oplus \rho_{\ell}}$, where $\dim M_0 = e_{\infty}$ and $\dim M_i = \delta$ otherwise. Since $e_{\infty}$ is a real root, $\XC_{{\bc}}(0)$ is a point. Therefore, Theorem~\ref{thm:quiverleafclosure} says that $	\widetilde{\LG}_{\rho}  \cong \mf{M}_{\bc}(\mathsf{G}(\Gamma),\delta)^{\ell} / \s(\rho)$, where $\ell = \ell(\rho)$. Since ${\bc}_{0}, {\bc}_1 \neq 0$, the quiver variety $\mf{M}_{\bc}(\mathsf{G}(\Gamma),\delta)$ is the generic deformation of the Kleinian singularity $\C^2 / \Z_2$. In other words, $\mf{M}_{{\bc}}(\mathsf{G}(\Gamma),\delta) \cong \ZCB_{c_{\gamma}}(\Z_2)$, where $c_{\gamma} = 1$. It follows that $\mf{M}_{{\bc}}(\mathsf{G}(\Gamma),\delta)^{\ell} / \s(\rho) \cong \ZCB_{c'}(\Z_2 \wr \s(\rho))$, where $c'(t) = 1$ for all reflections $t \in \Z_2^{\ell}$ and $c'(s) = 0$ for all reflections (= transpositions) in $\s(\rho)$. 
\end{proof}

\section*{Index of notation}

\begin{multicols}{2}

	$\Phi$ \ finite (irreducible) simply laced root system \hfill\pageref{rootsystemnotation}
    
	$\boldsymbol{\Lambda}$ \ root lattice of $\Phi$ \hfill\pageref{rootsystemnotation}
	
    $\Phi^+$ \ positive roots in $\Phi$ \hfill\pageref{rootsystemnotation}
	
    $\theta$ \ highest root in $\Phi^+$ \hfill\pageref{rootsystemnotation}
	
    $W$ \ Weyl group of $\Phi$ \hfill\pageref{sec:anotherpresentation}
	
    $R$ \ affine simply laced root system \hfill\pageref{rootsystemnotation}
	
    $R^+$ \ positive roots in $R$ \hfill\pageref{rootsystemnotation}
	
    $\Delta = \{ e_0, \ds, e_r \}$ \ set of simple roots in $R^+$ \hfill\pageref{rootsystemnotation}
	
    $\delta$ \ minimal imaginary root in $R$ \hfill\pageref{rootsystemnotation}
	
    $Q$ \ root lattice for $R$ \hfill\pageref{subs:graphs}
	
    $W^{\aff}$ \ affine Weyl group \hfill\pageref{rootsystemnotation}
	
    $F$ \ fundamental region for $\mathsf{G}$ \hfill\pageref{root-sys}
    
    $\mathsf{G}$ \ graph \hfill\pageref{subs:graphs}
    
    $\Qu$ \ quiver with underlying graph $\mathsf{G}$ \hfill\pageref{subs:graphs}
    
    $\mf{M}_{\bc}(\bw,\alpha)$ \ framed quiver variety \hfill\pageref{sec:quivernotation}
    
    $\mf{M}_{\bc}(\alpha)$ \ unframed ($\bw=0$) quiver variety \hfill\pageref{eq:quivernotation1} 
    
    $\XC_{\bc}(\alpha)$ \ framed quiver variety with $\bw=\Lambda_0$ \hfill\pageref{eq:quivernotation1} 
	
    $W_{\mathsf{G}}$ \ Weyl group associated to graph $\mathsf{G}$ \hfill\pageref{subs:graphs}
	
    $\Gamma$ \ finite subgroup of $SL(2,\C)$ \hfill\pageref{sec:symplc-refl-alg}
    
    $\mr{P}$\ parabolic subgroup of $\Gamma_n$  \hfill\pageref{sec:paraboliclabel}
    
    $\s_n$ \ symmetric group on $n$ letters \hfill\pageref{sec:symplc-refl-alg}
    
    $\Qu^\ell$ \ cyclic quiver of length $\ell$ (infinity linear quiver for $\ell=\infty)$ \hfill\pageref{sec:combinatoricsCM}
    
    $\Qov^\ell$ \ double quiver for $\Qu^\ell$ \hfill\pageref{sub:quiver-cycl}
    
    $\Qulinf$\ framed version of $\Qu^\ell$ \hfill\pageref{sub:quiver-cycl}
    
    $\oQulinf$ \ double quiver for $\Qulinf$ \hfill\pageref{sub:quiver-cycl}
    
    $\PC$ \ the set of partitions \hfill\pageref{sec:notation}
    
    $\CC_\ell$ \ the set of $\ell$-cores \hfill\pageref{sec:lcores}
    
    $\CC_J$ \ the set of $J$-cores \hfill\pageref{sec:Jcore}
    
    $\Sigma_\bc$ \ the set of dimension vectors of simple representations for the deformed preprojective algebra  \hfill\pageref{root-sys}
    
    $\Sigma\Sigma_\bc$ \ all possible sums of elements of $\Sigma_\bc$ \hfill\pageref{sec:sympleaf1}
     
    $E_\bc$ \ framed version of $\Sigma_\bc$ \hfill\pageref{sec:CBtrick}
    
    $\Hb_c$ \ symplectic reflection algebra  \hfill\pageref{sec:symplc-refl-alg}
    
    $\ZCB_c$ \ Calogero--Moser variety \hfill\pageref{sec:symplc-refl-alg}
    
    $\LG$ \ symplectic leaf \hfill\pageref{sec:sympleavesdefn}

\end{multicols}

\bibliographystyle{plain}
\bibliography{biblo}

\end{document}